\documentclass[11pt]{amsart}
\usepackage[T1]{fontenc}
\usepackage[utf8]{inputenc}
\usepackage{lmodern}
\usepackage[margin=1in]{geometry}
\usepackage{amsmath,amssymb,amsthm,mathtools}
\usepackage{tikz-cd}
\usepackage{xcolor}
\colorlet{auditred}{red}
\usepackage{microtype}
\usepackage{needspace}
\usepackage[numbers,square]{natbib}
\usepackage[hidelinks]{hyperref}
\newcommand{\doi}[1]{\url{https://doi.org/#1}}
\usepackage{enumitem}
\usepackage{longtable,booktabs,array}

\newcommand{\Q}{\mathbb Q}
\newcommand{\C}{\mathbb C}
\newcommand{\Z}{\mathbb Z}
\newcommand{\F}{\mathbb F}
\newcommand{\Ql}{\mathbb Q_\ell}
\newcommand{\Zl}{\mathbb Z_\ell}

\newcommand{\cM}{\mathcal M}
\newcommand{\M}{\textcolor{black}{\mathcal M}}

\newcommand{\cC}{\mathcal C}

\newcommand{\HH}{\mathbb{H}}

\newcommand{\et}{\mathrm{\acute et}}

\newcommand{\geom}{\mathrm{geom}}
\DeclareMathOperator{\End}{\operatorname{End}}

\DeclareMathOperator{\Sp}{Sp}
\DeclareMathOperator{\GSp}{GSp}

\DeclareMathOperator{\Lie}{Lie}
\DeclareMathOperator{\Spec}{Spec}

\DeclareMathOperator{\Conf}{Conf}
\DeclareMathOperator{\Jac}{Jac}

\newtheorem{theorem}{Theorem}[section]
\newtheorem*{theorem*}{Theorem}
\newtheorem{proposition}[theorem]{Proposition}
\newtheorem{lemma}[theorem]{Lemma}
\newtheorem{corollary}[theorem]{Corollary}
\theoremstyle{definition}

\theoremstyle{remark}

\definecolor{revisionred}{RGB}{0,0,0}

\title[Monodromy of divisors]{Monodromy of divisors in moduli spaces
of curves in positive characteristic}
\author{Ma Luo}
\address{School of Mathematical Sciences, Key Laboratory of MEA
(Ministry of Education), Shanghai Key Laboratory of PMMP,
East China Normal University, Shanghai, China}
\email{mluo@math.ecnu.edu.cn}
\author{Tatsunari Watanabe}
\address{Mathematics Department, Embry-Riddle Aeronautical University,
Prescott, Arizona, USA}
\email{watanabt@erau.edu}
\thanks{Corresponding author: Tatsunari Watanabe
(\texttt{watanabt@erau.edu}).}
\thanks{The first author is supported in part by Shanghai Pilot Program for Basic Research and by Science and Technology Commission of Shanghai Municipality (No. 22DZ2229014).}
\date{}

\begin{document}

\begin{abstract}
Let $D$ be an irreducible divisor in the moduli stack of smooth curves of genus at least three over an algebraically closed field of characteristic $p>0$. We prove that, after pullback to any smooth integral variety dominating $D$, the geometric monodromy of the universal curve on first $\ell$-adic cohomology is Zariski dense in the symplectic group for every prime $\ell\ne p$. The result holds with any number of marked points and extends to fine moduli spaces with prime-to-$p$ level structure and their connected finite \'etale covers.
Over the complex numbers, Hain proved that the closure of the monodromy image has finite index in the profinite symplectic group. Our theorem extends the resulting Zariski density to positive characteristic. For every divisor in the unpointed moduli space, we prove that the geometric generic Jacobian is simple. After a finite extension of the function field of the divisor, we construct a degeneration in which the Jacobian has semistable reduction of toric rank one. This degeneration supplies a rank-one inertia operator and, together with simplicity, implies that every geometric endomorphism of the generic Jacobian is multiplication by an integer.
\end{abstract}

\maketitle
\setcounter{tocdepth}{1}
\tableofcontents

\section{Introduction}
\label{sec:introduction}

How much of the monodromy of the universal curve is retained when
one restricts the family to a divisor in the moduli space of curves?
Over $\C$, Hain proved that if $g\geq3$ and a smooth irreducible
variety dominates a divisor in the moduli space of $n$-pointed
genus-$g$ curves with level structure, then the closure of its
monodromy image in \textcolor{black}{the profinite topology of}
$\Sp_{2g}(\widehat{\Z})$ has finite index
\cite[Thm.~1]{HainDensity}. Richard Hain suggested to the second
author the problem of studying this monodromy in positive
characteristic.

In \emph{Rational points of universal curves}, Hain used his
density theorem to study the universal curve over the function field
$K$ of a level moduli space in characteristic zero. For the
geometrically nontrivial irreducible coefficients of negative weight
in his argument, the residue of a degree-one class along an omitted
divisor is fixed by its geometric monodromy and hence vanishes
\cite[proof of Thm.~6.1]{HainDensity}. For $g\geq5$, if the
$\ell$-adic cyclotomic character of the base field has infinite
image, as it does for a number field, restriction from the moduli
space to its generic point induces an isomorphism on $H^1$ and an
injection on $H^2$ with these coefficients
\cite[Prop.~16.3(ii),(iv)]{HainRational}. Hain used this comparison
to compute the low-degree cohomology of the Lie algebra of the
weighted completion of the absolute Galois group of $K$
\cite[Prop.~17.1]{HainRational} and applied the computation to
sections of the generic curve \cite[\S18]{HainRational}.

Let $k$ be an algebraically closed field of characteristic $p>0$,
and let $g\geq3$. Write $\cM_{g,n/k}$ for the moduli stack of
smooth proper curves of genus $g$ with $n$ distinct ordered
markings, and put $\cM_{g/k}=\cM_{g,0/k}$. Our first result
concerns \textcolor{black}{the geometric generic Jacobian of a divisor.}

\medskip
\noindent\textbf{Theorem~\ref{thm:divisor-geometric-simplicity}
(Geometric simplicity).}
\emph{For every reduced irreducible divisor $D\subset\cM_{g/k}$,
\textcolor{black}{the geometric generic Jacobian of $D$ is simple.}}
\medskip

\textcolor{black}{If $C_{\bar\eta}$ is the curve represented by a
geometric generic point $\bar\eta$ of $D$, this means that
$\Jac(C_{\bar\eta})$ has no nonzero proper abelian subvariety.
In odd characteristic the assertion fails in genus two.} The genus-two curves
admitting a double cover of an elliptic curve form a surface
in $\cM_{2/k}$ \cite[\S\S1 and~4]{ShaskaVoelklein}.
Since $\dim\cM_{2/k}=3$, this locus is a divisor. For a double
cover $f:C\to A$ at a geometric generic point of this divisor,
\textcolor{black}{the composite $\operatorname{Nm}_f\circ f^*$ is
multiplication by two on $\Jac(A)$.} Here
$\operatorname{Nm}_f:\Jac(C)\to\Jac(A)$ is the norm map, given on
degree-zero divisor classes by $[Z]\mapsto[f_*Z]$.
It follows that $\ker(f^*)$ is
contained in the finite group $\Jac(A)[2]$. Since $\Jac(A)$ is an
elliptic curve and $f^*$ has finite kernel, its image is an elliptic subvariety
of $\Jac(C)$.

Our second result determines
the monodromy of the family over any smooth variety dominating a
divisor. Fix a prime $\ell\neq p$. For a geometric fiber
$C_{\bar e}$ of such a family, put
$H_{\Ql}=H^1_{\et}(C_{\bar e},\Ql(1))$, with its alternating
pairing.

\medskip
\noindent\textbf{Theorem~\ref{thm:divisor-monodromy}
(Symplectic monodromy).}
\emph{Let $n\geq0$, let $D\subset\cM_{g,n/k}$ be a reduced
irreducible divisor, and let $E\to D$ be a dominant morphism
from a smooth integral $k$-variety. For a geometric point
$\bar e$ of $E$, \textcolor{black}{let $\rho_{E,\ell}$ denote the action
of the \'etale fundamental group $\pi_1(E,\bar e)$ on $H_{\Ql}$.
Then}}
\[
 \overline{\rho_{E,\ell}\bigl(\pi_1(E,\bar e)\bigr)}^{\,\mathrm{Zar}}
 =\Sp(H_{\Ql}).
\]
\emph{The same conclusion holds for a divisor in the full
symplectic level cover $\cM_{g,n/k}[m]$, where $m\geq3$ is prime
to $p$, and in any connected finite \'etale cover of this level
scheme.}
\medskip

In positive characteristic, Achter and Pries determined the
integral $\ell$-adic monodromy of every irreducible component of
a $p$-rank stratum \cite[Thm.~4.5]{AchterPries}. Their work also
supplies a curve with absolutely simple Jacobian on each component
\cite[Application~5.7]{AchterPries}. We use this in
Section~\ref{sec:geometric-simplicity} for divisors contained in
the nonordinary locus. By the dimension formula of Faber and
van der Geer \cite[Thm.~2.3]{FaberVanDerGeer}, each such divisor
is the closure of a component of the $p$-rank $g-1$ stratum.
For a divisor meeting the ordinary locus, the argument is
different. If its geometric generic Jacobian were not simple,
a Rosati-symmetric endomorphism would give at least $g-1$
linearly independent quadratic differentials annihilating the
tangent space to the divisor at a general point. Since this
tangent space is a hyperplane, its annihilator is only
one-dimensional.

To pass from simplicity to monodromy, we show in
Section~\ref{sec:rank-one-degenerations} that the generic
Jacobian of every unpointed divisor admits, after a finite
extension, semistable reduction of toric rank one. This
degeneration has two consequences. Together with simplicity,
it gives $\End(\Jac(C_{\bar\eta}))=\Z$ by Zarhin's criterion
\cite[Cor.~3.2]{ZarhinBigMonodromy}. It also gives an element
of inertia whose action on the Tate module differs from the
identity by an operator of rank one. In
Section~\ref{sec:monodromy-theorem} we use the endomorphism
calculation to obtain an absolutely irreducible arithmetic
action. The conjugates of the rank-one operator under the
connected arithmetic monodromy group then generate the full
geometric symplectic Lie algebra.

\textcolor{black}{Related uses of boundary geometry and
semiabelian degeneration appear in Urbanik and Yang's study
of splitting loci of Jacobians
\cite[\S\S5.2--5.4]{UrbanikYangSpecialLoci}. Our argument
combines degeneration of toric rank one with geometric
simplicity to establish full symplectic geometric monodromy
for families dominating arbitrary divisors in $\cM_{g/k}$,
for $g\geq3$. The result holds in every characteristic $p>0$
and for every prime $\ell\neq p$.}

\section{Moduli spaces and monodromy}
\label{sec:moduli-monodromy}

We recall the standard facts about moduli spaces and their
symplectic local systems used below. We also give the reductions
to fine level schemes, normal bases, and unpointed moduli.
The details are included to specify the hypotheses needed
for the monodromy arguments in positive characteristic.
\subsection{Moduli stacks and universal curves}
\label{subsec:moduli-universal-curves}

Let $k$ be an algebraically closed field of characteristic
$p>0$, and let $g\geq3$ and $n\geq0$.
Write $\cM_{g,n/k}$ for the moduli stack of smooth proper
geometrically connected curves of genus $g$ with $n$
distinct ordered markings. An object over a $k$-scheme $T$
is a smooth proper morphism $f:C\to T$ whose geometric
fibers are connected curves of genus $g$, together with
pairwise disjoint sections $s_1,\ldots,s_n:T\to C$.
We denote the coarse
moduli space of $\cM_{g,n/k}$ by $M_{g,n/k}$.

Denote the universal complete curve by
\[
 \pi_{g,n}:\cC_{g,n/k}\longrightarrow\cM_{g,n/k}.
\]
This morphism is representable, smooth and proper, with
geometrically connected fibers of genus $g$. For $n\geq1$,
it has the tautological sections $s_1,\ldots,s_n$.
Its pullback along the classifying morphism $T\to\cM_{g,n/k}$
of a pointed family is that family, with its specified
sections. Allowing an additional point disjoint from the
markings gives the canonical identification
\[
 \cM_{g,n+1/k}
 \cong\cC_{g,n/k}\setminus
       \bigcup_{j=1}^n s_j(\cM_{g,n/k}).
\]
When $n=0$, we omit the index $0$ from the notation
for the stack, its coarse moduli space, and the universal curve.
In this case the union is empty, so
$\cC_{g/k}\cong\cM_{g,1/k}$.

The forgetful morphism
\[
 r_n:\cM_{g,n/k}\longrightarrow\cM_{g/k}
\]
is representable, separated and smooth of relative dimension $n$.
Indeed, its pullback along a family $C\to T$ is the relative
configuration space obtained from the $n$-fold fiber product
$C^n_T$ by removing the pairwise diagonals.
Its geometric fiber at a curve $C/\Omega$ is
\[
 \Conf_n(C)=C^n\setminus\bigcup_{i<j}\Delta_{ij},
\]
which is smooth and irreducible of dimension $n$;
we use $\Conf_0(C)=\Spec\Omega$.
The stack $\cM_{g/k}$ is a smooth irreducible separated
Deligne--Mumford stack of dimension $3g-3$
\cite[Secs.~1 and~3, Prop.~5.1 and Thm.~5.2]{DeligneMumford}.
The description of $r_n$ therefore gives the same
properties for $\cM_{g,n/k}$, with
$\dim\cM_{g,n/k}=3g-3+n$.

Fix an integer $m\geq3$ prime to $p$ and a primitive
$m$th root of unity $\zeta_m\in k$.
For a family $C\to T$, let $J=\Jac(C/T)$, with its
canonical principal polarization $\lambda$, and let
$e_m^\lambda:J[m]\times_T J[m]\to\mu_{m,T}$ be the
associated Weil pairing.
A full symplectic level-$m$ structure is an isomorphism
of finite \'etale group schemes
\[
 \alpha:(\Z/m\Z)^{2g}_T\xrightarrow{\sim}J[m]
\]
such that
\[
 e_m^\lambda(\alpha(v),\alpha(w))
   =\zeta_m^{\langle v,w\rangle_0},
\]
where $\langle\ ,\ \rangle_0$ is the standard alternating
form on $(\Z/m\Z)^{2g}$ and the identity is understood
for local sections $v,w$.
Write $\cM_{g,n/k}[m]$ for the resulting level moduli stack.
The choice of $\zeta_m$ is fixed throughout.

Forgetting the level structure gives a representable finite
\'etale surjection
\[
 q_m:\cM_{g,n/k}[m]\longrightarrow\cM_{g,n/k}.
\]
In the \'etale topology on the target stack, this is a torsor
under the constant group $\Sp_{2g}(\Z/m\Z)$: symplectic
bases exist \'etale locally, and any two differ by a unique
symplectic automorphism.
For $m\geq3$, the level moduli stack is represented by a
scheme \cite[Sec.~5.14, p.~108]{DeligneMumford}.
The unpointed level scheme is connected by
\cite[Thm.~5.15, pp.~108--109]{DeligneMumford}.
Since a level structure depends only on the underlying curve,
\[
 \cM_{g,n/k}[m]
 \cong\cM_{g,n/k}\times_{\cM_{g/k}}\cM_{g/k}[m].
\]
The configuration-space description above shows that
$\cM_{g,n/k}[m]$ is also connected. It is therefore a
smooth integral $k$-scheme of dimension $3g-3+n$.

The universal complete curve at level $m$ is the pullback
\[
\begin{tikzcd}
 \cC_{g,n/k}[m] \arrow[r] \arrow[d,"{\pi_{g,n}[m]}"'] &
 \cC_{g,n/k} \arrow[d,"{\pi_{g,n}}"] \\
 \cM_{g,n/k}[m] \arrow[r,"q_m"'] & \cM_{g,n/k}.
\end{tikzcd}
\]
We use the same notation $s_j$ for the pulled-back
tautological sections. All universal families below are
taken over the moduli stacks or their fine level covers.
Unless otherwise stated, a divisor is a reduced irreducible
closed substack, or subscheme, of codimension one.
We choose compatible geometric basepoints whenever
fundamental groups are compared, and write $\pi_1$ for
the \'etale fundamental group.
\subsection{The symplectic local system}
\label{subsec:symplectic-local-system}

Fix a prime $\ell\neq p$. Let $\M$ denote either
$\cM_{g,n/k}$ or its level cover $\cM_{g,n/k}[m]$, and
write $\pi:\cC\to\M$ for the universal complete curve.
Put
\[
 \HH_{\Zl}=R^1\pi_*\Zl(1),\qquad
 \HH_{\Ql}=\HH_{\Zl}\otimes_{\Zl}\Ql.
\]
Here $\HH_{\Zl}$ denotes the inverse system
$\bigl(R^1\pi_*\mu_{\ell^r}\bigr)_{r\geq1}$.
Proper smooth base change shows that it is a lisse
$\Zl$-sheaf of rank $2g$, whose fiber at a geometric
point $\bar x\to\M$ is
\[
 H_{\Zl}=H^1_{\et}(C_{\bar x},\Zl(1));\qquad
 H_{\Ql}=H_{\Zl}\otimes_{\Zl}\Ql.
\]
Here $C_{\bar x}$ is the fiber of the universal
curve at $\bar x$. We use proper base change for the fiber identification
\cite[Prop.~14.2, Thms.~17.7 and~20.2]{MilneEtale}.
When $\M$ is a stack, these constructions are made on
\'etale charts and descend along their overlaps.

Cup product and the relative trace give a perfect pairing
\begin{equation}
 \label{eq:symplectic-pairing}
 \HH_{\Zl}\otimes_{\Zl}\HH_{\Zl}
 \longrightarrow R^2\pi_*\Zl(2)
 \xrightarrow{\ \sim\ }\Zl(1).
\end{equation}
Perfection follows from Poincar\'e duality for curves
\cite[Thm.~14.7 and Ex.~14.8]{MilneEtale}.
This pairing is alternating, including when $\ell=2$:
on each geometric fiber, graded commutativity gives
$2(v\cup v)=0$, and the target
$H^2_{\et}(C_{\bar x},\Zl(2))\cong\Zl(1)$
is torsion-free.
If $J\to\M$ is the relative Jacobian, the Kummer sequence
gives compatible identifications
\[
 R^1\pi_*\mu_{\ell^r}\cong J[\ell^r],\qquad
 \HH_{\Zl}\cong T_\ell J.
\]
Under these identifications, \eqref{eq:symplectic-pairing}
is the inverse limit of the Weil pairings associated to
the canonical principal polarization of $J$
\cite[Ex.~7.9(a), Props.~14.1--14.2 and Ex.~14.8]{MilneEtale}.

Since $k$ is algebraically closed, a choice of compatible
primitive $\ell$-power roots of unity identifies the constant
sheaf $\Zl(1)$ with $\Zl$. When $m=\ell^r$, choose the system
so that its primitive $\ell^r$th root is the root $\zeta_m$ fixed
in the definition of the level structure. Fix such a system. We write
$\Sp(H_{\Zl})$ for the group of automorphisms of
$H_{\Zl}$ preserving the resulting alternating form, and
$\Sp(H_{\Ql})$ for the corresponding algebraic group
over $\Ql$. The local system defines a continuous
representation
\[
 \rho_{\M,\ell}:\pi_1(\M,\bar x)
 \longrightarrow\Sp(H_{\Zl}).
\]

More generally, let $f:S\to\M$ be a morphism from a
connected locally noetherian $k$-scheme, and choose
compatible geometric basepoints. When the map to
$\cM$ is understood, we also write $\HH_{\Zl}$ and $\HH_{\Ql}$ for
their pullbacks. For $\bar x=f(\bar s)$, we identify the stalk of
each pullback at $\bar s$ with the stalk of the original sheaf
at $\bar x$. Pullback gives
\[
 \rho_{S,\ell}:\pi_1(S,\bar s)
 \xrightarrow{\ f_*\ }\pi_1(\M,\bar x)
 \xrightarrow{\ \rho_{\M,\ell}\ }\Sp(H_{\Zl}).
\]
We use the same symbol $\rho_{S,\ell}$ for the induced action on
$H_{\Ql}=H_{\Zl}\otimes_{\Zl}\Ql$. Set
\[
 \Gamma_{S,\ell}=\rho_{S,\ell}(\pi_1(S,\bar s)),\qquad
 G_{S,\ell}=\overline{\Gamma_{S,\ell}}^{\,\mathrm{Zar}}
 \subseteq\Sp(H_{\Ql}).
\]
The integral image $\Gamma_{S,\ell}$ is compact and hence
closed in the $\ell$-adic topology. We call $G_{S,\ell}$
the algebraic monodromy group, and denote its identity
component by $G_{S,\ell}^{\circ}$. We say that the family
over $S$ has full symplectic monodromy at $\ell$ if
$G_{S,\ell}=\Sp(H_{\Ql})$.
These groups are well defined up to conjugacy under
change of compatible geometric basepoints. We use the
same notation when $S$ is one of the moduli stacks.

The following standard monodromy statement follows from
the connectedness of full symplectic level covers
\cite[Thm.~5.15]{DeligneMumford}. We give the argument
to include the integral assertion for every $\ell\neq p$.

\begin{lemma}
\label{lem:ambient-monodromy}
The image of $\rho_{\cM_{g,n/k},\ell}$ is
$\Sp(H_{\Zl})$. The image of
$\rho_{\cM_{g,n/k}[m],\ell}$ is an open subgroup of
$\Sp(H_{\Zl})$. In particular, both bases have
full symplectic monodromy at $\ell$.
\end{lemma}

\begin{proof}
For every $r$ with $\ell^r\geq3$, the torsor of symplectic
bases of $\HH_{\Zl}/\ell^r\HH_{\Zl}$ is the full
symplectic level-$\ell^r$ cover. It is connected by
Subsection~\ref{subsec:moduli-universal-curves}, using
\cite[Thm.~5.15]{DeligneMumford}. A connected torsor
under a finite constant group has surjective monodromy.
Thus $\rho_{\cM_{g,n/k},\ell}$ is surjective modulo
$\ell^r$ for every such $r$. These $r$ form a cofinal
set, so its image is dense in $\Sp(H_{\Zl})$.
The image is also compact, and therefore equals
$\Sp(H_{\Zl})$.

The connected finite \'etale cover
$q_m:\cM_{g,n/k}[m]\to\cM_{g,n/k}$ identifies the
fundamental group of the source with an open subgroup
of that of the target. Its monodromy image consequently
has finite index in $\Sp(H_{\Zl})$. This image is
compact, hence closed, and therefore open.
Finally, an open subgroup contains a principal congruence
subgroup and is Zariski dense in $\Sp(H_{\Ql})$:
after choosing a symplectic basis, its intersection with
each root subgroup contains
$\ell^a\Zl$ for some $a$, and the root subgroups
generate the symplectic group.
\end{proof}

The Tate twist in the definition of $\HH_{\Zl}$ is
useful when the family is descended to a non-algebraically
closed field. In that setting the same pairing is
$\Zl(1)$-valued, so the arithmetic monodromy lies in
the group of integral symplectic similitudes
$\GSp(H_{\Zl})$ with multiplier the $\ell$-adic
cyclotomic character. Its restriction to the geometric
fundamental group lies in $\Sp(H_{\Zl})$.
\subsection{Dominant morphisms and normalization}
\label{subsec:dominant-morphisms}

We recall a standard open-image property of dominant
morphisms to normal varieties. We include its proof to
make clear that no separability hypothesis is needed,
and then apply it to the normalizations of divisors.
In this subsection, a variety is an integral
separated scheme of finite type over $k$. These reductions
will be applied to divisors in the fine level schemes
of Subsection~\ref{subsec:moduli-universal-curves}.

\begin{lemma}
\label{lem:dominant-open-image}
Let $f:X\to Y$ be a dominant morphism of varieties,
with $Y$ normal. Then
\[
 f_*:\pi_1(X,\bar x)\longrightarrow\pi_1(Y,\bar y)
\]
has open image. If $U\subseteq Y$ is a nonempty open
subscheme, then $\pi_1(U)\to\pi_1(Y)$ is surjective.
\end{lemma}

\begin{proof}
Put $K=k(Y)$ and $L=k(X)$, and write
$G_F=\operatorname{Gal}(F^{\mathrm{sep}}/F)$ for the
absolute Galois group of a field $F$.
Choose compatible separable closures
$K^{\mathrm{sep}}\subseteq L^{\mathrm{sep}}$, and set
\[
 K'=L\cap K^{\mathrm{sep}}.
\]
The extension $L/K$ is finitely generated, so the
relative algebraic closure of $K$ in $L$ is finite
over $K$ \cite[Tag~037J]{Stacks}. In particular,
$K'/K$ is finite separable.
The image of the restriction map $G_L\to G_K$ is
closed, since $G_L$ is compact, and its fixed field
in $K^{\mathrm{sep}}$ is $K'$. Galois correspondence
therefore identifies this image with
$\operatorname{Gal}(K^{\mathrm{sep}}/K')$, an open
subgroup of $G_K$.

Normality of $Y$ gives a surjection
$G_K\to\pi_1(Y)$ \cite[Tag~0BQM]{Stacks}.
The commutative diagram
\[
\begin{tikzcd}
 G_L \arrow[r] \arrow[d] & G_K \arrow[d,two heads] \\
 \pi_1(X) \arrow[r,"f_*"'] & \pi_1(Y)
\end{tikzcd}
\]
shows that the image of $f_*$ contains a subgroup
of index at most $[K':K]$ in $\pi_1(Y)$.
Thus this image has finite index. It is also compact,
hence closed, and is therefore open.

For the final assertion, the generic point of $Y$
belongs to $U$, so the surjection
$G_K\to\pi_1(Y)$ factors through $\pi_1(U)$.
\end{proof}

This proof applies to inseparable dominant morphisms
as well. In particular, no separability hypothesis on
$k(X)/k(Y)$ is needed in the following consequence.

\begin{corollary}
\label{cor:dominant-monodromy}
Let $X\xrightarrow{f}Y\to\M$ be morphisms with
$f$ as in Lemma~\ref{lem:dominant-open-image}.
For the pullbacks of $\HH_{\Zl}$, the subgroup
$\Gamma_{X,\ell}$ has finite index in
$\Gamma_{Y,\ell}$, and
\[
 G_{X,\ell}^{\circ}=G_{Y,\ell}^{\circ}.
\]
Consequently, the families over $X$ and $Y$ have full
symplectic monodromy at $\ell$ simultaneously.
Moreover, $\Gamma_{X,\ell}$ is open in
$\Sp(H_{\Zl})$ if and only if $\Gamma_{Y,\ell}$ is.
For a nonempty open $U\subseteq Y$, one has the
stronger equality
$\Gamma_{U,\ell}=\Gamma_{Y,\ell}$.
\end{corollary}

\begin{proof}
The finite-index assertion follows by applying
$\rho_{Y,\ell}$ to the image in
Lemma~\ref{lem:dominant-open-image}. A finite coset
decomposition of $\Gamma_{Y,\ell}$ over
$\Gamma_{X,\ell}$ gives, on taking Zariski closures,
a finite union of cosets of $G_{X,\ell}$ covering
$G_{Y,\ell}$. Their identity components are therefore
equal. Since $\Sp(H_{\Ql})$ is connected, this proves
the assertion about full symplectic monodromy.

Both integral images are compact. Thus
$\Gamma_{X,\ell}$ is a closed finite-index subgroup
of $\Gamma_{Y,\ell}$, hence open in it. This proves
the equivalence of openness in $\Sp(H_{\Zl})$.
The last assertion follows from the surjectivity of
$\pi_1(U)\to\pi_1(Y)$.
\end{proof}

We now apply these observations to a possibly singular
divisor $D\subset\cM_{g,n/k}[m]$. Let
\[
 \nu:D^\nu\longrightarrow D
\]
be its normalization, and let $D^{\mathrm{sm}}$ denote
\textcolor{black}{the smooth locus of $D$}. Since $k$ is perfect,
$D^{\mathrm{sm}}$ is nonempty. The morphism $\nu$ is
finite and birational, and $D^\nu$ is a normal variety
\cite[Tag~0BXR]{Stacks}. Over the normal open
$D^{\mathrm{sm}}$, it is an isomorphism
\cite[Tag~0AB1]{Stacks}. We therefore regard
$D^{\mathrm{sm}}$ as a dense open subscheme of $D^\nu$.

\begin{corollary}
\label{cor:normalization-monodromy}
Let $E$ be a smooth integral $k$-variety and
$h:E\to D$ a dominant morphism. Then $h$ factors
uniquely as
\[
 E\xrightarrow{\widetilde h}D^\nu
 \xrightarrow{\nu}D,
\]
where $\widetilde h$ is dominant. With compatible
geometric basepoints, one has
\[
 \Gamma_{D^{\mathrm{sm}},\ell}=\Gamma_{D^\nu,\ell},
 \qquad
 [\Gamma_{D^\nu,\ell}:\Gamma_{E,\ell}]<\infty,
\]
and hence
\[
 G_{E,\ell}^{\circ}=G_{D^\nu,\ell}^{\circ}
   =G_{D^{\mathrm{sm}},\ell}^{\circ}.
\]
Full symplectic monodromy at $\ell$, as well as
openness of the integral monodromy image, is equivalent
for these three bases.
\end{corollary}

\begin{proof}
The variety $E$ is normal, so the factorization exists
and is unique by the universal property of normalization
\cite[Tag~035Q, part~(4)]{Stacks}.
It is dominant because $k(D^\nu)=k(D)$ and $h$ is
dominant. Apply Corollary~\ref{cor:dominant-monodromy}
to $E\to D^\nu$ and to the dense open immersion
$D^{\mathrm{sm}}\hookrightarrow D^\nu$.
\end{proof}

Thus, to prove full symplectic monodromy for every smooth
integral variety dominating a divisor, it suffices to
prove it over the normalization of that divisor, or
equivalently over its smooth locus.
\subsection{Finite covers and extension of constants}
\label{subsec:finite-covers-constants}

We first describe the effect of a finite \'etale cover on
the monodromy of $\HH_{\Zl}$. The following standard
consequence of the Galois correspondence will allow us
to pass from moduli stacks to fine level schemes.

\begin{lemma}
\label{lem:finite-etale-monodromy}
Let $S\to\M$ be a morphism from a connected scheme
of finite type over $k$, and let $v:S'\to S$ be a
connected finite \'etale cover of degree $d$. Then
\[
 [\Gamma_{S,\ell}:\Gamma_{S',\ell}]\leq d,
 \qquad
 G_{S',\ell}^{\circ}=G_{S,\ell}^{\circ}.
\]
Full symplectic monodromy at $\ell$ and openness of
the integral monodromy image are each equivalent for
$S'$ and $S$.
\end{lemma}

\begin{proof}
With compatible geometric basepoints, $v_*$ identifies
$\pi_1(S')$ with an open subgroup of index $d$ in
$\pi_1(S)$, by the Galois correspondence for finite
\'etale covers \cite[Tag~0BND]{Stacks}.
Taking images gives the index bound.
The finite-coset argument in the proof of
Corollary~\ref{cor:dominant-monodromy} gives equality
of the identity components. Connectedness of
$\Sp(H_{\Ql})$ then gives the assertion about full
symplectic monodromy. Finally, $\Gamma_{S',\ell}$
is compact, so it is a closed finite-index subgroup
of $\Gamma_{S,\ell}$ and hence open in it. This
proves the assertion about integral openness.
\end{proof}

\begin{corollary}
\label{cor:reduction-fine-level}
Let $D\subset\cM_{g,n/k}$ be a divisor and let
$E\to D$ be a dominant morphism from a smooth
integral $k$-variety. There exist a connected finite
\'etale cover $E'\to E$, a divisor
$D'\subset\cM_{g,n/k}[m]$, and a dominant morphism
$E'\to D'$ compatible with the maps to $\cM_{g,n/k}$.
Consequently, the assertions of full symplectic monodromy
and openness of the integral image for all such $E$
follow from their analogues for divisors in the fine
level scheme.
\end{corollary}

\begin{proof}
Choose a connected component $E'$ of
\[
 E\times_{\cM_{g,n/k}}\cM_{g,n/k}[m].
\]
Since $q_m$ is representable, finite \'etale and
surjective, this fiber product is a finite \'etale
scheme over $E$. The morphism $E'\to E$ is finite
\'etale and surjective: its image is nonempty, open
and closed in the connected scheme $E$.
The scheme $E'$ is smooth and connected, and hence
integral.

Let $D'$ be the reduced closure of the image of $E'$
in $\cM_{g,n/k}[m]$. Its generic point lies over
the generic point of $D$. Since $q_m^{-1}(D)\to D$
is finite \'etale, $D'$ is an irreducible component
of $q_m^{-1}(D)$. These assertions can be checked on
\'etale charts of $D$. Codimension is preserved by
\'etale morphisms, so $D'$ is a divisor.
The morphism $E'\to\cM_{g,n/k}[m]$ factors through
$D'$ because $E'$ is reduced, and $E'\to D'$ is
dominant by construction. Apply
Lemma~\ref{lem:finite-etale-monodromy} to $E'\to E$.
\end{proof}

By Corollary~\ref{cor:normalization-monodromy}, the
monodromy assertions on the fine level scheme may
in turn be checked on the normalizations of its divisors.

We next recall the standard base-change observation used
to compare geometric monodromy over different fields of
constants. We include the argument to specify the
surjectivity needed for equality of the monodromy images.

\begin{lemma}
\label{lem:extension-algebraically-closed-constants}
Let $k\subseteq K$ be an extension of algebraically
closed fields, and let $S\to\M$ be a morphism from
a connected scheme of finite type over $k$. Then
the natural homomorphism
\[
 \pi_1(S_K)\longrightarrow\pi_1(S)
\]
is surjective. Identifying the fibers of the local
systems by proper base change, one has
\[
 \Gamma_{S_K,\ell}=\Gamma_{S,\ell},\qquad
 G_{S_K,\ell}=G_{S,\ell}.
\]
\end{lemma}

\begin{proof}
Every connected finite \'etale cover $T\to S$ is
geometrically connected over the algebraically closed
field $k$. Thus $T_K$ is connected
\cite[Tag~0363]{Stacks}. The connected-cover criterion
for surjectivity, obtained from the Galois correspondence
\cite[Tag~0BND]{Stacks}, gives the first assertion.
The monodromy representation on $S_K$ is the composite
of this surjection with $\rho_{S,\ell}$, so its image
is exactly $\Gamma_{S,\ell}$.
\end{proof}

For a nonproper base in characteristic $p$, the
fundamental-group homomorphism in this lemma need
not be an isomorphism \cite[Tag~0A47]{Stacks}.
Surjectivity gives the equality of the monodromy
images needed here.

We will also use models over finitely generated fields.
After passing to a fine level scheme, any finite
collection of schemes and morphisms of finite
presentation under consideration, including the curve
families, divisors and chosen finite \'etale covers, descends to a common
subfield $k_0\subset k$ finitely generated over
$\F_p$. We may require $\zeta_m\in k_0$.
This follows from the limit result for objects of
finite presentation and descent of the relevant
morphism properties under faithfully flat base change
\cite[Tags~01ZM and~02YJ]{Stacks}.
The local system $\HH_{\Zl}$ is then obtained by
pullback from the descended curve family, by proper
base change.

Choose an algebraic closure $\overline{k_0}\subset k$.
If $S_0$ is a model of a connected base $S$, then
$S_0$ is geometrically connected: a disconnection of
$S_0\times_{k_0}\overline{k_0}$ would persist after
base change to $k$. Let
\textcolor{black}{$\rho_{S_0,\ell}:\pi_1(S_0,\bar s)\to\GSp(H_{\Zl})$
be the representation of the descended family. Its multiplier is the
$\ell$-adic cyclotomic character. Define} its geometric image by
\[
 \Gamma^{\geom}_{S_0,\ell}
 =\rho_{S_0,\ell}
   \bigl(\pi_1(S_0\times_{k_0}\overline{k_0})\bigr).
\]
Here the geometric fundamental group is viewed as
the subgroup in the arithmetic exact sequence
\cite[Tag~0BTX]{Stacks}. \textcolor{black}{This subgroup acts trivially
on $\Zl(1)$, so its image lies in $\Sp(H_{\Zl})$.} Applying
Lemma~\ref{lem:extension-algebraically-closed-constants}
over $\overline{k_0}\subset k$ gives
\[
 \Gamma^{\geom}_{S_0,\ell}=\Gamma_{S,\ell}.
\]
If $k_1/k_0$ is a finite extension inside $k$, then
$\overline{k_0}$ is also an algebraic closure of $k_1$.
Thus replacing $S_0$ by $S_0\times_{k_0}k_1$ leaves
its geometric monodromy image unchanged.

\subsection{Reduction to unpointed moduli}
\label{subsec:reduction-unpointed}

We now reduce the monodromy assertions for pointed curves
to the unpointed case. The geometric reduction is a
standard dimension argument. We include its proof to
identify precisely the divisors to which the unpointed
monodromy statement will be applied.
Over $\C$, Hain likewise distinguishes the case
in which the divisor dominates the base from the case in which the
closure of its image is a divisor \cite[proof of Lem.~3.10]{HainDensity}.

We use the same notation $r_n$ for the forgetful morphism
at level $m$:
\[
 r_n:\cM_{g,n/k}[m]\longrightarrow\cM_{g/k}[m].
\]
By Subsection~\ref{subsec:moduli-universal-curves},
it is smooth and surjective with geometrically integral
fibers of dimension $n$. The underlying universal complete curve
over $\cM_{g,n/k}[m]$ is the pullback of the universal
curve over $\cM_{g/k}[m]$. Proper base change therefore
identifies the corresponding local systems $\HH_{\Zl}$
under pullback by $r_n$.

\begin{lemma}
\label{lem:divisor-forgetful}
Let $D\subset\cM_{g,n/k}[m]$ be a divisor, and let
$B$ be the reduced closure of $r_n(D)$ in
$\cM_{g/k}[m]$. Then either $B=\cM_{g/k}[m]$,
or $B$ is a divisor and
\[
 D=r_n^{-1}(B),
\]
where the inverse image is scheme-theoretic.
\end{lemma}

\begin{proof}
Let $\eta$ be the generic point of $B$.
Since $D$ is reduced, $r_n|_D$ factors through
$B$, and $D\to B$ is dominant. Its generic fiber
$D_\eta$ is a closed subscheme of the $n$-dimensional
fiber of $r_n$ over $\eta$, so
$\dim D_\eta\leq n$. The dimension formula gives
\[
 3g-4+n=\dim D=\dim B+\dim D_\eta
       \leq\dim B+n.
\]
Here the middle equality follows by expressing each
dimension as the transcendence degree of the corresponding
function field and using the tower law
\cite[Tag~0A21, part~(6)]{Stacks}.
Thus $\dim B\geq3g-4$. If $B\neq\cM_{g/k}[m]$,
it follows that $B$ is a divisor.

In this case, $r_n^{-1}(B)$ is flat over the
integral scheme $B$ and has integral generic fiber.
It is therefore integral \cite[Tag~0BCM]{Stacks},
and its dimension is $\dim B+n=3g-4+n$.
The divisor $D$ is a closed integral subscheme of
$r_n^{-1}(B)$ of the same dimension. Since
$r_n^{-1}(B)$ is integral, $D=r_n^{-1}(B)$.
\end{proof}

We follow Hain's argument for the pullback of monodromy along a lean
morphism \cite[Lem.~3.10]{HainDensity}. We give the argument for the
forgetful morphism $r_n$ to track both Zariski density and openness of
the integral $\ell$-adic image when passing from unpointed to pointed
curves.

\begin{proposition}
\label{prop:reduction-unpointed}
Suppose that, for every divisor $B\subset\cM_{g/k}[m]$,
the family over its normalization $B^\nu$ has full
symplectic monodromy at $\ell$. Then, for every
$n\geq0$, every divisor $D\subset\cM_{g,n/k}[m]$, and every
dominant morphism $E\to D$ from a smooth integral
$k$-variety, the family over $E$ has full symplectic
monodromy at $\ell$.

The same implication holds with full symplectic
monodromy replaced throughout by openness of the
integral monodromy image.
\end{proposition}

\begin{proof}
Let $B$ be as in Lemma~\ref{lem:divisor-forgetful}.
If $B=\cM_{g/k}[m]$, the composite
$E\to\cM_{g/k}[m]$ is dominant. By
Lemma~\ref{lem:ambient-monodromy}, the family over
$\cM_{g/k}[m]$ has full symplectic monodromy and
open integral image. Since this base is normal,
Corollary~\ref{cor:dominant-monodromy} gives both
assertions for the family over $E$.

If $B\neq\cM_{g/k}[m]$, then $B$ is a divisor and the
composite $E\to B$ is dominant. Applying
Corollary~\ref{cor:normalization-monodromy} in the
unpointed case, this morphism factors through a
dominant morphism $E\to B^\nu$, and
\[
 [\Gamma_{B^\nu,\ell}:\Gamma_{E,\ell}]<\infty.
\]
\textcolor{black}{After forgetting the markings, the curve over $E$ is
the pullback of the universal curve restricted to $B^\nu$. Its
symplectic local system is therefore pulled back from $B^\nu$.}
The same corollary shows
that each of the two monodromy assertions over $E$
is equivalent to the corresponding assertion over
$B^\nu$.
\end{proof}

Together with Corollary~\ref{cor:reduction-fine-level},
this proposition reduces the assertions for divisors
in all the pointed moduli stacks to those for the
normalizations of divisors in the unpointed fine
level scheme.

\section{\texorpdfstring{\textcolor{black}{Endomorphisms and the Kodaira--Spencer map}}{Endomorphisms and the Kodaira--Spencer map}}
\label{sec:endomorphisms-kodaira-spencer}

\textcolor{black}{We retain the notation and hypotheses of
Section~\ref{sec:moduli-monodromy}. We show that a Rosati-symmetric
endomorphism of a family of Jacobians produces quadratic
differentials annihilating the image of the Kodaira--Spencer map
(Corollary~\ref{cor:ks-endomorphism-annihilation}).}
For complex Jacobians isogenous to a product,
Ciliberto and van der Geer use products of differentials from the two
factors to study the cotangent space to the locus on which the splitting
persists and deduce a lower bound on its codimension
\cite[Prop.~2.2 and Cor.~2.3, p.~555]{CilibertoVanDerGeer1992}.
See also the tangent-space calculation of Ciliberto, van der Geer,
and Teixidor i Bigas \cite[Thm.~1.2]{CilibertoVanDerGeerTeixidor1992}. In characteristic $p$, the estimate
in Lemma~\ref{lem:products-differentials} applies to every nonscalar
operator on regular differentials.

\subsection{Products of canonical differentials}
\label{subsec:products-differentials}
We now give an estimate for quadratic differentials arising
from a nonscalar action on canonical differentials.
It will be used in the Kodaira--Spencer argument below.
\textcolor{black}{Products of regular differentials are taken in the
tensor square of the sheaf of differentials; thus they are
quadratic differentials.}

\begin{lemma}
\label{lem:products-differentials}
Let $C/k$ be a smooth proper connected curve of genus $g$,
\textcolor{black}{and let}
{\color{black}
\[
 T:H^0(C,\Omega^1_{C/k})\longrightarrow H^0(C,\Omega^1_{C/k})
\]
}
\textcolor{black}{be a nonscalar $k$-linear endomorphism. Then the
$k$-linear map}
{\color{black}
\[
 \begin{aligned}
 H^0(C,\Omega^1_{C/k})\otimes_k H^0(C,\Omega^1_{C/k})
 &\longrightarrow H^0(C,(\Omega^1_{C/k})^{\otimes2}),\\
 u\otimes v&\longmapsto u\,T(v)-T(u)\,v
 \end{aligned}
\]
}
\textcolor{black}{has image of dimension at least $g-1$.
Equivalently, this is the dimension bound for the subspace spanned
by the quadratic differentials $u\,T(v)-T(u)\,v$.}
\end{lemma}

\begin{proof}
Since $k$ is algebraically closed, choose an eigenvalue
$c\in k$ of $T$, and put
\[
 U=\ker(T-c\,\mathrm{id}),\qquad
 V=\operatorname{im}(T-c\,\mathrm{id}).
\]
\textcolor{black}{Both are subspaces of $H^0(C,\Omega^1_{C/k})$.}
The choice of $c$ gives $U\neq0$, and nonscalarity
of $T$ gives $V\neq0$. By rank--nullity,
\[
 \dim_k U+\dim_k V=g.
\]
For $u\in U$ and $v\in V$, choose
$w\in H^0(C,\textcolor{black}{\Omega^1_{C/k}})$ with
$v=(T-c\,\mathrm{id})(w)$. Since $T(u)=cu$, we have
\[
 uv=u\,T(w)-T(u)\,w.
\]
Thus every product of an element of $U$ with an
element of $V$ belongs to the \textcolor{black}{image} in the
statement. It remains to exhibit $g-1$ linearly
independent such products.

Choose a point $P\in C(k)$ \textcolor{black}{and} a local parameter \textcolor{black}{$t$} at $P$.
\textcolor{black}{Since $C$ is smooth, $dt$ is a local basis of
$\Omega^1_{C/k}$ near $P$. For $u\neq0$, write $u=t^m a\,dt$
with $a$ regular and $a(P)\neq0$; by definition,
$\operatorname{ord}_P(u)=m$. We set $\operatorname{ord}_P(0)=+\infty$.
Changing the local parameter multiplies the local basis by a unit,
so the order is independent of this choice.}
\textcolor{black}{For $n\geq0$, define the vanishing filtrations by}
{\color{black}
\[
 U_n=\{u\in U:\operatorname{ord}_P(u)\geq n\},\qquad
 V_n=\{v\in V:\operatorname{ord}_P(v)\geq n\}.
\]
}
\textcolor{black}{For $u\in U_n$, the quotient $u/(t^n\,dt)$ is
regular at $P$. Its value at $P$ vanishes precisely when
$u\in U_{n+1}$. Hence the leading coefficient defines an injection}
{\color{black}
\[
 U_n/U_{n+1}\hookrightarrow k,\qquad
 [u]\longmapsto\left(\frac{u}{t^n\,dt}\right)(P).
\]
}
\textcolor{black}{The same construction applies to $V_n/V_{n+1}$.}
The filtrations on $U$ and $V$ by order of vanishing
at $P$ have graded pieces of dimension at most one,
since the leading coefficient identifies each graded
piece with a subspace of $k$. These filtrations are
separated, because a nonzero differential has finite
order of vanishing at $P$.
\textcolor{black}{Since $U$ and $V$ are finite-dimensional, the
filtrations eventually reach zero. Choosing a lift from each
nonzero graded piece gives a basis with strictly increasing
vanishing orders.}

Put $r=\dim_k U$ and $\textcolor{black}{d}=\dim_k V$. We may therefore
choose bases $u_1,\ldots,u_r$ of $U$ and
$v_1,\ldots,v_{\textcolor{black}{d}}$ of $V$ whose vanishing orders at $P$
are strictly increasing. The $r+\textcolor{black}{d}-1$ products
\[
 u_1v_1,\ldots,u_1v_{\textcolor{black}{d}},\quad
 u_2v_{\textcolor{black}{d}},\ldots,u_rv_{\textcolor{black}{d}}
\]
also have strictly increasing vanishing orders,
because orders add under multiplication.
They are consequently linearly independent:
in any nontrivial linear combination, the term
of least vanishing order cannot cancel.
There are $r+\textcolor{black}{d}-1=g-1$ of them, proving the assertion.
\end{proof}

\subsection{\texorpdfstring{\textcolor{black}{Hodge bundles and endomorphisms}}{Hodge bundles and endomorphisms}}
\label{subsec:endomorphisms-kodaira-spencer}

An endomorphism of a family of Jacobians must be compatible with the
variation of its regular differentials. We first express this
compatibility using the Gauss--Manin connection. We then describe the
same variation in terms of the Kodaira--Spencer class of the family
of curves. Combining the two descriptions will show that certain
quadratic differentials from Lemma~\ref{lem:products-differentials}
annihilate the tangent directions of the family.

\subsubsection*{The Hodge bundle and its quotient}

Let $f:C\to S$ be a smooth proper family of geometrically connected
curves of genus $g$ over a smooth $k$-scheme $S$. Let
$q:J=\operatorname{Jac}(C/S)\to S$ be its relative Jacobian, with
canonical principal polarization $\lambda:J\to J^\vee$. Put
\[
 \mathcal H=R^1q_*\Omega^\bullet_{J/S},
 \qquad E=q_*\Omega^1_{J/S}.
\]
The notation $R^1q_*\Omega^\bullet_{J/S}$ denotes the first
hypercohomology sheaf of the relative de Rham complex. Its fiber at
$s\in S(k)$ is $H^1_{\mathrm{dR}}(J_s/k)$, whereas
$E_s=H^0(J_s,\Omega^1_{J_s/k})$.
The Hodge--de Rham spectral sequence for an abelian scheme degenerates
at $E_1$ \cite[Proposition~2.5.2]{BBM}. In degree one, this gives
\begin{equation}
\label{eq:jacobian-hodge-sequence}
 0\longrightarrow E\longrightarrow\mathcal H
 \longrightarrow R^1q_*\mathcal O_J\longrightarrow0.
\end{equation}
Thus the Hodge quotient is $R^1q_*\mathcal O_J$. The ranks of the three
bundles are $g$, $2g$, and $g$, respectively.

We will use two identifications, and specify their origins.
Write $\operatorname{Lie}J=\operatorname{Lie}(J/S)$
for the relative Lie algebra of $J$.
First, evaluating an invariant differential at the identity identifies
\[
 E\simeq(\operatorname{Lie}J)^\vee
 \quad\text{\cite[proof of Lemma~2.5.3, p.~100]{BBM}}.
\]
The tangent space of the Picard variety of a curve is its first
cohomology group with coefficients in the structure sheaf. In families,
this is the identification
$\operatorname{Lie}J\simeq R^1f_*\mathcal O_C$
\cite[\S8.4, Theorem~1(a)]{BLR}.
For example, on a fiber a first-order deformation of the trivial
line bundle has transition functions $1+\epsilon h_{ij}$, with
$\epsilon^2=0$. The gluing condition is
$h_{jk}-h_{ik}+h_{ij}=0$, and changing trivializations changes
$(h_{ij})$ by a coboundary. Its tangent vector is therefore a class
in $H^1(C_s,\mathcal O_{C_s})$.
Relative Serre duality now gives
\[
 E\simeq(R^1f_*\mathcal O_C)^\vee
   \simeq f_*\Omega^1_{C/S}
 \quad\text{\cite[Tag~0E6N]{Stacks}}.
\]
We use this identification to regard elements of $E_s$ as regular
differentials on $C_s$.

Second, the same Picard tangent-space description gives
$R^1q_*\mathcal O_J\simeq\operatorname{Lie}(J^\vee)$
\cite[\S5.1.1]{BBM}. We identify the Hodge quotient with $E^\vee$
by using the polarization:
\[
 \begin{aligned}
 \mathcal H/E
 &\simeq R^1q_*\mathcal O_J
 \simeq\operatorname{Lie}(J^\vee)\\
 &\xrightarrow{(d\lambda)^{-1}}\operatorname{Lie}J
 \simeq E^\vee.
 \end{aligned}
\]
Thus a quotient class is first a tangent vector of $J^\vee$, and
$(d\lambda)^{-1}$ makes it a tangent vector of $J$. It then acts on
$E$ by evaluation. This particular choice of identification determines
the sign in Lemma~\ref{lem:second-fundamental-form-ks}.

\subsubsection*{Differentiating a regular differential}

The bundle $\mathcal H$ carries the Gauss--Manin connection
\cite[Theorem~1]{KatzOda1968}
\[
 \nabla:\mathcal H\longrightarrow
 \mathcal H\otimes_{\mathcal O_S}\Omega^1_{S/k}.
\]
For a local section $u$ and a local function $h$ on $S$, it satisfies
\[
 \nabla(hu)=h\nabla u+u\otimes d_{S/k}h.
\]
If $\xi\in T_sS$, evaluating the cotangent factor on $\xi$ gives
$(\nabla u)_s(\xi)\in H^1_{\mathrm{dR}}(J_s/k)$.
This is the derivative of the de Rham section in the direction $\xi$.
We recall the cochain construction of this connection below, before
using it for a local computation.

A section of $E$ is a family of regular differentials. Its derivative
need not lie in $E$, so we project that derivative to the Hodge
quotient. This defines the second fundamental form:
\begin{equation}
\label{eq:hodge-second-fundamental-form}
 \theta:E\longrightarrow
 (\mathcal H/E)\otimes_{\mathcal O_S}\Omega^1_{S/k}
 \simeq E^\vee\otimes_{\mathcal O_S}\Omega^1_{S/k}.
\end{equation}
This map is $\mathcal O_S$-linear. Indeed, if $u\in E$, then
$u\otimes d_{S/k}h$ maps to zero in the quotient, and the Leibniz
rule therefore gives $\theta(hu)=h\theta(u)$.

Let $\theta_s$ denote the fiber of $\theta$ at $s$. For $\xi\in T_sS$,
put
\[
 \theta_\xi
 =(\operatorname{id}_{E_s^\vee}\otimes\operatorname{ev}_\xi)\theta_s:
 E_s\longrightarrow E_s^\vee,
 \qquad
 \operatorname{ev}_\xi(\beta)=\beta(\xi).
\]
Concretely, for $x\in E_s$ choose a local extension $\widetilde x$ in
$E$, differentiate, evaluate on $\xi$, and take the quotient:
\[
 (\nabla\widetilde x)_s(\xi)
 \in H^1_{\mathrm{dR}}(J_s/k)
 \longmapsto \theta_\xi(x)\in E_s^\vee.
\]
We make explicit why this construction is independent of the chosen
local extension.

\begin{lemma}[Independence of the local extension]
\label{lem:projected-derivative-independent}
Let $s\in S(k)$, $\xi\in T_sS$, and $x\in E_s$.
The image of $(\nabla\widetilde x)_s(\xi)$ in the Hodge quotient,
identified with $E_s^\vee$, is independent of the local section
$\widetilde x$ of $E$ extending $x$.
\end{lemma}

\begin{proof}
Let $\widetilde x$ and $\widetilde x'$ be two such extensions.
After restricting to a common neighborhood of $s$, choose a local
frame $e_1,\ldots,e_g$ of $E$ and write
\[
 \widetilde x-\widetilde x'=\sum_r h_re_r.
\]
Since the two sections have the same value at $s$ and the vectors
$e_r(s)$ form a basis of $E_s$, we have $h_r(s)=0$ for every $r$.
The $\mathcal O_S$-linearity of $\theta$ therefore gives
\[
 \bigl(\theta(\widetilde x)-\theta(\widetilde x')\bigr)_s
 =\sum_r h_r(s)\bigl(\theta(e_r)\bigr)_s=0.
\]
Evaluating the cotangent factor on $\xi$ proves that the two
projected derivatives agree.
\end{proof}

We may therefore pair the resulting functional with a second
regular differential. Define the bilinear form
\[
 b_\xi(x,y)=\theta_\xi(x)(y),\qquad x,y\in E_s.
\]
Here $\theta_\xi(x)$ is a functional on $E_s$, and the right-hand
side means its value at $y$.

\subsubsection*{Compatibility with a Rosati-symmetric endomorphism}

Suppose that $a:J\to J$ is an $S$-endomorphism satisfying
\[
 a^\vee\circ\lambda=\lambda\circ a.
\]
This is the Rosati-symmetry condition. Pullback induces an
endomorphism $a^*$ of $\mathcal H$ preserving $E$; write
$T=a^*|_E$. To understand the effect on $\theta$, we need the actions
on both $E$ and \textcolor{black}{$\mathcal H/E$}.

An invariant differential $u$ is a functional on $\operatorname{Lie}J$.
Write $da$ for the differential of $a$ at the
identity section.
The definition of pullback gives
\[
 (Tu)(v)=u(da(v)),\qquad T=(da)^\vee.
\]
On the Hodge quotient $R^1q_*\mathcal O_J$, the action is
$d(a^\vee)$: the morphism $a^\vee$ of the dual abelian scheme is
pullback of degree-zero line bundles by $a$. After transporting
this action through the polarization, Rosati symmetry gives
\[
 (d\lambda)^{-1}\circ d(a^\vee)\circ d\lambda=da.
\]
Consequently, under $\mathcal H/E\simeq E^\vee$, the quotient action
is $da=T^\vee$. Recall that, for an open subset $U\subseteq S$, a section
$\varphi\in E^\vee(U)$ is an $\mathcal O_U$-linear map
$E|_U\to\mathcal O_U$. The dual map acts by precomposition:
\[
 (T^\vee\varphi)(u)=\varphi(Tu),\qquad u\in E(U).
\]

Pullback by an $S$-morphism commutes with the Gauss--Manin connection
\cite[\S2]{KatzOda1968}. Thus $a^*$ is horizontal, meaning
\[
 \nabla\circ a^*=(a^*\otimes1)\circ\nabla.
\]
This compatibility implies the following identity.

\begin{lemma}
\label{lem:bilinear-endomorphism}
Let $s\in S(k)$ and $\xi\in T_sS$, and let $T_s:E_s\to E_s$
be the fiber of $T$. Then
\begin{equation}
\label{eq:bilinear-endomorphism}
 b_\xi(T_s(x),y)=b_\xi(x,T_s(y))
 \qquad\text{for all }x,y\in E_s.
\end{equation}
\end{lemma}

\begin{proof}
For a local section $u$ of $E$, horizontality gives
$\nabla(Tu)=(a^*\otimes1)\nabla u$.
Taking the Hodge quotient, on which $a^*$ acts as $T^\vee$, yields
\begin{equation}
\label{eq:theta-endomorphism}
 \theta\circ T=(T^\vee\otimes1)\circ\theta.
\end{equation}
On the fiber at $s$, this says
$\theta_s(T_s(x))=(T_s^\vee\otimes1)\theta_s(x)$.

Evaluation of the cotangent factor on $\xi$ commutes with
$T_s^\vee$ on the first factor. To see this, take
$\varphi\in E_s^\vee$ and $\beta\in\Omega^1_{S/k}|_s$.
Applying $T_s^\vee$ first and then evaluating on $\xi$ sends
$\varphi\otimes\beta$ to $\beta(\xi)T_s^\vee\varphi$.
Evaluating first and then applying $T_s^\vee$ gives
$T_s^\vee(\beta(\xi)\varphi)$.
Since $\beta(\xi)\in k$ and $T_s^\vee$ is $k$-linear,
\[
 T_s^\vee\bigl(\beta(\xi)\varphi\bigr)
 =\beta(\xi)\,T_s^\vee\varphi.
\]
Thus the two orders agree. Consequently,
\[
 \theta_\xi(T_s(x))=T_s^\vee\theta_\xi(x).
\]
Evaluating this functional at $y$ gives
\[
 \begin{aligned}
 b_\xi(T_s(x),y)
 &=\theta_\xi(T_s(x))(y)\\
 &=(T_s^\vee\theta_\xi(x))(y)\\
 &=\theta_\xi(x)(T_s(y))
 =b_\xi(x,T_s(y)).
 \end{aligned}
\]
\end{proof}

\subsection{\texorpdfstring{\textcolor{black}{The Kodaira--Spencer map and Serre duality}}{The Kodaira--Spencer map and Serre duality}}
\label{subsec:ks-serre-duality}

We now describe the deformation of $C_s$ in the direction $\xi$.
Dualizing the relative cotangent sequence and restricting to the
fiber gives
\begin{equation}
\label{eq:curve-tangent-sequence}
 0\longrightarrow\mathcal T_{C_s/k}
 \longrightarrow\mathcal T_{C/k}|_{C_s}
 \longrightarrow\mathcal O_{C_s}\otimes_k T_sS
 \longrightarrow0
 \quad\text{\cite[Tag~02K4]{Stacks}}.
\end{equation}
The last map sends a tangent vector of the total space to its
direction on the base; its kernel consists of vertical tangent
vectors. Since $H^0(C_s,\mathcal O_{C_s})=k$, the connecting
homomorphism gives the Kodaira--Spencer map
\[
 \operatorname{KS}_s:T_sS
       \longrightarrow H^1(C_s,\mathcal T_{C_s/k}).
\]

We recall explicitly how this connecting homomorphism is computed.
Choose local lifts $\xi_i$ of $1\otimes\xi$ in
$\mathcal T_{C/k}|_{C_s}$ on an affine cover of $C_s$.
On a pairwise overlap put $v_{ij}=\xi_j-\xi_i$.
Both lifts have the same image on the base, so $v_{ij}$ is vertical.
Moreover, on a triple overlap,
\[
 v_{jk}-v_{ik}+v_{ij}
 =(\xi_k-\xi_j)-(\xi_k-\xi_i)+(\xi_j-\xi_i)=0.
\]
Thus $(v_{ij})$ is a \v{C}ech $1$-cocycle with values in
$\mathcal T_{C_s/k}$. Changing the lifts by vertical sections
$w_i$ changes $v_{ij}$ by $w_j-w_i$, a coboundary. By definition,
\[
 \operatorname{KS}_s(\xi)=[(v_{ij})].
\]
The brackets mean the cohomology class of the entire collection.

For a regular differential $x\in H^0(C_s,\Omega^1_{C_s/k})$,
evaluation on vector fields defines the sheaf map
\[
 \mathrm{ev}_x:\mathcal T_{C_s/k}\longrightarrow\mathcal O_{C_s},
 \qquad v\longmapsto x(v).
\]
The induced $k$-linear map on cohomology is
\[
 \begin{aligned}
 H^1(\mathrm{ev}_x):\quad
 H^1(C_s,\mathcal T_{C_s/k})
 &\longrightarrow H^1(C_s,\mathcal O_{C_s}),\\
 \kappa&\longmapsto\kappa\smile x.
 \end{aligned}
\]
In terms of \v{C}ech representatives, this means
\[
 H^1(\mathrm{ev}_x)\bigl([(v_{ij})]\bigr)
 =[(x(v_{ij}))].
\]
Thus $\kappa\smile x$ denotes cup product followed by contraction.
The next standard compatibility explains why products of two regular
differentials enter the argument.

\begin{lemma}[Cup product and Serre duality]
\label{lem:cup-product-serre-duality}
Denote the perfect Serre duality pairings by
\[
 \begin{aligned}
 \langle\, ,\,\rangle_{\mathcal O_{C_s}}:\;&
 H^1(C_s,\mathcal O_{C_s})\times H^0(C_s,\Omega^1_{C_s/k})
 \longrightarrow k,\\
 \langle\, ,\,\rangle_{\mathcal T_{C_s/k}}:\;&
 H^1(C_s,\mathcal T_{C_s/k})\times
 H^0(C_s,(\Omega^1_{C_s/k})^{\otimes2})\longrightarrow k.
 \end{aligned}
\]
These are the pairings of \cite[Remark~48.27.6]{Stacks}
with $X=C_s$, $d=1$, $i=0$, and
$\mathcal E=\mathcal O_{C_s}$ or
$\mathcal T_{C_s/k}$, respectively.
For $\kappa\in H^1(C_s,\mathcal T_{C_s/k})$ and
$x,y\in H^0(C_s,\Omega^1_{C_s/k})$, one has
\[
 \langle\kappa\smile x,y\rangle_{\mathcal O_{C_s}}
 =\langle\kappa,xy\rangle_{\mathcal T_{C_s/k}}.
\]
\end{lemma}

\begin{proof}
We first recall how the trace enters Serre duality. Since
$H^0(C_s,\mathcal O_{C_s})=k$, the trace is the composite
\[
 \operatorname{Tr}_{C_s}:H^1(C_s,\Omega^1_{C_s/k})
 \xrightarrow{\sim}H^0(C_s,\mathcal O_{C_s})^\vee
 \xrightarrow{\ell\mapsto\ell(1)}k,
\]
where the first arrow is Serre duality. Thus the trace pairs a class
with the constant section $1$.

For $y\in H^0(C_s,\Omega^1_{C_s/k})$, multiplication defines
\[
 m_y:\mathcal O_{C_s}\longrightarrow\Omega^1_{C_s/k},
 \qquad h\longmapsto hy,
\]
and induces the $k$-linear map
\[
 H^1(m_y):H^1(C_s,\mathcal O_{C_s})
 \longrightarrow H^1(C_s,\Omega^1_{C_s/k}).
\]
Under Serre duality, its dual sends $1$ to $y$: the section $1$
represents the identity of $\Omega^1_{C_s/k}$, and precomposing this
identity with $m_y$ gives $m_y$, which represents $y$.
Consequently, naturality of Serre duality gives, for
$c\in H^1(C_s,\mathcal O_{C_s})$,
\[
 \begin{aligned}
 \langle c,y\rangle_{\mathcal O_{C_s}}
 &=\langle H^1(m_y)(c),1\rangle_{\Omega^1_{C_s/k}}\\
 &=\operatorname{Tr}_{C_s}\bigl(H^1(m_y)(c)\bigr).
 \end{aligned}
\]
Here the middle pairing is the Serre pairing between
$H^1(C_s,\Omega^1_{C_s/k})$ and $H^0(C_s,\mathcal O_{C_s})$;
see \cite[Remarks~48.27.2--48.27.3]{Stacks}.

Now represent $\kappa$ by a \v{C}ech cocycle $(v_{ij})$.
Since $\kappa\smile x=[(x(v_{ij}))]$, multiplication by $y$ gives
\[
 \langle\kappa\smile x,y\rangle_{\mathcal O_{C_s}}
 =\operatorname{Tr}_{C_s}\bigl([(x(v_{ij})y)]\bigr).
\]
On the other hand, the quadratic differential $xy=x\otimes y$
corresponds to the sheaf map
\[
 \mathcal T_{C_s/k}\longrightarrow\Omega^1_{C_s/k},
 \qquad v\longmapsto x(v)y.
\]
Its induced map on $H^1$ sends $\kappa$ to the same class
$[(x(v_{ij})y)]\in H^1(C_s,\Omega^1_{C_s/k})$.
By the same naturality, pairing $\kappa$ with $xy$ applies the trace
to this class. Hence the two pairings agree.
\end{proof}

It remains to compare $\theta_\xi(x)$ with the class obtained by
evaluating $x$ on the Kodaira--Spencer cocycle. For this comparison we
need a cochain description of the Gauss--Manin connection. We briefly
review the construction of Katz and Oda
\cite[\S\S~2--3, pp.~202--206]{KatzOda1968} before making the local
calculation.

\subsection{\texorpdfstring{\textcolor{black}{The \v{C}ech--de Rham construction of the Gauss--Manin connection}}{The Cech--de Rham construction of the Gauss--Manin connection}}
\label{subsec:cech-gauss-manin}

\textcolor{black}{To compare the projected Gauss--Manin derivative
with the Kodaira--Spencer class, we need a cochain description of
both the connection and the Hodge quotient. We review these
constructions before applying them to the family of curves.}
Write $\nabla^C$ for the Gauss--Manin connection of
$f:C\to S$ and $\nabla^J=\nabla$ for the connection of $q:J\to S$
introduced in Subsection~\ref{subsec:endomorphisms-kodaira-spencer}.

\subsubsection*{The \v{C}ech--de Rham complex and the Hodge quotient}

Work over an affine neighborhood in $S$ and choose an affine open
cover $\mathcal U=(U_i)$ of $C$. Since $C$ is separated over this
affine base, finite intersections of these opens are affine.
The sheaves of differentials are quasi-coherent, so their
\v{C}ech complexes on this cover compute the sheaf cohomology needed
for de Rham hypercohomology.

There are two degrees in the \v{C}ech--de Rham complex.
An element of $\check C^r(\mathcal U,\Omega^b_{C/k})$ is a collection
of $b$-forms on intersections of $r+1$ opens. Its \emph{total degree}
is $r+b$. Thus the total complex is
\[
 \operatorname{Tot}^n
 =\bigoplus_{r+b=n}\check C^r(\mathcal U,\Omega^b_{C/k}).
\]
In particular,
\[
 \begin{aligned}
 \operatorname{Tot}^1
 &=\check C^0(\mathcal U,\Omega^1_{C/k})
       \oplus\check C^1(\mathcal U,\mathcal O_C),\\
 \operatorname{Tot}^2
 &=\check C^0(\mathcal U,\Omega^2_{C/k})
       \oplus\check C^1(\mathcal U,\Omega^1_{C/k})
       \oplus\check C^2(\mathcal U,\mathcal O_C).
 \end{aligned}
\]
A degree-one cochain is therefore a pair $((\omega_i),(h_{ij}))$:
one-forms on the individual opens and functions on their pairwise
overlaps. A degree-two cochain has three components: two-forms on
opens, one-forms on pairwise overlaps, and functions on triple overlaps.
Each component is an entire indexed collection.

Let $\delta$ be the usual alternating \v{C}ech differential.
In the degrees we need,
\[
 (\delta\omega)_{ij}=\omega_j-\omega_i,\qquad
 (\delta h)_{ijk}=h_{jk}-h_{ik}+h_{ij},
\]
with restrictions to the common overlap understood.
Let $d_{C/k}:\Omega^b_{C/k}\to\Omega^{b+1}_{C/k}$ denote the
$k$-linear exterior differential. For local functions $f_0,\ldots,f_b$,
it is given by
\[
 d_{C/k}(f_0\,df_1\wedge\cdots\wedge df_b)
 =df_0\wedge df_1\wedge\cdots\wedge df_b,
 \qquad df_i=d_{C/k}f_i.
\]
We denote its componentwise action on \v{C}ech cochains by
\[
 \begin{aligned}
 d_{\mathrm{dR}}:\check C^r(\mathcal U,\Omega^b_{C/k})
 &\longrightarrow\check C^r(\mathcal U,\Omega^{b+1}_{C/k}),\\
 (\alpha_{i_0\ldots i_r})
 &\longmapsto(d_{C/k}\alpha_{i_0\ldots i_r}).
 \end{aligned}
\]
Thus $d_{\mathrm{dR}}$ raises the form degree $b$ and preserves the
\v{C}ech degree $r$, whereas $\delta$ raises $r$ and preserves $b$.
On cochains of form degree $b$, set
\[
 d_{\mathrm{tot}}=d_{\mathrm{dR}}+(-1)^b\delta.
\]
The two unsigned differentials commute because exterior
differentiation commutes with restriction. The factor $(-1)^b$
makes their mixed terms cancel, so $d_{\mathrm{tot}}^2=0$.
This is the convention of \cite[\S3, pp.~205--206]{KatzOda1968}:
their \v{C}ech operator already contains this factor, and they
write the sum as $d+\delta$.

For a degree-one cochain, the definition gives
\begin{equation}
\label{eq:cech-derham-degree-one}
 d_{\mathrm{tot}}\bigl((\omega_i),(h_{ij})\bigr)
 =\bigl((d_{C/k}\omega_i),\,
        (d_{C/k}h_{ij}-(\omega_j-\omega_i)),\,
        (h_{jk}-h_{ik}+h_{ij})\bigr).
\end{equation}
The middle component is a sum because both differentiating the
overlap functions and taking differences of the local one-forms
produce one-forms on pairwise overlaps. The last component is
the \v{C}ech differential of the overlap functions.

We describe the Hodge quotient on the fiber for use
below. After applying $\nabla^C$ to a local section
and evaluating at $s$ in the direction $\xi$, we obtain a class in
$H^1_{\mathrm{dR}}(C_s/k)$. To compute its Hodge quotient, we use the
induced affine cover $\mathcal U_s=(U_i\cap C_s)$.
Affineness of all finite intersections and quasi-coherence of the
sheaves of forms give
\[
 \begin{aligned}
 H^1\!\left(\check C^\bullet(\mathcal U_s,\mathcal O_{C_s})\right)
 &\cong H^1(C_s,\mathcal O_{C_s}),\\
 H^1\!\left(\operatorname{Tot}\check C^\bullet
       (\mathcal U_s,\Omega^\bullet_{C_s/k})\right)
 &\cong H^1_{\mathrm{dR}}(C_s/k).
 \end{aligned}
\]
The first identification is the \v{C}ech comparison theorem
\cite[Lemma~30.2.6]{Stacks}; the second is its hypercohomology version
\cite[\S3, pp.~205--206]{KatzOda1968}.
Thus de Rham classes are represented by total cocycles modulo total
coboundaries. On the fiber, a pair $((\alpha_i),(g_{ij}))$ is a
cocycle precisely when
\[
 d_{C_s/k}\alpha_i=0,\qquad
 \alpha_j-\alpha_i=d_{C_s/k}g_{ij},\qquad
 g_{jk}-g_{ik}+g_{ij}=0.
\]
The first condition is automatic for a smooth curve, since
$\Omega^2_{C_s/k}=0$. The last condition says that $(g_{ij})$ is a
function cocycle. If the total pair is changed by the differential
of a degree-zero cochain $(h_i)$, its components change by
\[
 \bigl((d_{C_s/k}h_i),(h_j-h_i)\bigr).
\]
The function cocycle therefore changes by a \v{C}ech coboundary.
It follows that taking the function component defines a map on
cohomology:
\begin{equation}
\label{eq:cech-hodge-projection}
 \begin{aligned}
 H^1_{\mathrm{dR}}(C_s/k)&\longrightarrow H^1(C_s,\mathcal O_{C_s}),\\
 [(\alpha_i),(g_{ij})]&\longmapsto[(g_{ij})].
 \end{aligned}
\end{equation}
This is the Hodge quotient map: it is induced by projection of the
de Rham complex onto its degree-zero term $\mathcal O_{C_s}$.
By the long exact sequence of the two-term de Rham
complex, the image of \eqref{eq:cech-hodge-projection} is the kernel of
$d:H^1(C_s,\mathcal O_{C_s})\to H^1(C_s,\Omega^1_{C_s/k})$.
This differential vanishes for a smooth proper curve
\cite[Tag~0FW6]{Stacks}, so \eqref{eq:cech-hodge-projection} is surjective.
In particular, we must first obtain a total cocycle representing the
derivative, and only then take its function component.

\subsubsection*{How the total differential computes the connection}

We now explain which cochain represents that derivative.
Let $F^r\Omega^\bullet_{C/k}$ be the subcomplex of absolute forms
with at least $r$ differential factors pulled back from $S$.
For example,
\[
 F^1\Omega^1_{C/k}=f^*\Omega^1_{S/k},\qquad
 F^1\Omega^2_{C/k}=f^*\Omega^1_{S/k}\wedge\Omega^1_{C/k},
\]
and $F^2\Omega^2_{C/k}$ consists of sums of wedges of two base
one-forms. There is an exact sequence of complexes
\begin{equation}
\label{eq:gm-complex-sequence}
 0\longrightarrow
 \frac{F^1\Omega^\bullet_{C/k}}{F^2\Omega^\bullet_{C/k}}
 \longrightarrow
 \frac{\Omega^\bullet_{C/k}}{F^2\Omega^\bullet_{C/k}}
 \longrightarrow
 \Omega^\bullet_{C/S}\longrightarrow0.
\end{equation}
Its connecting homomorphism is
\[
 \partial:R^1f_*\Omega^\bullet_{C/S}
 \longrightarrow
 R^2f_*\left(F^1\Omega^\bullet_{C/k}/F^2\Omega^\bullet_{C/k}\right).
\]
Here these are hypercohomology sheaves of complexes, rather than
higher direct images of a single sheaf.

In $F^1/F^2$ exactly one base differential factor remains.
Writing it first gives the degreewise identification
\[
 F^1\Omega^{n+1}_{C/k}/F^2\Omega^{n+1}_{C/k}
 \simeq f^*\Omega^1_{S/k}\otimes_{\mathcal O_C}\Omega^n_{C/S}.
\]
The differential on the remaining relative forms has a minus sign,
because the base factor has degree one. Thus this is the relative
de Rham complex shifted by one degree, with a base cotangent factor.
Consequently, after interchanging the tensor factors in the target,
\begin{equation}
\label{eq:gm-connecting-identification}
 R^2f_*\left(F^1\Omega^\bullet_{C/k}/F^2\Omega^\bullet_{C/k}\right)
 \simeq R^1f_*\Omega^\bullet_{C/S}
              \otimes_{\mathcal O_S}\Omega^1_{S/k}.
\end{equation}
This isomorphism follows from \cite[\S2, p.~202, (7)]{KatzOda1968}
with $p=q=1$, followed by interchanging the tensor factors.
Under this isomorphism the connecting homomorphism
$\partial$ is the Gauss--Manin connection $\nabla^C$
\cite[\S\S2--3, pp.~202--206]{KatzOda1968}.
For the family $q:J\to S$, the same construction is the connection
\textcolor{black}{$\nabla^J$} used above.

We spell out the formula for this connecting homomorphism.
Let $z$ be a relative total cocycle of degree one and lift it to an
absolute cochain $\widehat z$. Since $z$ is closed, the image of
$d_{\mathrm{tot}}\widehat z$ in the relative complex is zero.
Hence $d_{\mathrm{tot}}\widehat z$, modulo $F^2$, is a cochain in
$F^1/F^2$. It is closed because $d_{\mathrm{tot}}^2=0$, and the
definition of the connecting homomorphism is
\begin{equation}
\label{eq:gm-cochain-formula}
 \partial[z]=[d_{\mathrm{tot}}\widehat z\bmod F^2].
\end{equation}
Modulo $F^2$, any two lifts differ by a cochain in
$F^1/F^2$, so their differentials differ
by a coboundary in that complex. This proves independence of the lift.

The cohomology class in \eqref{eq:gm-cochain-formula} is taken in
$F^1/F^2$. Although the expression is a differential in the larger
absolute complex, the image of its primitive
$\widehat z$ modulo $F^2$ need not lie in
$F^1/F^2$, so this class can be nonzero. Under
\eqref{eq:gm-connecting-identification} it is the connection
$\nabla^C[z]$. This is the precise relation between the total
differential on cochains and the connection on cohomology.

After restricting to $s$ and evaluating the base cotangent factor
on $\xi$, the class becomes an element of $H^1_{\mathrm{dR}}(C_s/k)$.
For a local section $h$ of
$R^1f_*\Omega^\bullet_{C/S}$ and a local one-form $\beta$ on $S$,
evaluation on $\xi$ sends
\[
 (h\otimes\beta)_s\longmapsto\beta_s(\xi)\,h_s.
\]
We then apply the Hodge quotient
\eqref{eq:cech-hodge-projection}. These are two successive operations:
first compute the connection class, then take its image in
$H^1(C_s,\mathcal O_{C_s})$.

This construction is functorial for morphisms over $S$.
Such a pullback preserves the filtration by base differentials and
commutes with $d_{\mathrm{tot}}$, hence with the connecting
homomorphism. It also commutes with projection of the de Rham
complex onto its degree-zero term. We will use both compatibilities
for the relative Abel map.

\subsection{\texorpdfstring{\textcolor{black}{The comparison formula and the annihilation corollary}}{The comparison formula and the annihilation corollary}}
\label{subsec:comparison-annihilation}

For a smooth proper family of curves, Katz identifies
the map induced by the Gauss--Manin connection on the graded pieces
of the Hodge filtration with cup product by the Kodaira--Spencer class
\cite[Thm.~3.5, part~(3.5.2)]{KatzNilpotent1970}. We include the
calculation to explain the minus sign that appears when we identify
the Hodge quotient of the Jacobian using its principal polarization.

\begin{lemma}[The second fundamental form and the Kodaira--Spencer class]
\label{lem:second-fundamental-form-ks}
Let $s\in S(k)$, $\xi\in T_sS$, and
$x\in H^0(C_s,\Omega^1_{C_s/k})$. Under the identifications
$E_s\simeq H^0(C_s,\Omega^1_{C_s/k})$ and
$E_s^\vee\simeq H^1(C_s,\mathcal O_{C_s})$ fixed above,
\[
 \theta_\xi(x)=-\operatorname{KS}_s(\xi)\smile x.
\]
\end{lemma}

\begin{proof}
We compute the projected derivative on the curve and then convert
Abel pullback into the polarization identification.

\smallskip
\noindent\emph{Step 1: Transfer the derivative to the curve.}
Choose $P\in C_s(k)$. Smoothness gives an immersion $j:Z\to C$
through $P$, with $Z\to S$ \'etale \cite[Tag~057G]{Stacks}.
The graph
\[
 (j,\operatorname{id}_Z):Z\longrightarrow C\times_S Z
\]
is a section of the base-changed family. Thus, after an \'etale
base change preserving the fiber and its base tangent space, we
\textcolor{black}{retain the notation $S,s,C,J$. Let $e:S\to C$
be the resulting section through $P$ and $i_e:C\to J$ its Abel map.
Write $i_{P,s}:C_s\to J_s$ for its fiber at $s$.}

Abel pullback gives the identification
\[
 i_{P,s}^*:H^0(J_s,\Omega^1_{J_s/k})
 \xrightarrow{\sim}H^0(C_s,\Omega^1_{C_s/k})
\]
used above \cite[Proposition~2.2 and its proof]{MilneJacobian1986}.
Put $x_J=(i_{P,s}^*)^{-1}(x)$, choose a local extension
$\widetilde x_J$ in $q_*\Omega^1_{J/S}$, and let
$\widetilde x=\textcolor{black}{i_e^*}\widetilde x_J$ be its pullback, a local section
of $f_*\Omega^1_{C/S}$.
Define $\eta$ by
\[
 \begin{aligned}
 H^1_{\mathrm{dR}}(J_s/k)&\longrightarrow H^1(J_s,\mathcal O_{J_s}),\\
 (\nabla^J\widetilde x_J)_s(\xi)&\longmapsto\eta.
 \end{aligned}
\]
Then $\theta_\xi(\textcolor{black}{x_J})=(d\lambda_s)^{-1}\eta$, independently of the
extension by Lemma~\ref{lem:projected-derivative-independent}.

Functoriality of the connection and the Hodge quotient gives
\[
 \begin{aligned}
 H^1_{\mathrm{dR}}(C_s/k)&\longrightarrow H^1(C_s,\mathcal O_{C_s}),\\
 (\nabla^C\widetilde x)_s(\xi)&\longmapsto i_{P,s}^*\eta.
 \end{aligned}
\]
We will compute this last class.

\smallskip
\noindent\emph{Step 2: Choose local lifts of the tangent direction.}
Over an affine neighborhood of $s$, choose an affine cover $(U_i)$
with relative \'etale coordinates $t_i$
\cite[Theorem~41.13.1]{Stacks}. The \'etale identification of
differentials gives
\[
 \Omega^1_{C/k}|_{U_i}
 \simeq (f^*\Omega^1_{S/k})|_{U_i}
       \oplus\mathcal O_{U_i}\,d_{C/k}t_i
 \quad\text{\cite[Tag~0H93]{Stacks}}.
\]
After restriction to $C_s$, evaluation on $\xi$ in the first
summand and zero on the second define
\[
 \xi_i\in
 \Gamma(U_i\cap C_s,\mathcal T_{C/k}|_{C_s}).
\]
Thus $\xi_i$ is a local lift of $1\otimes\xi$ in the tangent
sequence \eqref{eq:curve-tangent-sequence}, and
\[
 \xi_i(f^*h)=\xi(h),\qquad \xi_i(t_i)=0.
\]
Here $h$ is a local function on $S$. Each $\xi_i$ is a vector
field along the fiber, acting on functions on $U_i$ with values in
$\mathcal O_{U_i\cap C_s}$.
The differences $\xi_j-\xi_i$ annihilate base functions and are
therefore vertical. Our convention for the Kodaira--Spencer map gives
\[
 v_{ij}=\xi_j-\xi_i,
 \qquad [(v_{ij})]=\operatorname{KS}_s(\xi).
\]

\smallskip
\noindent\emph{Step 3: Compute the function cocycle in the Hodge quotient.}
Write
\[
 \widetilde x|_{U_i}=a_i\,d_{C/S}t_i,\qquad
 \widehat x_i=a_i\,d_{C/k}t_i.
\]
When evaluating absolute forms on $\xi_i$, we implicitly restrict
them to $C_s$. Since $(d_{C/k}t_i)(\xi_i)=\xi_i(t_i)=0$, we have
\[
 \widehat x_i(\xi_i)=a_i|_{C_s}\,\xi_i(t_i)=0.
\]

The global relative one-form $\widetilde x$ is closed because
$\Omega^2_{C/S}=0$. For $\mathcal U=(U_i)$, it therefore defines
the relative total cocycle
\[
 z=\bigl((\widetilde x|_{U_i})_i,(0)_{ij}\bigr)
 \in\check C^0(\mathcal U,\Omega^1_{C/S})
       \oplus\check C^1(\mathcal U,\mathcal O_C)
   =\operatorname{Tot}^1\check C^\bullet
       (\mathcal U,\Omega^\bullet_{C/S}).
\]
Its absolute lift is the total degree-one cochain
\[
 \widehat z=\bigl((\widehat x_i)_i,(0)_{ij}\bigr)
 \in\check C^0(\mathcal U,\Omega^1_{C/k})
       \oplus\check C^1(\mathcal U,\mathcal O_C)
   =\operatorname{Tot}^1\check C^\bullet
       (\mathcal U,\Omega^\bullet_{C/k}).
\]
Our total-differential convention
\eqref{eq:cech-derham-degree-one} gives
\[
 d_{\mathrm{tot}}\widehat z
 =\bigl((d_{C/k}\widehat x_i)_i,\,
        (-(\widehat x_j-\widehat x_i))_{ij},\,
        (0)_{ijk}\bigr).
\]
By the connecting-homomorphism formula
\eqref{eq:gm-cochain-formula}, its class modulo $F^2$ represents the
Gauss--Manin derivative under
\eqref{eq:gm-connecting-identification}, before evaluation on $\xi$.

Here a \emph{base differential factor} means a pulled-back one-form
$f^*\beta$, where $\beta$ is a local section of $\Omega^1_{S/k}$.
Modulo $F^2$, the components above are locally sums of terms
$f^*\beta\wedge\gamma$. Evaluating this factor on $\xi$ gives, for
$n=0,1$, the maps
\[
 \begin{aligned}
 \left.(F^1\Omega^{n+1}_{C/k}/F^2\Omega^{n+1}_{C/k})\right|_{C_s}
 &\longrightarrow\Omega^n_{C_s/k},\\
 [f^*\beta\wedge\gamma]|_{C_s}
 &\longmapsto\beta_s(\xi)\,\gamma_s,
 \end{aligned}
\]
where $\gamma$ is a local section of $\Omega^n_{C/k}$, $\gamma_s$
is its pullback to $C_s$, and
$\Omega^0_{C_s/k}=\mathcal O_{C_s}$. The restriction on the left
is restriction of the quotient sheaf. For example, if $u$ is a local
function on $S$ and $h$ a local function on $C$, then
\[
 \begin{aligned}
 f^*(d_{S/k}u)\wedge d_{C/k}t_i
 &\longmapsto \xi(u)\,d_{C_s/k}(t_i|_{C_s}),\\
 h\,f^*(d_{S/k}u)
 &\longmapsto h|_{C_s}\,\xi(u).
 \end{aligned}
\]
Thus the first component becomes local one-forms and the second
becomes overlap functions. The Hodge quotient
\eqref{eq:cech-hodge-projection} retains the cohomology class of
these functions.
Since $\widehat x_j-\widehat x_i$ has zero relative image, it is a
\textcolor{black}{section of $f^*\Omega^1_{S/k}$ on the overlap.
It vanishes on vertical tangent vectors, so its evaluation on a lift of
$1\otimes\xi$ is independent of the lift. Taking $\xi_i$ gives}
\[
 g_{ij}=-(\widehat x_j-\widehat x_i)(\xi_i),
 \qquad [(g_{ij})]=i_{P,s}^*\eta.
\]

Using $\widehat x_i(\xi_i)=\widehat x_j(\xi_j)=0$, we now calculate
\[
 \begin{aligned}
 g_{ij}
 &=-\widehat x_j(\xi_i)\\
 &=\widehat x_j(\xi_j-\xi_i)\\
 &=x(v_{ij}).
 \end{aligned}
\]
The last equality uses verticality of $\xi_j-\xi_i$, on which the
absolute lift and its relative image agree. Consequently,
\[
 i_{P,s}^*\eta=[(x(v_{ij}))]
             =\operatorname{KS}_s(\xi)\smile x.
\]

\smallskip
\noindent\emph{Step 4: Convert to the polarization identification.}
The Abel map $i_{P,s}:C_s\to J_s$, $Q\mapsto[Q-P]$, induces
pullback of algebraically trivial line bundles:
\[
 \operatorname{Pic}^0(i_{P,s}):
 \operatorname{Pic}^0(J_s)=J_s^\vee
 \longrightarrow \operatorname{Pic}^0(C_s)=J_s.
\]
\textcolor{black}{On $k$-points it sends} $[\mathcal L]$
to $[i_{P,s}^*\mathcal L]$.
Here the two identifications are the definitions of the dual
abelian variety and the Jacobian. The canonical principal
polarization $\lambda_s:J_s\xrightarrow{\sim}J_s^\vee$ satisfies
\[
 \operatorname{Pic}^0(i_{P,s})\circ\lambda_s=[-1]_{J_s}
 \quad\text{\cite[Lemma~6.9 and Remark~6.10(c)]{MilneJacobian1986}}.
\]
Thus $\operatorname{Pic}^0(i_{P,s})$ is also an isomorphism.

Recall the Picard tangent-space identifications
\[
 \operatorname{Lie}(J_s)\simeq H^1(C_s,\mathcal O_{C_s}),
 \qquad
 \operatorname{Lie}(J_s^\vee)\simeq H^1(J_s,\mathcal O_{J_s}).
\]
Under these identifications, the tangent maps at the identity are
\[
 \begin{aligned}
 d\lambda_s:\ H^1(C_s,\mathcal O_{C_s})
 &\xrightarrow{\sim}H^1(J_s,\mathcal O_{J_s}),\\
 d\bigl(\operatorname{Pic}^0(i_{P,s})\bigr)=i_{P,s}^*:\ 
 H^1(J_s,\mathcal O_{J_s})
 &\xrightarrow{\sim}H^1(C_s,\mathcal O_{C_s}).
 \end{aligned}
\]
For the second map, pulling back first-order transition functions
$1+\varepsilon h_{ab}$, with $\varepsilon^2=0$, gives
$1+\varepsilon i_{P,s}^*h_{ab}$; hence the tangent map is precisely
pullback on $H^1$. Differentiating the composition identity, and
using $d[-1]=-\operatorname{id}$, gives
\[
 i_{P,s}^*\circ d\lambda_s
 =-\operatorname{id}_{H^1(C_s,\mathcal O_{C_s})},
 \qquad (d\lambda_s)^{-1}=-i_{P,s}^*.
\]
Finally, $H^1(C_s,\mathcal O_{C_s})\simeq E_s^\vee$ by the Serre
duality identification fixed above. Applying the inverse map to
$\eta\in H^1(J_s,\mathcal O_{J_s})$ therefore gives
\[
 \theta_\xi(\textcolor{black}{x_J})
 =(d\lambda_s)^{-1}\eta
 =-i_{P,s}^*\eta
 =-\operatorname{KS}_s(\xi)\smile x.
\]
\textcolor{black}{Since $i_{P,s}^*x_J=x$, this is the stated formula
under the identification of $E_s$ with $H^0(C_s,\Omega^1_{C_s/k})$.}
The argument applies in positive characteristic, including
characteristic two, where the two signs coincide.
\end{proof}

Retain the family $f:C\to S$, its Jacobian
$J=\operatorname{Jac}(C/S)$, and the Rosati-symmetric
$S$-endomorphism $a:J\to J$ considered above.
Recall the canonical identifications fixed earlier:
\[
 E\simeq(\operatorname{Lie}J)^\vee
 \simeq(R^1f_*\mathcal O_C)^\vee
 \simeq f_*\Omega^1_{C/S}.
\]
They arise, respectively, from evaluation of invariant differentials
at the identity, the Picard tangent-space identification, and
relative Serre duality. Through their composite, we regard
$T=a^*|_E:E\to E$ as an endomorphism of $f_*\Omega^1_{C/S}$,
retaining the same notation. For $s\in S(k)$, its fiber is therefore
viewed as
\[
 T_s:H^0(C_s,\Omega^1_{C_s/k})
 \longrightarrow H^0(C_s,\Omega^1_{C_s/k}).
\]
The preceding lemmas now give the following constraint.

\begin{corollary}[Annihilation of the Kodaira--Spencer image]
\label{cor:ks-endomorphism-annihilation}
With this notation, for every $x,y\in H^0(C_s,\Omega^1_{C_s/k})$
and $\xi\in T_sS$,
\begin{equation}
\label{eq:ks-endomorphism-annihilation}
 \left\langle\operatorname{KS}_s(\xi),
          x\,T_s(y)-T_s(x)\,y\right\rangle_{\mathcal T_{C_s/k}}=0.
\end{equation}
Thus the subspace of $H^0(C_s,(\Omega^1_{C_s/k})^{\otimes2})$
spanned by these quadratic differentials is contained in the
annihilator of $\operatorname{im}(\operatorname{KS}_s)$ under
Serre duality. If $T_s$ is nonscalar, this subspace has dimension
at least $g-1$.
\end{corollary}

\begin{proof}
By definition, $b_\xi(x,y)=\theta_\xi(x)(y)$. Under the Serre
duality identification $E_s^\vee\simeq H^1(C_s,\mathcal O_{C_s})$,
this evaluation is the Serre duality pairing. Hence
\[
 \begin{aligned}
 b_\xi(x,y)
 &=\langle\theta_\xi(x),y\rangle_{\mathcal O_{C_s}}
 &&\text{(definition)}\\
 &=-\langle\operatorname{KS}_s(\xi)\smile x,y\rangle_{\mathcal O_{C_s}}
 &&\text{(Lemma~\ref{lem:second-fundamental-form-ks})}\\
 &=-\langle\operatorname{KS}_s(\xi),xy\rangle_{\mathcal T_{C_s/k}}
 &&\text{(Lemma~\ref{lem:cup-product-serre-duality})}.
 \end{aligned}
\]
Lemma~\ref{lem:bilinear-endomorphism} gives
$b_\xi(T_s(x),y)=b_\xi(x,T_s(y))$. Substituting the formula above for
each side, we obtain
\[
 -\langle\operatorname{KS}_s(\xi),T_s(x)\,y\rangle_{\mathcal T_{C_s/k}}
 =
 -\langle\operatorname{KS}_s(\xi),x\,T_s(y)\rangle_{\mathcal T_{C_s/k}}.
\]
By bilinearity of the pairing, this is equivalent to
\eqref{eq:ks-endomorphism-annihilation}. Since $\xi\in T_sS$ is
arbitrary, each of these quadratic differentials, and hence their
span, annihilates $\operatorname{im}(\operatorname{KS}_s)$.
If $T_s$ is nonscalar, Lemma~\ref{lem:products-differentials}
applied to $C_s$ and $T_s$ gives the lower bound $g-1$ for the
dimension of this span.
\end{proof}

\section{Geometric simplicity of the generic Jacobian}
\label{sec:geometric-simplicity}

We retain the notation and hypotheses of
Section~\ref{sec:moduli-monodromy}. In particular, $k$ is
algebraically closed of characteristic $p>0$ and $g\geq3$.
In this section we prove
Theorem~\ref{thm:divisor-geometric-simplicity}, which asserts
that the geometric generic Jacobian over every reduced
irreducible divisor in $\cM_{g/k}$ is simple. For the proof,
we fix an integer $m\geq3$ prime to $p$ and work on the
fine level scheme $\cM_{g/k}[m]$.

\textcolor{black}{For divisors meeting the ordinary locus, we combine
Corollary~\ref{cor:ks-endomorphism-annihilation} with a
Rosati-symmetric endomorphism acting nonscalarly on regular
differentials. The first subsection gives the criterion used
to obtain this endomorphism.} For divisors contained in the nonordinary locus, we use
the dimension of the $p$-rank strata
\cite[Thm.~2.3]{FaberVanDerGeer} and the existence of
an absolutely simple Jacobian in each irreducible component
of these strata \cite[Application~5.7]{AchterPries}, and
then pass to the geometric generic fibre by specialization.
We give the details of the endomorphism arguments to make
clear that they also apply in characteristic two.

\subsection{Rosati-symmetric endomorphisms of ordinary abelian varieties}
\label{subsec:ordinary-endomorphisms}

To apply
Corollary~\ref{cor:ks-endomorphism-annihilation} in the ordinary
case, we need a Rosati-symmetric endomorphism whose action on
regular differentials is nonscalar. The following results
supply this operator.

Let $(A,\lambda)$ be a principally polarized abelian variety
of dimension $h>0$ over $k$.
Here $\lambda:A\xrightarrow{\sim}A^\vee$ is the
given principal polarization.
Put $\End^0(A)=\End(A)\otimes_{\Z}\Q$.
The Rosati involution associated to $\lambda$ is
\[
 \End^0(A)\longrightarrow\End^0(A),\qquad
 a\longmapsto a^\dagger:=\lambda^{-1}a^\vee\lambda,
\]
where $a^\vee$ denotes the dual rational endomorphism
of the dual abelian variety $A^\vee$.
Since $\lambda$ is principal, this involution preserves
$\End(A)$. Write
\[
 R^+=\{a\in\End(A):a^\dagger=a\}.
\]
For $a\in\End(A)$, write
\[
 da:=d_0a:\operatorname{Lie}A\longrightarrow\operatorname{Lie}A
\]
for its $k$-linear tangent map at the identity $0\in A$.
The map $a\mapsto da$ is a homomorphism of abelian groups.

The group $\End(A)$ is finite free over $\Z$
\cite[Thm.~12.5]{MilneAbelian}. Thus $R^+$ is also
finite free, and clearing denominators gives
\[
 R^+\otimes_{\Z}\Q=\End^0(A)^{\dagger=1}.
\]
Here $\End^0(A)^{\dagger=1}$ denotes the
$\Q$-subspace of elements fixed by the Rosati involution.
For the standard properties of the Rosati involution,
see \cite[Sec.~17]{MilneAbelian}.

The following elementary calculation describes differentials
of endomorphisms of the connected factor of ordinary $p$-torsion.

\begin{lemma}[Endomorphisms of $\mu_p^h$ and their differentials]
\label{lem:mu-p-matrix-differential}
Every homomorphism of $k$-group schemes
$\varphi:\mu_p^h\to\mu_p^h$ is uniquely represented by a matrix
$(m_{ij})\in\mathrm{M}_h(\F_p)$ through the formula
\[
 \varphi:\mu_p^h\longrightarrow\mu_p^h,\qquad
 (z_j)_{j=1}^h\longmapsto
 \left(\prod_{j=1}^h z_j^{m_{ij}}\right)_{i=1}^h.
\]
Here $1\leq i,j\leq h$, with $j$ indexing the input coordinates
and $i$ the output coordinates; exponents are read modulo $p$.
The formula holds on points over every $k$-algebra.
Under $\operatorname{Lie}(\mu_p^h)\cong k^h$, its differential
at the identity $(1,\ldots,1)$ is
\[
 d\varphi:k^h\longrightarrow k^h,\qquad
 (v_j)_{j=1}^h\longmapsto
 \left(\sum_{j=1}^h m_{ij}v_j\right)_{i=1}^h.
\]
In particular, $d\varphi=0$ implies that $\varphi$ is the zero
group homomorphism, and $d\varphi=c\,\mathrm{id}$ for $c\in k$
implies $c\in\F_p$.
\end{lemma}

\begin{proof}
Cartier duality identifies
\[
 \operatorname{Hom}_{k\text{-}\mathrm{gp}}(\mu_p^h,\mu_p)
 \simeq
 \operatorname{Hom}_{k\text{-}\mathrm{gp}}
 \bigl(\Z/p\Z,(\Z/p\Z)^h\bigr).
\]
A homomorphism on the right is determined by the image of $1$,
a tuple $(m_j)_{j=1}^h$; its Cartier dual is
$(z_j)_{j=1}^h\mapsto\prod_{j=1}^h z_j^{m_j}$.
Applying this to each coordinate of $\varphi$ gives the
unique matrix $(m_{ij})$.

A tangent vector $(v_1,\ldots,v_h)\in k^h$ at the identity
is represented by $(1+\varepsilon v_1,\ldots,1+\varepsilon v_h)$
over $k[\varepsilon]/(\varepsilon^2)$.
These coordinates lie in $\mu_p$ because
$(1+\varepsilon v_j)^p=1$ in characteristic $p$.
Applying $\varphi$ to this infinitesimal point gives,
in the $i$th coordinate,
\[
 \prod_{j=1}^h(1+\varepsilon v_j)^{m_{ij}}
 =1+\varepsilon\sum_{j=1}^hm_{ij}v_j,
 \qquad 1\leq i\leq h.
\]
Taking the coefficients of $\varepsilon$ proves the displayed
formula for $d\varphi$. Its matrix is $(m_{ij})$, now viewed
over $k$ through $\F_p\subset k$.
If $d\varphi=0$, all entries vanish, so $\varphi$ is the zero
group homomorphism. If $d\varphi=c\,\mathrm{id}$, then
$m_{ii}=c$ for each $1\leq i\leq h$; since $h>0$ and
$m_{ii}\in\F_p$, we have $c\in\F_p$.
\end{proof}

We next recall the pairings used to apply Rosati symmetry.
Write
\[
 A[p]^D=\underline{\operatorname{Hom}}(A[p],\mathbb G_m)
\]
for the Cartier dual of $A[p]$. Cartier duality for
$[p]:A\to A$ gives a canonical identification
$A^\vee[p]\cong A[p]^D$
\cite[Thm.~11.1]{MilneAbelian}.
The resulting perfect evaluation pairing takes values
in $\mu_p$, since $A[p]$ is killed by $p$:
\[
 e_p:A[p]\times A^\vee[p]\longrightarrow\mu_p,
 \qquad(x,z)\longmapsto z(x).
\]
Here $z$ is viewed as a character through
$A^\vee[p]\simeq A[p]^D$. All pairing formulas below are
identities on points over arbitrary $k$-schemes.
Since $\lambda$ is a principal polarization, it induces a perfect pairing
\[
 e_p^\lambda:A[p]\times A[p]\longrightarrow\mu_p,
 \qquad (x,y)\longmapsto e_p(x,\lambda(y)).
\]
\begin{lemma}[Vanishing on the connected factor]
\label{lem:ordinary-connected-vanishing}
Suppose that $A$ is ordinary, and fix an isomorphism of
finite group schemes
\[
 A[p]\cong\mu_p^h\oplus(\Z/p\Z)^h.
\]
If $a\in R^+$ vanishes on the connected factor $\mu_p^h$,
then it vanishes on $A[p]$.
\end{lemma}

\begin{proof}
Such a decomposition exists by the standard structure
of ordinary $p$-torsion \cite[Sec.~2.1]{PriesPTorsion}.
There are no nonzero homomorphisms between the two factors.
Indeed, consider homomorphisms of $k$-group schemes
\[
 \varphi:\mu_p^h\longrightarrow(\Z/p\Z)^h,\qquad
 \psi:(\Z/p\Z)^h\longrightarrow\mu_p^h.
\]
The image of $\varphi$ is connected and contains the identity
$0$ of the \'etale target, so $\varphi$ factors through that
identity point. Thus $\varphi=0$.
Since the source of $\psi$ is reduced, $\psi$ factors through
the reduced subscheme
$(\mu_p^h)_{\mathrm{red}}=\operatorname{Spec}k$,
which is the identity subgroup $(1,\ldots,1)$.
Thus $\psi=0$ as well.

For $a\in\End(A)$, the identity $[p]\circ a=a\circ[p]$
shows that $a$ restricts to an endomorphism of $A[p]$.
The two composites
\[
 \begin{aligned}
 \mu_p^h&\hookrightarrow A[p]
 \xrightarrow{a|_{A[p]}}A[p]\twoheadrightarrow(\Z/p\Z)^h,\\
 (\Z/p\Z)^h&\hookrightarrow A[p]
 \xrightarrow{a|_{A[p]}}A[p]\twoheadrightarrow\mu_p^h
 \end{aligned}
\]
are therefore zero by the preceding vanishing.
Here the first and last arrows are the inclusions and
projections of the fixed decomposition.
The first vanishing says that $a$ maps $\mu_p^h$ into itself,
and the second says the same for $(\Z/p\Z)^h$.
Thus every endomorphism of $A$ preserves both factors of $A[p]$.

The restrictions of $e_p^\lambda$ to the two self-pairings
induce homomorphisms
\[
 \begin{aligned}
 \mu_p^h&\longrightarrow(\mu_p^h)^D
       \simeq(\Z/p\Z)^h,\\
 (\Z/p\Z)^h&\longrightarrow\bigl((\Z/p\Z)^h\bigr)^D
       \simeq\mu_p^h.
 \end{aligned}
\]
In each line the map sends $y$ to the character
$x\mapsto e_p^\lambda(x,y)$, with $x$ and $y$ in that same
factor. The displayed identifications follow by taking products of
the Cartier duality between $\mu_p$ and $\Z/p\Z$
\cite[Example~11.13]{MilneAlgebraicGroups}.
Both homomorphisms are zero by the vanishing proved above.
Equivalently,
\[
 e_p^\lambda(x,y)=1
 \quad\text{if $x,y$ both lie in $\mu_p^h$
 or both lie in $(\Z/p\Z)^h$.}
\]
Thus each factor pairs trivially with itself.

Perfectness of $e_p^\lambda$ means that
$A[p]\to A[p]^D$, sending $y$ to
$\bigl(x\mapsto e_p^\lambda(x,y)\bigr)$, is an isomorphism.
Since the two self-pairings vanish, this map is, after
interchanging the two target summands, the direct sum of
\[
 \begin{aligned}
 \mu_p^h&\longrightarrow\bigl((\Z/p\Z)^h\bigr)^D,
 &x&\longmapsto\bigl(y\longmapsto e_p^\lambda(y,x)\bigr),\\
 (\Z/p\Z)^h&\longrightarrow(\mu_p^h)^D,
 &y&\longmapsto\bigl(x\longmapsto e_p^\lambda(x,y)\bigr).
 \end{aligned}
\]
Their direct sum is an isomorphism, so each of these maps
is an isomorphism. This proves that the restricted pairing
\[
 e_p^\lambda:\mu_p^h\times(\Z/p\Z)^h\longrightarrow\mu_p
\]
is perfect.

Compatibility of evaluation with Cartier duality gives
$e_p(ax,z)=e_p(x,a^\vee z)$.
Indeed, the canonical identification $A^\vee[p]\simeq A[p]^D$
identifies $a^\vee|_{A^\vee[p]}$ with the Cartier dual of
$a|_{A[p]}$. By definition, the latter sends a character
$z$ to $z\circ(a|_{A[p]})$. Thus, viewing $z$ as a character,
we obtain
\[
 e_p(ax,z)=z(ax)=(a^\vee z)(x)=e_p(x,a^\vee z).
\]
Recall that $a^\dagger=\lambda^{-1}a^\vee\lambda$, so
$a^\vee\circ\lambda=\lambda\circ a^\dagger$ as maps
$A\to A^\vee$. Thus, for $a\in R^+$ and
$x,y\in A[p](T)$ over any $k$-scheme $T$, we have
\[
\begin{aligned}
 e_p^\lambda(ax,y)
 &=e_p(ax,\lambda(y))
 &&\text{(definition of $e_p^\lambda$)}\\
 &=e_p(x,a^\vee(\lambda(y)))
 &&\text{(preceding identity)}\\
 &=e_p(x,\lambda(a^\dagger y))
 &&\text{($a^\vee\circ\lambda=\lambda\circ a^\dagger$)}\\
 &=e_p^\lambda(x,a^\dagger y)
 &&\text{(definition of $e_p^\lambda$)}\\
 &=e_p^\lambda(x,ay)
 &&\text{($a^\dagger=a$ for $a\in R^+$).}
\end{aligned}
\]
In particular, if $a$ vanishes on $\mu_p^h$,
take $x$ in the connected factor and $y$ in the
\'etale factor. Then
\[
 e_p^\lambda(x,ay)=e_p^\lambda(ax,y)=1.
\]
By perfectness of the pairing between the two factors,
$a$ also vanishes on $(\Z/p\Z)^h$, and hence on $A[p]$.

\end{proof}

The preceding lemmas give the divisibility by $p$ needed below.

\begin{lemma}[Kernel of the differential map]
\label{lem:ordinary-symmetric-differential}
Suppose that $A$ is ordinary. The differential map
\[
 R^+\longrightarrow\End_k(\Lie A),\qquad a\longmapsto da,
\]
has kernel $pR^+$. If $a\in R^+$ and
$da=c\,\mathrm{id}$ for some $c\in k$, then $c\in\F_p$.
\end{lemma}

\begin{proof}
Let $a\in R^+$. Choose an ordinary decomposition of $A[p]$ as above.
The inclusions induce isomorphisms
\[
 \operatorname{Lie}(\mu_p^h)
 \xrightarrow{\sim}\operatorname{Lie}(A[p])
 \xrightarrow{\sim}\operatorname{Lie}A.
\]
Indeed, the \'etale factor has zero tangent
space, and
$\operatorname{Lie}(A[p])=\ker(d[p])=\operatorname{Lie}A$
because $d[p]=p\,\mathrm{id}=0$.
Under these identifications, the differential of the restriction
$a|_{\mu_p^h}:\mu_p^h\to\mu_p^h$ corresponds to $da$.
If $da=0$, Lemma~\ref{lem:mu-p-matrix-differential} gives
$a|_{\mu_p^h}=0$, and Lemma~\ref{lem:ordinary-connected-vanishing}
then gives $a|_{A[p]}=0$.
Thus $a$ kills $A[p]$ as a group scheme. The quotient
property of $[p]:A\to A$ gives a unique
$b\in\End(A)$ with
$a=b\circ[p]=[p]\circ b=pb$.
Since $\lambda$ is a principal polarization,
$\lambda^{-1}:A^\vee\to A$ is a morphism, so
$b^\dagger=\lambda^{-1}b^\vee\lambda$ belongs to $\End(A)$.
The Rosati involution is $\Z$-linear, so applying it to
$a=pb$ gives $a^\dagger=p b^\dagger$.
As $a\in R^+$, we also have $a^\dagger=a$, and hence
\[
 p(b^\dagger-b)=p b^\dagger-pb=a^\dagger-a=0.
\]
The abelian group $\End(A)$ is torsion-free
\cite[Lem.~12.2]{MilneAbelian}, so multiplication by $p$
is injective on this group. It follows that
$b^\dagger-b=0$, and therefore $b\in R^+$ by definition.
Thus $a=pb$ lies in $pR^+$.
Conversely, every element of
$pR^+$ has zero differential,
since $d(pb)=p\,db=0$ for $b\in R^+$.

If $da=c\,\mathrm{id}$, the same tangent-space identifications give
$d(a|_{\mu_p^h})=c\,\mathrm{id}$.
Lemma~\ref{lem:mu-p-matrix-differential} therefore gives
$c\in\F_p$.
\end{proof}

\begin{corollary}
\label{cor:ordinary-nonscalar-endomorphism}
If $A$ is ordinary and
$\dim_{\Q}\End^0(A)^{\dagger=1}>1$, then some
Rosati-symmetric endomorphism of $A$ acts nonscalarly
on $H^0(A,\Omega^1_{A/k})$.
\end{corollary}

\begin{proof}
Since $R^+$ is finite free over $\Z$, the hypothesis gives
\[
 \dim_{\F_p}(R^+/pR^+)
 =\operatorname{rank}_{\Z}R^+
 =\dim_{\Q}\End^0(A)^{\dagger=1}>1.
\]
By Lemma~\ref{lem:ordinary-symmetric-differential},
the differential map induces an injection
\[
 R^+/pR^+\hookrightarrow\End_k(\Lie A),
 \qquad[a]\longmapsto da
\]
of $\F_p$-vector spaces. The same lemma shows that
if $da=c\,\mathrm{id}$ for $a\in R^+$, then
$c\in\F_p$. Thus the scalar elements in the image
lie in the one-dimensional $\F_p$-subspace
$\F_p\,\mathrm{id}$, whereas the image has dimension
greater than one. Some $a\in R^+$ therefore has
nonscalar differential.

Evaluation at the identity gives the canonical
isomorphism
\[
 H^0(A,\Omega^1_{A/k})
 \xrightarrow{\sim}(\operatorname{Lie}A)^\vee,
 \qquad u\longmapsto\bigl(v\longmapsto u_0(v)\bigr).
\]
Here $u_0$ is the value of $u$ at the identity,
and $v\in\operatorname{Lie}A=T_0A$.
Pullback of differential forms along $a:A\to A$
defines the map
\[
 a^*:H^0(A,\Omega^1_{A/k})
 \longrightarrow H^0(A,\Omega^1_{A/k}).
\]
Under the evaluation isomorphism, it is dual to $da$, since
\[
 (a^*u)_0(v)=u_0(da(v)),\qquad v\in\operatorname{Lie}A.
\]
Since a linear map is scalar if and only if its dual is scalar,
$a^*$ is also nonscalar.
\end{proof}
\subsection{Divisors meeting the ordinary locus}
\label{subsec:ordinary-divisors}

We now turn to a divisor $D\subset\cM_{g/k}[m]$ whose
geometric generic Jacobian is ordinary. To use the
annihilation statement of
Corollary~\ref{cor:ks-endomorphism-annihilation}, we need
a Rosati-symmetric endomorphism with nonscalar action on
regular differentials. Corollary~\ref{cor:ordinary-nonscalar-endomorphism}
provides such an endomorphism whenever the Rosati-fixed
subspace of the rational endomorphism algebra has dimension
greater than one. The following lemma shows that failure
of simplicity implies this condition.

\begin{lemma}
\label{lem:nonsimple-symmetric-endomorphism}
Let $(A,\lambda)$ be a principally polarized abelian variety
of positive dimension over $k$. If $A$ is not simple, then
\[
 \dim_{\Q}\End^0(A)^{\dagger=1}\geq2.
\]
Consequently, if $A$ is also ordinary, some
Rosati-symmetric endomorphism of $A$ acts nonscalarly
on $H^0(A,\Omega^1_{A/k})$.
\end{lemma}

\begin{proof}
Choose a nonzero proper abelian subvariety $i:B\hookrightarrow A$.
The polarization $\lambda$ induces a polarization
\[
 \lambda_B:=i^\vee\circ\lambda\circ i:B\longrightarrow B^\vee,
\]
where $i^\vee:A^\vee\to B^\vee$ is pullback of
algebraically trivial line bundles. Indeed, if $\lambda$
is induced by an ample line bundle $L$ on $A$, then
$\lambda_B$ is induced by the ample line bundle $i^*L$
on $B$; see the proof of
\cite[Prop.~12.1]{MilneAbelian}.
In particular, $\lambda_B$ is an isogeny, so it has an inverse
\[
 \lambda_B^{-1}\in
 \operatorname{Hom}(B^\vee,B)\otimes_{\Z}\Q.
\]
This inverse is a rational homomorphism; the induced
polarization need not be principal.
Define
\[
 e:=i\lambda_B^{-1}i^\vee\lambda\in\End^0(A).
\]
The identity $i^\vee\lambda i=\lambda_B$ gives
\[
 e^2
 =i\lambda_B^{-1}(i^\vee\lambda i)\lambda_B^{-1}i^\vee\lambda
 =i\lambda_B^{-1}i^\vee\lambda=e,
 \qquad ei=i.
\]
Thus $e$ is an idempotent, and $ei=i$ shows that $e\neq0$.
We also have $e\neq\mathrm{id}_A$. Indeed, choose a positive
integer $N$ such that $u=N\lambda_B^{-1}$ is an actual
homomorphism $B^\vee\to B$. Then $Ne=i\,u\,i^\vee\lambda$,
so the following diagram of homomorphisms commutes.
\[
\begin{tikzcd}[column sep=large,row sep=large]
 A \arrow[r,"{u\,i^\vee\lambda}"] \arrow[d,"Ne"']
   & B \arrow[d,"i"] \\
 A \arrow[r,"\mathrm{id}_A"'] & A
\end{tikzcd}
\]
In particular, $Ne$ factors through the inclusion $i:B\hookrightarrow A$.
If $e=\mathrm{id}_A$, then $Ne=[N]_A$, the multiplication-by-$N$
map on $A$. Since this map is surjective, the factorization
would give
\[
 A=\operatorname{im}([N]_A)\subseteq i(B)\subsetneq A,
\]
contradicting the assumption that $B$ is a proper abelian
subvariety of $A$.

We next check that $e$ is Rosati symmetric.
Under the canonical identifications with the double duals,
polarizations satisfy $\lambda^\vee=\lambda$ and
$\lambda_B^\vee=\lambda_B$; see
\cite[Examples~1.9 and~2.2]{ConradPolarizations}.
Duality extends $\Q$-linearly and reverses composition, so
\[
 (\lambda_B^{-1})^\vee
 =(\lambda_B^\vee)^{-1}
 =\lambda_B^{-1}.
\]
Using this identity and $(i^\vee)^\vee=i$, we obtain
\[
\begin{aligned}
 e^\vee
 &=(i\lambda_B^{-1}i^\vee\lambda)^\vee\\
 &=\lambda^\vee(i^\vee)^\vee
   (\lambda_B^{-1})^\vee i^\vee\\
 &=\lambda i\lambda_B^{-1}i^\vee.
\end{aligned}
\]
Therefore
\[
 e^\vee\lambda
 =\lambda i\lambda_B^{-1}i^\vee\lambda
 =\lambda e,
 \qquad
 e^\dagger=\lambda^{-1}e^\vee\lambda=e.
\]
A scalar idempotent in $\End^0(A)$ is either $0$ or
$\mathrm{id}_A$, since $c^2=c$ for $c\in\Q$ implies
$c=0$ or $1$. Therefore $\mathrm{id}_A$ and $e$ are
linearly independent elements of
$\End^0(A)^{\dagger=1}$, proving the dimension bound.
If $A$ is ordinary, the final assertion follows from
Corollary~\ref{cor:ordinary-nonscalar-endomorphism}.
\end{proof}

Applied to a nonsimple ordinary geometric generic Jacobian,
this lemma supplies the required endomorphism. To use it
in the Kodaira--Spencer argument, we must next extend it
over a family after passing to a suitable finite \'etale
cover of a smooth dense open of $D$.

We first record the separable descent needed to construct
this cover.
For an integral $k$-scheme $T$, we take its
geometric generic point $\bar\eta$ to be $\Spec\overline{k(T)}$,
for a chosen algebraic closure of its function field $k(T)$.

\begin{lemma}
\label{lem:separable-descent-endomorphism}
Let $A$ be an abelian variety over a field $K$, and let
$\overline K$ be an algebraic closure of $K$. Every
endomorphism of $A_{\overline K}$ is defined over a finite
separable extension of $K$.
\end{lemma}

\begin{proof}
Let $a\in\End(A_{\overline K})$, and let
$K^{\mathrm{sep}}\subseteq\overline K$ be the separable
closure of $K$. Since $\overline K/K^{\mathrm{sep}}$ is
purely inseparable, Chow's theorem
\cite[Thm.~3.19]{ConradChowTrace} gives an isomorphism
\[
 \End_{K^{\mathrm{sep}}}(A_{K^{\mathrm{sep}}})
 \xrightarrow{\sim}
 \End_{\overline K}(A_{\overline K}).
\]
Thus $a$ is defined over $K^{\mathrm{sep}}$. By finite
presentation, it is already defined over a finite
subextension $L/K$ of $K^{\mathrm{sep}}/K$, which is
necessarily separable.
\end{proof}

We now pass from the finite separable extension of function
fields to a finite \'etale cover of a dense open subset of $D$.

\begin{lemma}
\label{lem:etale-extension-endomorphism}
Let $\mathcal J\to D$ be the relative Jacobian, with its
canonical principal polarization $\lambda$, and let
$\bar\eta$ be a geometric generic point of $D$.
Every endomorphism
\[
 a\in\End(\mathcal J_{\bar\eta})
\]
extends, after choosing a compatible geometric generic point,
to an endomorphism $a_S$ of $\mathcal J_S$ over a connected
finite \'etale cover
\[
 S\longrightarrow U
\]
of a smooth dense open subset $U\subseteq D$.
If $a$ is Rosati symmetric, then $a_S$ is Rosati symmetric.
\end{lemma}

\begin{proof}
Put $K=k(D)$. By
Lemma~\ref{lem:separable-descent-endomorphism}, $a$ descends
to an endomorphism $a_L$ of $\mathcal J_L$ over a finite
separable extension $L/K$.

Choose a smooth dense open subset $U\subseteq D$, and let
$S$ be its normalization in $L$. The morphism $S\to U$
is finite \cite[Tag~032L]{Stacks} and surjective.
Since $L/K$ is separable, it is \'etale over a dense open
subset of $U$ \cite[Tag~07ND]{Stacks}. Replacing $U$ by
this open subset and $S$ by its inverse image, we obtain
a connected finite \'etale cover.

In particular, $S$ is smooth and hence normal, and its
function field is $L$. By finite presentation, $a_L$
extends to a homomorphism
$a_V:\mathcal J_V\to\mathcal J_V$ over some dense open
subset $V\subseteq S$. Over a Noetherian normal base,
a homomorphism between abelian schemes extends uniquely
from any dense open subset
\cite[Prop.~I.2.7]{FaltingsChaiDegeneration}.
Applying this result to $a_V$ gives
\[
 a_S:\mathcal J_S\longrightarrow\mathcal J_S,
 \qquad (a_S)_L=a_L.
\]
Let $\lambda_S$ denote the pullback of $\lambda$ to $S$.
If $a$ is Rosati symmetric, the homomorphisms
$\textcolor{black}{a_S^\vee\lambda_S}$ and $\textcolor{black}{\lambda_S a_S}$ agree on the generic
fibre. Uniqueness of extension gives
\[
 \textcolor{black}{a_S^\vee\lambda_S=\lambda_S a_S},
\]
which is the required Rosati symmetry.
\end{proof}

We adapt the tangent-space argument for complex
Jacobians isogenous to a product used by Ciliberto and van der Geer
\cite[Prop.~2.2 and Cor.~2.3, p.~555]{CilibertoVanDerGeer1992} and by Ciliberto, van der
Geer, and Teixidor i Bigas \cite[Thm.~1.2]{CilibertoVanDerGeerTeixidor1992}.
In the ordinary case, Lemma~\ref{lem:nonsimple-symmetric-endomorphism}
shows that a nonsimple Jacobian admits a Rosati-symmetric endomorphism
whose action on regular differentials is nonscalar.
After extending this endomorphism to the family, Corollary~\ref{cor:ks-endomorphism-annihilation}
shows that the Kodaira--Spencer image has codimension at least $g-1$
wherever the action on regular differentials is nonscalar.

\Needspace{10\baselineskip}
\begin{proposition}
\label{prop:generic-jacobian-simple}
Let $k$ be an algebraically closed field of characteristic
$p>0$, let $g\geq3$, and let $m\geq3$ be an integer prime
to $p$. Let $D\subset\cM_{g/k}[m]$ be a divisor whose
geometric generic Jacobian is ordinary. Then its
geometric generic Jacobian is simple.
\end{proposition}

\begin{proof}
Let $\bar\eta$ be a geometric generic point of $D$, and
write $\mathcal J\to D$ for the relative Jacobian.
Suppose, for contradiction, that $\mathcal J_{\bar\eta}$
is not simple.
By Lemma~\ref{lem:nonsimple-symmetric-endomorphism}, there
is a Rosati-symmetric endomorphism $a$ whose action on
regular differentials is nonscalar.
Lemma~\ref{lem:etale-extension-endomorphism} extends $a$
to a Rosati-symmetric endomorphism $a_S$ over a connected
finite \'etale cover $S\to U$ of a smooth dense open
subset $U\subseteq D$.

Let $f:C_S\to S$ be the pullback of the universal curve
and $q:\mathcal J_S\to S$ its relative Jacobian.
Recall the identification
\[
 E=q_*\Omega^1_{\mathcal J_S/S}
 \simeq f_*\Omega^1_{C_S/S}.
\]
The endomorphism $a_S$ induces $T=a_S^*|_E$.
Its action on the geometric generic fibre is nonscalar.

To see that nonscalar action persists on an open subset,
choose an affine open $W\subseteq S$ on which the rank-$g$
bundle $E$ has a basis $e_1,\ldots,e_g$. In this basis,
write
\[
 T(e_j)=\sum_{i=1}^g b_{ij}e_i,
 \qquad b_{ij}\in\Gamma(W,\mathcal O_W).
\]
For a point $t\in W$, let $i_t:\Spec\kappa(t)\to W$ denote the
canonical morphism. Write $E_t=i_t^*(E|_W)$ and $T_t=i_t^*(T|_W)$.
If $b_{ij}(t)$ denotes the image of $b_{ij}$ in $\kappa(t)$, then
the matrix of $T_t$ on $E_t$ is $(b_{ij}(t))$. It is scalar if and
only if
\[
 b_{ij}(t)=0\quad(i\neq j),
 \qquad
 b_{ii}(t)-b_{11}(t)=0\quad(2\leq i\leq g).
\]
Indeed, these equations say precisely that the matrix
is $b_{11}(t)I_g$. Thus the scalar locus in $W$ is the
common zero locus of the regular functions $b_{ij}$ for
$i\neq j$ and $b_{ii}-b_{11}$ for $2\leq i\leq g$, so it
is closed in $W$.

Being scalar is independent of the chosen basis.
As such open sets $W$ cover $S$, the scalar locus is
therefore a closed subset of $S$. Its complement is
nonempty because the action of $T$ on the geometric
generic fibre is nonscalar. Hence $T_t$ is nonscalar
for every point $t$ in this nonempty open subset.
Since $S$ is of finite type over the algebraically closed
field $k$, we may choose a point $s\in S(k)$ in this
open subset.

Let $u\in U(k)$ be the image of $s$. Since $S\to U$ is
\'etale, its differential gives an isomorphism
\[
 T_sS\xrightarrow{\sim}T_uU.
\]
Since $C_S\to S$ is the pullback of the universal
curve, its fibre $C_s$ is canonically identified with the fibre of
the universal curve at $u$.
The Kodaira--Spencer map of the universal curve gives
the isomorphism
\[
 T_u\cM_{g/k}[m]
 \xrightarrow[\sim]{\operatorname{KS}^{\mathrm{univ}}_u}
 H^1(C_s,\mathcal T_{C_s/k}).
\]
This map sends a tangent vector at $u$ to the first-order
deformation class obtained by pulling back the universal
curve \cite[Prop.~6.11]{VistoliDeformations}.
Indeed, first-order deformations of $C_s$ are classified
by this cohomology group
\cite[Sec.~1, pp.~79--80]{DeligneMumford}, and the
prime-to-$p$ level structure introduces no additional
infinitesimal deformation data, since forgetting it is
\'etale.

\textcolor{black}{Functoriality} of the Kodaira--Spencer map under base
change \cite[Prop.~6.10(b)]{VistoliDeformations}
identifies $\operatorname{KS}_s$ with the composite
\[
 T_sS
 \xrightarrow{\sim}T_uU
 \lhook\joinrel\longrightarrow T_u\cM_{g/k}[m]
 \xrightarrow[\sim]{\operatorname{KS}^{\mathrm{univ}}_u}
 H^1(C_s,\mathcal T_{C_s/k}).
\]
Here $U$ is a smooth open subset of the divisor $D$, so
\[
 \dim_k T_uU=3g-4,
 \qquad
 \dim_k T_u\cM_{g/k}[m]=3g-3.
\]
Consequently, $\operatorname{KS}_s$ is injective and its
image has codimension one.

\Needspace{21\baselineskip}
Recall the perfect Serre duality pairing introduced in
Lemma~\ref{lem:cup-product-serre-duality},
\[
 H^1(C_s,\mathcal T_{C_s/k})
 \times
 H^0\!\left(C_s,(\Omega^1_{C_s/k})^{\otimes2}\right)
 \longrightarrow k.
\]
Under this pairing, the annihilator of
$\operatorname{im}(\operatorname{KS}_s)$ identifies with
the linear functionals on $H^1(C_s,\mathcal T_{C_s/k})$
that vanish on this image. These are precisely the
linear functionals that factor through the quotient,
so the annihilator is naturally isomorphic to
\[
 \left(
 H^1(C_s,\mathcal T_{C_s/k})
 /\operatorname{im}(\operatorname{KS}_s)
 \right)^\vee.
\]
Since $\operatorname{im}(\operatorname{KS}_s)$ has
codimension one, this quotient and its dual are
one-dimensional. On the other hand,
Corollary~\ref{cor:ks-endomorphism-annihilation} shows
that this annihilator contains the subspace spanned by
\[
 x\,T_s(y)-T_s(x)\,y,
 \qquad
 x,y\in H^0(C_s,\Omega^1_{C_s/k}).
\]
Because $T_s$ is nonscalar, the same corollary gives
dimension at least $g-1$ for this subspace. Hence
\[
 \dim_k\!\left[
 \left(
 H^1(C_s,\mathcal T_{C_s/k})
 /\operatorname{im}(\operatorname{KS}_s)
 \right)^\vee
 \right]
 =1\geq g-1\geq2,
\]
a contradiction. Therefore
$\mathcal J_{\bar\eta}$ is simple.
\end{proof}

\subsection{Divisors in the nonordinary locus}
\label{subsec:nonordinary-divisors}

We now treat divisors contained in the nonordinary locus.
For $0\leq f\leq g$, let
\[
 V_f\subseteq\cM_{g/k}[m]
\]
denote the reduced closed locus where the Jacobian has
$p$-rank at most $f$. By \cite[Thm.~2.3]{FaberVanDerGeer},
restricted to the smooth-curve locus and pulled back along
the finite \'etale level cover, $V_f$ is pure of
codimension $g-f$.

\begin{lemma}
\label{lem:nonordinary-divisor-prank}
Let $D\subseteq\cM_{g/k}[m]$ be a divisor whose geometric
generic Jacobian is nonordinary. Then $D$ is an
irreducible component of $V_{g-1}$, and its geometric
generic Jacobian has $p$-rank $g-1$.
\end{lemma}

\begin{proof}
Let $\eta$ be the generic point of $D$. Since the
geometric generic Jacobian is nonordinary, its $p$-rank
is less than $g$, hence at most $g-1$. Thus
$\eta\in V_{g-1}$. Since $V_{g-1}$ is closed, it contains
the closure of $\{\eta\}$, which is $D$. Hence
$D\subseteq V_{g-1}$.

The locus $V_{g-1}$ is pure of dimension $3g-4$.
Since $D$ is irreducible and
\[
 \dim D=3g-4=\dim V_{g-1},
\]
$D$ is an irreducible component of $V_{g-1}$. Moreover,
\[
 \dim V_{g-2}=3g-5<\dim D,
\]
so $D$ is not contained in $V_{g-2}$. Thus
$D\setminus V_{g-2}$ is a nonempty dense open subset
on which the Jacobians have $p$-rank exactly $g-1$.
In particular, the geometric generic Jacobian has
$p$-rank $g-1$.
\end{proof}

We next use a result of Achter--Pries to find a point
of $D$ with simple Jacobian. We will then pass from
this fibre to the geometric generic fibre by specialization.

\begin{lemma}
\label{lem:nonordinary-simple-fibre}
Let $D$ be as in Lemma~\ref{lem:nonordinary-divisor-prank}.
There exists
\[
 s\in(D\setminus V_{g-2})(k)
\]
whose Jacobian is simple.
\end{lemma}

\begin{proof}
Put $D^\circ=D\setminus V_{g-2}$.
By Lemma~\ref{lem:nonordinary-divisor-prank},
$D^\circ$ is an irreducible component of
$V_{g-1}\setminus V_{g-2}$, and hence is closed
in this stratum. The finite \'etale level cover restricts
to a finite \'etale morphism between the corresponding
$p$-rank-$(g-1)$ strata with and without level structure.
Thus the image $Z$ of $D^\circ$ is closed in the
$p$-rank-$(g-1)$ stratum of $\cM_{g/k}$, and we obtain
a finite surjective morphism
\[
 D^\circ\longrightarrow Z.
\]
In particular, $Z$ is irreducible of dimension $3g-4$,
the dimension of every
component of this stratum. Hence $Z$ is an irreducible
component of the stratum.

The stratum of $p$-rank exactly $g-1$ is defined
over $\F_p$, as follows from Newton polygon
semicontinuity \cite[Thm.~2.3.1 and Cor.~2.3.2]{KatzSlope1979}.
Consequently, $Z$ is obtained by base change from a
geometrically irreducible component defined over a finite
extension of $\F_p$.
Let $\overline{\F_p}$ denote the algebraic
closure of $\F_p$ in $k$. By
\cite[Application~5.7]{AchterPries}, the descended component has
an $\overline{\F_p}$-point with absolutely simple Jacobian.
Its base change to $k$ gives $z\in Z(k)$
with simple Jacobian, since simplicity is preserved under
extension of algebraically closed fields
\cite[Cor.~3.21]{ConradChowTrace}.

Since $D^\circ\to Z$ is finite and surjective, its fibre
over $z$ is a nonempty finite $k$-scheme. As $k$ is
algebraically closed, this fibre has a $k$-point, giving
a lift $s\in D^\circ(k)$.
The corresponding curve has the same Jacobian,
which is therefore simple.
\end{proof}

We now pass from the simple fibre supplied by
Lemma~\ref{lem:nonordinary-simple-fibre} to the geometric
generic fibre using the following standard specialization
argument.

\begin{lemma}
\label{lem:simple-specialization}
Let $A\to T$ be an abelian scheme over an integral scheme
of finite type over $k$, and let $\bar\eta$ be a geometric
generic point of $T$. If $A_t$ is simple for some
$t\in T(k)$, then $A_{\bar\eta}$ is simple.
\end{lemma}

\begin{proof}
Put $K=k(T)$. Suppose that $A_{\bar\eta}$ is not simple.
By Poincar\'e complete reducibility
\cite[Prop.~12.1]{MilneAbelian}, its rational
endomorphism algebra contains an element
\[
 e^2=e,\qquad e\ne0,1.
\]
\textcolor{black}{For an abelian variety $B$ over a field $F$, write
$\End_F^0(B)=\End_F(B)\otimes_{\Z}\Q$. After clearing denominators
in $e$, Lemma~\ref{lem:separable-descent-endomorphism} gives a finite
separable extension $L/K$ and an element $e_L\in\End_L^0(A_L)$
whose base change is $e$.}

Let $T'$ be the normalization of $T$ in $L$, and put
$A'=A\times_TT'$. Since $T$ is of finite type over a
field, $T'\to T$ is finite and surjective.
Choose $t'\in T'(k)$ above $t$. Then
\[
 A'_{t'}\cong A_t.
\]

By \cite[Prop.~I.2.7]{FaltingsChaiDegeneration}, as used
in the proof of Lemma~\ref{lem:etale-extension-endomorphism},
every endomorphism of $A_L$ extends uniquely to $A'$.
Moreover, restriction to $t'$ gives an injection
\[
 \operatorname{End}_{T'}(A')
 \hookrightarrow
 \operatorname{End}_k(A'_{t'}).
\]
Indeed, the rigidity lemma \cite[Prop.~20.1]{MilneAbelian}
implies that an endomorphism $\varphi$ which is zero on the
fibre at $t'$ factors through the structure morphism
$A'\to T'$, since $T'$ is connected. As $\varphi$ preserves
the zero section, the resulting section $T'\to A'$
is the zero section. Hence $\varphi=0$.

Tensoring with $\Q$, extension and restriction therefore
give an injective homomorphism of $\Q$-algebras
\[
 \textcolor{black}{\End_L^0(A_L)\hookrightarrow\End_k^0(A_t)}.
\]
The image $e_t$ of $\textcolor{black}{e_L}$ satisfies $e_t^2=e_t$ and
$e_t\ne0,1$. Since $A_t$ is simple,
$\textcolor{black}{\End_k^0(A_t)}$ is a division algebra
\cite[Sec.~12, p.~20]{MilneAbelian}, so
$e_t\ne0$ is invertible. Multiplying the equality
$e_t(e_t-1)=0$ on the left by $e_t^{-1}$ gives
$e_t-1=0$, a contradiction.
\end{proof}

We can now combine the ordinary and nonordinary cases to
obtain geometric simplicity over every divisor in $\cM_{g/k}$.

\Needspace{10\baselineskip}
\begin{theorem}[Geometric simplicity over divisors]
\label{thm:divisor-geometric-simplicity}
Let $k$ be an algebraically closed field of characteristic
$p>0$, and let $g\geq3$. For every reduced irreducible
divisor
\[
 D\subseteq\cM_{g/k},
\]
the geometric generic Jacobian is simple.
\end{theorem}

\begin{proof}
Choose an integer $m\geq3$ prime to $p$, and recall the
finite \'etale level cover
\[
 q_m:\cM_{g/k}[m]\longrightarrow\cM_{g/k}.
\]
Let $D'$ be an irreducible component of $q_m^{-1}(D)$,
with its reduced structure. Since $q_m$ is finite
\'etale, $D'$ is a divisor mapping surjectively onto $D$.
A geometric generic point of $D'$ maps to a geometric
generic point of $D$, so it suffices to prove the
assertion over $D'$.

Let $\mathcal J\to D'$ be the relative Jacobian, and
let $\bar\eta$ be a geometric generic point of $D'$.
If $\mathcal J_{\bar\eta}$ is ordinary, the assertion
follows from Proposition~\ref{prop:generic-jacobian-simple}.
Otherwise, Lemma~\ref{lem:nonordinary-simple-fibre}
gives a point $s\in D'(k)$ such that $\mathcal J_s$
is simple. Since $D'$ is integral of finite type over
$k$ and $\mathcal J\to D'$ is an abelian scheme,
Lemma~\ref{lem:simple-specialization} implies that
$\mathcal J_{\bar\eta}$ is simple.
\end{proof}

Over $\overline{\mathbf F}_p$, Urbanik and Yang
prove that every irreducible component of the locus of curves whose
Jacobians admit a prime-to-$p$ isogeny from a product of an elliptic
curve and a $(g-1)$-dimensional abelian variety has codimension $g-1$
\cite[Thm.~5.14]{UrbanikYangSpecialLoci}. For $g\geq3$, this already
excludes such a splitting at the geometric generic point of a divisor.
Theorem~\ref{thm:divisor-geometric-simplicity} excludes isogeny
decompositions with factors of arbitrary positive dimensions and
with no restriction on the degree of the isogeny.

\section{Rank-one degenerations along divisors}
\label{sec:rank-one-degenerations}

Section~\ref{sec:geometric-simplicity} establishes geometric
simplicity of the generic Jacobian over every divisor.
We now show that this Jacobian also admits, after a finite
field extension, semistable reduction of toric rank one at
a discrete valuation trivial on $k$. These two properties
will be used in Section~\ref{sec:monodromy-theorem} to control
its endomorphisms and local monodromy.

We first recall how the toric rank of the generalized
Jacobian of a stable curve is determined by its dual graph.
We then construct a compactification of the coarse moduli
space on which the closure of every divisor meets the
required boundary locus.

\subsection{Generalized Jacobians and toric rank}
\label{subsec:generalized-jacobians-toric-rank}

Let $C/k$ be a stable curve of arithmetic genus $g$, and write
\[
 \nu:\widetilde C=\coprod_{i=1}^{v}C_i\longrightarrow C
\]
for its normalization. Each $C_i$ is smooth and connected;
put $g_i=g(C_i)$.

The dual graph $\Gamma_C$ has one vertex for each
irreducible component and one edge for each node.
The edge joins the vertices containing the two branches
of the node, with a loop when both branches belong to
the same component.
Write $V(\Gamma_C)$ and $E(\Gamma_C)$ for its vertex
and edge sets, respectively, and use $|\cdot|$ to denote
their cardinalities. The first Betti number
\[
 b_1(\Gamma_C)=\dim_{\Q}H_1(\Gamma_C,\Q)
\]
counts the independent cycles in the graph.
Since $C$ is connected, so is $\Gamma_C$, and hence
\[
 b_1(\Gamma_C)
 =|E(\Gamma_C)|-|V(\Gamma_C)|+1.
\]

The generalized Jacobian
\[
 J_C=\operatorname{Pic}^0_{C/k}
\]
parametrizes line bundles of degree zero on every
irreducible component.
By \cite[Thm.~2.3]{AlexeevCompactifiedJacobians},
pullback along $\nu$ fits into an exact sequence
of algebraic groups
\[
 0\longrightarrow T_C
 \longrightarrow J_C
 \xrightarrow{\nu^*}
 \prod_{i=1}^{v}\operatorname{Jac}(C_i)
 \longrightarrow0,
 \qquad
 T_C\simeq\mathbb G_m^{\,b_1(\Gamma_C)}.
\]
Thus $J_C$ is semiabelian, and its \emph{toric rank}
is the dimension of $T_C$.

Since $\dim J_C=g$, taking dimensions gives
\[
 g=b_1(\Gamma_C)+\sum_{i=1}^{v}g_i.
\]
In particular, $J_C$ has toric rank one precisely when
\[
 b_1(\Gamma_C)=1,
 \qquad\text{equivalently}\qquad
 \sum_{i=1}^{v}g_i=g-1.
\]

\subsection{Coarse moduli spaces and the Torelli map}
\label{subsec:coarse-torelli}

Recall that $M_{g/k}$ denotes the coarse moduli space
of $\cM_{g/k}$. Let $\overline{\cM}_{g/k}$ be the moduli
stack of stable curves of genus $g$, and let
$\overline M_{g/k}$ be its coarse moduli space.
The coarse space $\overline M_{g/k}$ is a normal integral
projective $k$-scheme, by Deligne--Mumford
\cite[Sec.~3 and Thm.~5.2]{DeligneMumford} and the GIT
construction \cite[Thms.~2.0.1--2.0.2]{GiesekerModuliCurves}.
The smooth-curve locus gives a dense open immersion
\[
 M_{g/k}\hookrightarrow\overline M_{g/k}.
\]
Consequently, $M_{g/k}$ is a normal integral
quasi-projective $k$-scheme, and
\[
 \dim M_{g/k}=\dim\overline M_{g/k}=3g-3.
\]
We will use these coarse spaces for the compactification
and divisor arguments.

Let $\mathcal A_{g/k}$ be the moduli stack of principally
polarized abelian varieties of dimension $g$, with coarse
moduli space $A_{g/k}$. Sending a smooth curve to its
canonically polarized Jacobian defines a morphism
\[
 \cM_{g/k}\longrightarrow\mathcal A_{g/k},
 \qquad
 C\longmapsto\bigl(\operatorname{Jac}(C),\lambda_C\bigr).
\]
The universal property of the coarse moduli space induces
a morphism $t$ making the diagram
\[
\begin{tikzcd}[column sep=large,row sep=large]
 \cM_{g/k} \arrow[r] \arrow[d]
   & \mathcal A_{g/k} \arrow[d] \\
 M_{g/k} \arrow[r,"t"'] & A_{g/k}
\end{tikzcd}
\]
commute, where the vertical arrows are the coarse moduli
morphisms. Thus
\[
 t([C])=\bigl[\operatorname{Jac}(C),\lambda_C\bigr].
\]
This is the \emph{Torelli map}. By Torelli's theorem
\cite[Thm.~12.1]{MilneJacobian1986}, an isomorphism between
canonically polarized Jacobians determines an isomorphism
between the corresponding curves, so $t$ is injective
on geometric points.

For a stable curve $C$ with $b_1(\Gamma_C)>0$, the
generalized Jacobian $J_C$ has a nontrivial torus part
and therefore does not define a point of $A_{g/k}$.
Its abelian quotient has dimension
\[
 \sum_i g_i=g-b_1(\Gamma_C).
\]
To study the limit of the Torelli map at such a curve,
we will use a projective target whose boundary records
principally polarized abelian varieties of these smaller
dimensions.

Write
\[
 \Z_{(p)}=\left\{\frac{a}{b}\in\Q
 \;\middle|\; a,b\in\Z,\ p\nmid b\right\}
\]
for the localization of $\Z$ at the prime ideal $(p)$.
This is a local ring with maximal ideal $p\Z_{(p)}$ and
residue field $\F_p$. Its $p$-adic completion is the ring
of $p$-adic integers $\Z_p$.
Let $S=\operatorname{Spec}\Z_{(p)}$, and let
$\mathcal A_{g/S}$ denote the moduli stack of principally
polarized abelian schemes of relative dimension $g$ over $S$.
Write $A_{g/S}$ for its coarse moduli space,
formed over $S$. By Faltings--Chai
\cite[Thm.~V.2.3]{FaltingsChaiDegeneration}, there is a
normal projective $S$-scheme $A^*_{g/S}$ containing
$A_{g/S}$ as a dense open subscheme.
We set
\[
 P=A^*_{g/S}\times_S\operatorname{Spec}k,
\]
which is a projective $k$-scheme.

The base change of the coarse moduli morphism over $S$ gives
\[
 \mathcal A_{g/k}\longrightarrow
 (A_{g/S})_k,
\]
where the subscript $k$ denotes base change from $S$ to
$\operatorname{Spec}k$.
Recall that every morphism from $\mathcal A_{g/k}$ to an
algebraic space factors uniquely through its coarse moduli
space $A_{g/k}$. Applying this universal property to the map
above gives the comparison morphism
\[
 \beta_g:A_{g/k}\longrightarrow(A_{g/S})_k.
\]
Thus we obtain the composite
\[
 t_P\colon M_{g/k}\xrightarrow{\ t\ }A_{g/k}
 \xrightarrow{\ \beta_g\ }(A_{g/S})_k
 \hookrightarrow P.
\]
We retain this distinction because the open subscheme of $P$
supplied by the integral compactification is $(A_{g/S})_k$,
whereas the Torelli map over $k$ naturally takes values in
$A_{g/k}$.

By Faltings--Chai \cite[Thm.~IV.5.7]{FaltingsChaiDegeneration},
the stack $\mathcal A_{g/S}$ admits a proper smooth toroidal
compactification equipped with a semiabelian scheme extending
the universal abelian scheme. \textcolor{black}{Fix such a compactification
$\mathcal T_S$ and put $\mathcal T=(\mathcal T_S)_k$.} Let
\[
 \mathcal G\longrightarrow\mathcal T
\]
be the corresponding semiabelian scheme, extending the
universal abelian scheme over $\mathcal A_{g/k}$.

Faltings--Chai \cite[Thm.~V.2.3(1)]{FaltingsChaiDegeneration}
constructs a canonical morphism from the chosen toroidal
compactification to $A^*_{g/S}$. Its base change to $k$ gives
\[
 q\colon\mathcal T\longrightarrow P.
\]
For $0\leq h\leq g$, let $\mathcal A_{h/S}$ and $A_{h/S}$
denote the moduli stack of principally polarized abelian
schemes of relative dimension $h$ over $S$ and its coarse
space. \textcolor{black}{By part~(3) of the same theorem, the locally
closed subschemes $(A_{h/S})_k$, for $0\leq h\leq g$, stratify $P$.
The stratum for $h=0$ is a point.}

Part~(4) describes $q$ on geometric points.
If the semiabelian fibre $\mathcal G_x$ has abelian quotient
$B_x$ of dimension $h$, then $q(x)$ is the point of
$(A_{h/S})_k$ corresponding to $B_x$ with its induced
principal polarization. This polarization is part of the
boundary data described in
\cite[Cor.~IV.5.11(b)(i)]{FaltingsChaiDegeneration}.
These objects will be used below to determine limits of
Jacobians and compare Hodge lines.

As above, there is a comparison morphism
\[
 \beta_h\colon A_{h/k}\longrightarrow
 (A_{h/S})_k,
\]
where $A_{h/k}$ is the coarse space formed over $k$.
We regard geometric points of each stratum as points of $P$.

The comparison morphism $\beta_h$ to the stratum is a finite
universal homeomorphism. In the commutative diagram
\[
\begin{tikzcd}
 \mathcal A_{h/k} \arrow[r] \arrow[d] & (A_{h/S})_k \\
 A_{h/k} \arrow[ur,"\beta_h"'] &
\end{tikzcd}
\]
the vertical arrow is the coarse moduli morphism formed
over $k$, and the horizontal arrow is the base change
of the coarse moduli morphism over $S$. Both are proper
universal homeomorphisms
\cite[Thm.~1.1 and its proof]{ConradKeelMori}.
After pulling this diagram back along any morphism
$U\to(A_{h/S})_k$, the horizontal and vertical arrows
still induce homeomorphisms. Hence the base change of
$\beta_h$ does as well.
Thus $\beta_h$ is a universal homeomorphism.
Since it is separated and of finite type, it is proper
with finite geometric fibres, hence finite
\cite[Tag~02LS]{Stacks}.
In particular, \textcolor{black}{$\beta_h$ induces a bijection between
the geometric points of $A_{h/k}$ and those of $(A_{h/S})_k$,
and it preserves dimensions.}

\subsection{The extended Torelli morphism}
\label{subsec:extension-by-graph}

\textcolor{black}{By a theorem of Mumford and Namikawa, the Torelli
map extends to the second Voronoi compactification \cite[\S18]{Namikawa1976II}. Alexeev and Brunyate
give a proof valid in arbitrary characteristic
\cite[\S4, Thm.~4.4, pp.~182--183]{AlexeevBrunyateTorelli}.}
We record the extension to the target $P$ fixed above, together
with its boundary values. At the boundary, the image is determined
by the Jacobians of the normalization components with their
canonical polarizations.

We recall the objects used to compare boundary limits. A
semiabelian scheme is a smooth separated group scheme whose
geometric fibres are extensions of abelian varieties by tori;
the toric rank may vary between fibres
\cite[Def.~I.2.3 and Rem.~I.2.4]{FaltingsChaiDegeneration}.
A family of stable semiabelic pairs consists of a proper flat
family $\mathcal X$, acted on by a semiabelian scheme
$\mathcal J$, together with a relative effective Cartier
divisor $\Theta$. Each geometric fibre of $\mathcal X$ is
connected, reduced, and seminormal, has the same dimension as
the acting group, and has finitely many orbits with stabilizers
contained as subtori in the toric part. The divisor $\Theta$
is fibrewise ample and contains no entire orbit
\cite[\S1.1]{AlexeevCompleteModuli}.

Let $R$ be a complete discrete valuation \textcolor{black}{$k$-algebra}
with fraction field $K$ and residue field $k$ \textcolor{black}{via its structure map}. For a stable
curve family $Y/R$ with smooth generic fibre and
special fibre $C$, write
\[
 \mathcal J=\operatorname{Pic}^0_{Y/R},
 \qquad
 \mathcal X=\overline{\operatorname{Jac}}^{\,g-1}(Y/R).
\]
Here $\mathcal X$ is the relative compactified Jacobian
parametrizing semistable rank-one torsion-free sheaves of
degree $g-1$, up to \emph{S-equivalence}. Two such sheaves on
the same fibre are S-equivalent if their associated graded
sheaves for Jordan--H\"older filtrations are isomorphic
\cite[\S1.1]{AlexeevCompactifiedJacobians}.
Tensor product gives an action
\[
 \mathcal J\times_R\mathcal X\longrightarrow\mathcal X.
\]
Together with its relative theta divisor, $\mathcal X$ forms a
family of stable semiabelic pairs
\cite[Thms.~5.1--5.3 and Cor.~5.4]{AlexeevCompactifiedJacobians}.
Its generic and special fibres are
\[
 \begin{aligned}
 (\mathcal X_K,\Theta_K)
 &=\bigl(\textcolor{black}{\operatorname{Pic}^{g-1}_{Y_K/K}},\Theta_K\bigr),\\
 (\mathcal X_k,\Theta_k)
 &\cong\bigl(\overline{\operatorname{Jac}}^{\,g-1}(C),\Theta_C\bigr).
 \end{aligned}
\]
In particular, the special theta pair is determined by $C$,
independently of the chosen family.

The Mumford--Faltings--Chai degeneration data of the polarized
generic Jacobian include an auxiliary abelian scheme $A/R$,
an ample line bundle $M$ on $A$, and torus, period, and gluing
data. One has
\[
 A_k\cong J_C/T_C.
\]
The scheme $A/R$ arises from the Raynaud extension of
$\mathcal J/R$; it is not generally obtained by taking the
abelian quotient of each fibre of $\mathcal J/R$
\cite[Ch.~II, \S1]{FaltingsChaiDegeneration}.
After a finite extension $R\subset R'$ of discrete valuation rings,
put $K'=\operatorname{Frac}(R')$. From the full degeneration data,
Alexeev constructs a family $(\mathcal X',\Theta')$ of stable
semiabelic pairs over $R'$ whose generic fibre satisfies
\[
 (\mathcal X'_{K'},\Theta'_{K'})
 \cong (\mathcal X_K,\Theta_K)\times_{\Spec K}\Spec K'.
\]
\textcolor{black}{The polarized abelian quotient of the semiabelian group
acting on $\mathcal X'_k$ is $(A_k,\lambda_{M_k})$, where
$\lambda_{M_k}\colon A_k\to A_k^\vee$ is the polarization induced by $M_k=M|_{A_k}$}
\cite[Thm.~5.6.1 and proof of Thm.~5.7.1, Case~2]{AlexeevCompleteModuli}.

For a stable curve $C/k$ of genus $g$, recall the normalization
components $C_i$ from
Section~\ref{subsec:generalized-jacobians-toric-rank}.
Put $h=\sum_i g(C_i)$. We define a map of sets
\[
 b\colon\overline M_{g/k}(k)\longrightarrow P(k),
 \qquad [C]\longmapsto
 \beta_h\!\left(
 \left[\prod_i\operatorname{Jac}(C_i),
       \prod_i\lambda_{C_i}\right]\right),
\]
where $\lambda_{C_i}$ is the canonical principal polarization.
The polarized product depends only on the isomorphism class
of $C$, so this map is well defined. We write $b(C)$ for $b([C])$.
A genus-zero component contributes the trivial abelian variety.
For any algebraically closed extension $\Omega/k$, the same
formula after base change defines
$b\colon\overline M_{g/k}(\Omega)\to P(\Omega)$.
We first compute the limit in $P$ along a stable family.

\begin{lemma}[Torelli limit of a stable family]
\label{lem:torelli-stable-limit}
Let $k$ be algebraically closed of characteristic $p>0$ and let $g\geq2$.
Let $R$ be a complete discrete valuation \textcolor{black}{$k$-algebra} with
fraction field $K$ and residue field $k$ \textcolor{black}{via its structure map}, and let $s$ be the closed
point of $\operatorname{Spec}R$. Let $Y\to\operatorname{Spec}R$ be
a stable family of curves of genus $g$ with smooth generic fiber $Y_K$
and special fiber $C=Y_k$. The Torelli point of $Y_K$ extends uniquely
to a morphism $\tau_R\colon\operatorname{Spec}R\to P$, and
$\tau_R(s)=b(C)$.
\end{lemma}

\begin{proof}
Properness and separatedness of $P$ give the existence and uniqueness
of $\tau_R$. It remains to compute its value at $s$.

The normalization sequence recalled in
Section~\ref{subsec:generalized-jacobians-toric-rank} gives
\[
 J_C/T_C\cong\prod_i\operatorname{Jac}(C_i).
\]
To determine the limit in $P$, we also need the polarization
on this quotient. In Alexeev's description, the compactified Jacobian
has a closed $J_C$-orbit that is a torsor under this quotient
\cite[Thm.~3.7, Cor.~5.4, and \S\S2.4, 5.5(3)(a)--(b)]{AlexeevCompactifiedJacobians}.
Caporaso and Viviani identify this orbit with
$\prod_i\operatorname{Pic}^{g_i-1}(C_i)$ and describe the restriction
of the theta divisor as the union of the pullbacks of the theta
divisors of the factors
\cite[Rem.~1.2.11, eq.~(1.7)]{CaporasoVivianiTorelli}.
This restriction therefore induces the product of the canonical principal
polarizations, so the polarized abelian quotient is
\[
 \left(\prod_i\operatorname{Jac}(C_i),\prod_i\lambda_{C_i}\right).
\]

Use the notation $\mathcal J$, $(\mathcal X,\Theta)$,
$(\mathcal X',\Theta')$, and $(A,M)$ recalled above for the
given family $Y/R$. Over $R'$, the base change of
$(\mathcal X,\Theta)$ and the pair $(\mathcal X',\Theta')$
extend isomorphic generic theta pairs over $K'$.
Uniqueness of stable semiabelic reduction identifies these families
after a further finite base change
\cite[Thm.~5.7.1]{AlexeevCompleteModuli}.
\textcolor{black}{The semiabelian groups acting on the two special fibres}
therefore have isomorphic polarized abelian quotients. For Alexeev's family $(\mathcal X',\Theta')$, the polarized abelian quotient of the group acting on $\mathcal X'_k$ is $(A_k,\lambda_{M_k})$. For the compactified Jacobian family $(\mathcal X,\Theta)$, the group acting on $\mathcal X_k$ is $J_C$, and its abelian quotient $J_C/T_C$ carries the product polarization computed above. Thus
\[
 (A_k,\lambda_{M_k})
 \cong\left(\prod_i\operatorname{Jac}(C_i),\prod_i\lambda_{C_i}\right).
\]
Thus the polarized abelian quotient arising from the given
family depends only on $C$.

Recall the proper toroidal stack $\mathcal T$ and the contraction
\[
 q\colon\mathcal T\longrightarrow P
\]
from Section~\ref{subsec:coarse-torelli}.
The map $\tau_R$ was obtained directly from properness of $P$.
To determine its special value, we use the toroidal stack
$\mathcal T$. The contraction $q$ sends a boundary point to
the moduli point of the polarized abelian quotient of its
semiabelian fibre. Identifying this fibre with $J_C$ will
therefore let us apply the polarization calculation above.
The required extension to $\mathcal T$ may involve a finite
base change.

By the valuative criterion for properness
\cite[Thm.~4.19]{DeligneMumford}, there is a finite extension
$K'/K$ such that the generic polarized Jacobian extends over
the integral closure $R'$ of $R$ in $K'$, giving a morphism
\[
 \widetilde\tau\colon\operatorname{Spec}R'\longrightarrow\mathcal T.
\]
Since $R$ is a complete discrete valuation ring, $R'$ is again
a complete discrete valuation ring and is finite over $R$
\cite[Tags~032W and~09EM]{Stacks}.
The pulled-back semiabelian family agrees with the relative
generalized Jacobian of $Y_{R'}/R'$.
Indeed, both extend the generic Jacobian, and the generic
isomorphism and its inverse extend uniquely
\cite[\S9.4, Thm.~1]{BLR},
\cite[Prop.~I.2.7]{FaltingsChaiDegeneration}.

Let $s'$ be the closed point of $\operatorname{Spec}R'$.
The description of $q$ in Section~\ref{subsec:coarse-torelli}
gives
\[
 q\bigl(\widetilde\tau(s')\bigr)
 =\beta_h\!\left(
 \left[\prod_i\operatorname{Jac}(C_i),\prod_i\lambda_{C_i}\right]\right)=b(C).
\]
Consider the diagram
\[
\begin{tikzcd}[column sep=large,row sep=large]
 \operatorname{Spec}R' \arrow[r,"\widetilde\tau"] \arrow[d]
   & \mathcal T \arrow[d,"q"] \\
 \operatorname{Spec}R \arrow[r,"\tau_R"'] & P,
\end{tikzcd}
\]
where the left vertical map is induced by $R\subset R'$.
The two composites from $\operatorname{Spec}R'$ to $P$ agree
on $\operatorname{Spec}K'$, where both give the Torelli point
of $Y_{K'}$. Since $P$ is separated and
$\operatorname{Spec}R'$ is integral, they agree everywhere,
so the diagram commutes. Since the residue field of $R'$ is
still $k$, we obtain
\begin{equation}\label{eq:torelli-stable-limit}
 \tau_R(s)=b(C),
\end{equation}
where $s$ is the closed point of $\operatorname{Spec}R$.
In particular, this value is independent of the chosen smoothing.
\end{proof}

To apply Lemma~\ref{lem:torelli-stable-limit} to a morphism
$\operatorname{Spec}R\to\overline M_{g/k}$ whose generic point maps
into $M_{g/k}$, we need a stable family inducing the given morphism.
Stable reduction supplies such a family after a finite extension of
the fraction field of $R$.

\begin{lemma}[Realizing a \textcolor{black}{map to the coarse space} by a stable family]
\label{lem:stable-family-coarse-map}
Let $k$ be algebraically closed of characteristic $p>0$ and let $g\geq2$.
Let $R$ be a discrete valuation ring over $k$ with fraction field $K$.
Suppose that \textcolor{black}{the $k$-morphism} $h\colon\operatorname{Spec}R\to\overline M_{g/k}$ sends
the generic point into $M_{g/k}$. There exist an extension of discrete
valuation rings $R\subset R_L$, with $L=\operatorname{Frac}(R_L)$
finite over $K$, \textcolor{black}{and a stable family
$Y^{\mathrm{st}}\to\operatorname{Spec}R_L$ with smooth generic fiber.
The classifying morphism of this family, followed by the coarse moduli morphism,
is the composite}
\[
 \operatorname{Spec}R_L\longrightarrow\operatorname{Spec}R
 \xrightarrow{h}\overline M_{g/k}.
\]
If $h$ sends the closed point to a $k$-point $[C]$, then
$Y^{\mathrm{st}}_{\Omega}\simeq C\times_k\Omega$, where $\Omega$ is
an algebraic closure of the residue field of $R_L$.
\end{lemma}

\begin{proof}
Over an algebraic closure of $K$, the image of the generic point
under $h$ is represented by a smooth projective curve of genus $g$.
Its projective equations involve only finitely many coefficients,
so the curve is defined over a finite extension $K'/K$. Extend the
valuation of $R$ to $K'$. By the stable reduction theorem of Deligne
and Mumford \cite[Cor.~2.7]{DeligneMumford}
(see also \cite[Tag~0E98]{Stacks}), after a further finite extension
$L/K'$ the curve extends to a stable family $Y^{\mathrm{st}}$ over
a discrete valuation ring $R_L$ of $L$ extending the chosen valuation
ring of $K'$.

By construction, \textcolor{black}{the induced map
$\operatorname{Spec}R_L\to\overline M_{g/k}$} agrees
with the displayed composite on $\operatorname{Spec}L$. Since
$\overline M_{g/k}$ is separated, the two maps agree on
$\operatorname{Spec}R_L$ \cite[Tag~01KZ]{Stacks}. If the closed point
maps to $[C]$, then $Y^{\mathrm{st}}_{\Omega}$ and $C\times_k\Omega$
have the same $\Omega$-point of the coarse moduli space and are
therefore isomorphic.
\end{proof}

We now use these two lemmas to determine the fibers of the closure
of the graph of $t_P$ over $\overline M_{g/k}$. This will show that
the pointwise formula $[C]\mapsto b(C)$ defines a morphism of schemes.

\begin{proposition}[Extended Torelli morphism]
\label{prop:extended-torelli}
Let $k$ be algebraically closed of characteristic $p>0$ and
let $g\geq2$. With the notation above, the morphism
$t_P\colon M_{g/k}\to P$ extends uniquely to a morphism
\[
 \overline t_P\colon\overline M_{g/k}\longrightarrow P.
\]
For every stable curve $C$ over an algebraically closed
extension $\Omega/k$, one has
\begin{equation}\label{eq:extended-torelli-boundary}
 \overline t_P([C])=b(C),
\end{equation}
where the maps and strata defining $b(C)$ are understood
after base change to $\Omega$.
\end{proposition}

\begin{proof}
Let
\[
 \Gamma\subset\overline M_{g/k}\times_k P
\]
be the reduced closure of the graph of $t_P$, and write
\[
 \pi\colon\Gamma\longrightarrow\overline M_{g/k},
 \qquad
 r\colon\Gamma\longrightarrow P
\]
for the projections. We seek the dashed morphism in
\[
\begin{tikzcd}[column sep=large,row sep=large]
 \Gamma \arrow[d,"\pi"'] \arrow[r,"r"] & P \\
 \overline M_{g/k} \arrow[ur,dashed,"\overline t_P"'] &
\end{tikzcd}
\]
by showing that $\pi$ is an isomorphism.
The scheme $\Gamma$ is integral. The morphism $\pi$ is proper
and surjective, and it is an isomorphism over $M_{g/k}$.
We will show that it is an isomorphism everywhere.

Fix a closed point $[C]\in\overline M_{g/k}$. Choose a closed point
$z\in\Gamma$ over $[C]$ and an affine neighbourhood $V\subset\Gamma$ of $z$.
Since $\pi^{-1}(M_{g/k})$ is dense in $\Gamma$,
$V\cap\pi^{-1}(M_{g/k})$ is a nonempty open subset of $\Gamma$.
As $\Gamma$ is of finite type over $k$, the Nullstellensatz
then gives a closed point
\[
 w\in V\cap\pi^{-1}(M_{g/k}).
\]
There is an integral curve $E\subset V$ containing both $z$
and $w$ \cite[Lemma~10.1]{JardenRazonRumely}.
Thus $E\cap\pi^{-1}(M_{g/k})$ is a nonempty open subset of $E$,
so it contains the generic point of $E$.

Let $\widetilde E\to E$ be the normalization, and choose a
closed point $\widetilde z$ above $z$.
Since $\widetilde E$ is a normal curve,
$\mathcal O_{\widetilde E,\widetilde z}$ is a discrete valuation ring.
Its residue field is $k$, because $k$ is algebraically closed.
Completing this local ring gives
\[
 R=\widehat{\mathcal O}_{\widetilde E,\widetilde z}
\]
and an induced morphism
\[
 \gamma\colon\operatorname{Spec}R\longrightarrow\Gamma.
\]
The closed point maps to $z$, while the generic point maps to
the generic point of $E$, hence lies over $M_{g/k}$.

Apply Lemma~\ref{lem:stable-family-coarse-map} to the composite
$\pi\circ\gamma\colon\operatorname{Spec}R\to\overline M_{g/k}$, and
let $Y^{\mathrm{st}}\to\operatorname{Spec}R_L$ be the resulting stable
family. Since $R$ is complete with algebraically closed residue field
$k$, the ring $R_L$ is again complete with residue field $k$
\cite[Tags~032W and~09EM]{Stacks}. The closed point maps to
$\pi(z)=[C]$, so $Y^{\mathrm{st}}_k\simeq C$.

Let $\rho\colon\operatorname{Spec}R_L\to\operatorname{Spec}R$ be the natural
morphism, and let $s_L$ be the closed point of $\operatorname{Spec}R_L$.
By properness of $P$, the Torelli point of $Y^{\mathrm{st}}_L$
extends uniquely to a morphism
\[
 \tau_{R_L}\colon\operatorname{Spec}R_L\longrightarrow P.
\]
By \eqref{eq:torelli-stable-limit} in Lemma~\ref{lem:torelli-stable-limit},
applied to
$Y^{\mathrm{st}}/R_L$ with special fiber isomorphic to $C$, we have
$\tau_{R_L}(s_L)=b(C)$.
We compare $\tau_{R_L}$ with the composite
\[
 u:=r\circ\gamma\circ\rho\colon\operatorname{Spec}R_L\longrightarrow P.
\]
On the generic point these maps agree: over $M_{g/k}$, the scheme $\Gamma$
is the graph of $t_P$, so $u$ also gives the Torelli point of the same
smooth curve over $L$. Since $P$ is separated and
$\operatorname{Spec}R_L$ is integral, the two maps agree everywhere.
As $\gamma(\rho(s_L))=z$, evaluating at the closed point gives
\[
 r(z)=u(s_L)=\tau_{R_L}(s_L)=b(C).
\]
Thus every closed point of the fiber $\Gamma_{[C]}$ has
the same second coordinate. Since $\Gamma\hookrightarrow\overline M_{g/k}\times_k P$
is a closed immersion, base change to the $k$-point $[C]$ gives
\[
 \Gamma_{[C]}\hookrightarrow\{[C]\}\times_k P\simeq P.
\]
This map is the restriction of $r$, so it is injective on points.
Every closed point of the fiber maps to $b(C)$, and the fiber
is nonempty because $\pi$ is surjective. Thus it has exactly
one closed point. Closed points
are dense in every finite-type $k$-scheme, so its support
is a singleton.

For any closed point $x\in\Gamma$, its image $\pi(x)$ is closed
because $\pi$ is proper. Applying the preceding argument with
$[C]=\pi(x)$ shows that $x$ is the only point of its fiber,
so it is isolated there.
Consequently, the quasi-finite locus of $\pi$ contains every
closed point of $\Gamma$. Its closed complement must be empty,
since any nonempty closed subset of a finite-type $k$-scheme
contains a closed point \cite[Tags~02J6 and~01TI]{Stacks}.
Hence $\pi$ is proper and quasi-finite, and therefore finite.
It is also birational, and $\overline M_{g/k}$ is normal,
so $\pi$ is an isomorphism
\cite[Tags~02LS and~0AB1]{Stacks}.

We obtain the required morphism
\[
 \overline t_P=r\circ\pi^{-1}
 \colon\overline M_{g/k}\longrightarrow P.
\]
Its uniqueness follows from density of $M_{g/k}$ and
separatedness of $P$.

Finally, the same argument applies after any algebraically
closed extension $\Omega/k$, since coarse spaces commute
with this flat base change \cite[Thm.~1.1]{ConradKeelMori}.
Uniqueness identifies the resulting extension with
$(\overline t_P)_\Omega$, giving the asserted formula
for every geometric stable curve.
\end{proof}

\subsection{Compatibility with the Hodge line}
\label{subsec:torelli-hodge-line}

Let
\[
 \epsilon\colon\overline{\cM}_{g/k}\longrightarrow\overline M_{g/k}
 \quad\text{and}\quad
 f\colon\mathcal C\longrightarrow\overline{\cM}_{g/k}
\]
be the coarse moduli morphism and the universal stable curve, respectively.
The Hodge line on the stack of stable curves is $\lambda=\det(f_*\omega_f)$,
where $\omega_f$ is the relative dualizing sheaf. For the divisor argument,
we will identify a power of $\lambda$ with the pullback of an ample line
bundle on $P$ along $\overline t_P\circ\epsilon$.

Write $q_S\colon\mathcal T_S\to A^*_{g/S}$ for the toroidal contraction
over $S=\operatorname{Spec}\mathbb Z_{(p)}$ from
Section~\ref{subsec:coarse-torelli}. Faltings--Chai define $\omega$ on
$\mathcal T_S$ as the determinant of the bundle of invariant differentials
of the universal semiabelian scheme
\cite[Ch.~V, \S1, p.~137]{FaltingsChaiDegeneration}. Their theorem
\cite[Thm.~V.2.3(1), p.~150]{FaltingsChaiDegeneration} gives a very ample
line bundle $\mathcal L$ on $A^*_{g/S}$ and an integer $N>0$ such that
$q_S^*\mathcal L\simeq\omega^{\otimes N}$.

After base change to $k$, the contraction becomes $q\colon\mathcal T\to P$,
and $\mathcal L$ pulls back to a very ample line bundle $H$ on $P$.
Let $e\colon\mathcal T\to\mathcal G$ be the identity section.
Since relative differentials commute with base change, $\omega$ pulls back
to $\det(e^*\Omega^1_{\mathcal G/\mathcal T})$. Here
$e^*\Omega^1_{\mathcal G/\mathcal T}$ is the cotangent bundle at the
identity, dual to $\operatorname{Lie}(\mathcal G/\mathcal T)$.
We fix the resulting isomorphism
\begin{equation}\label{eq:toroidal-hodge-identity}
 q^*H\simeq
 \bigl(\det(\operatorname{Lie}(\mathcal G/\mathcal T)^\vee)\bigr)^{\otimes N}.
\end{equation}

\begin{corollary}[The Hodge line]\label{cor:torelli-hodge-line}
Let $H$ and $N>0$ be the ample line bundle on $P$ and the
integer in the fixed toroidal Hodge identity. Then
\[
 \epsilon^*\overline t_P^*H\simeq\lambda^{\otimes N}.
\]
In particular, $\lambda^{\otimes N}$ descends to the line bundle
$\overline t_P^*H$ on $\overline M_{g/k}$.
\end{corollary}

\begin{proof}
Put $h=\overline t_P\circ\epsilon$, and let
\[
 j\colon\cM_{g/k}\longrightarrow\mathcal A_{g/k}\hookrightarrow\mathcal T
\]
be the morphism sending a smooth curve to its polarized Jacobian.
Let $i$ denote the inclusion of the smooth-curve locus. In the diagram
\[
\begin{tikzcd}[column sep=large,row sep=large]
 \cM_{g/k} \arrow[r,"j"] \arrow[d,hook,"i"']
   & \mathcal T \arrow[d,"q"] \\
 \overline{\cM}_{g/k} \arrow[r,"h"'] & P,
\end{tikzcd}
\]
we have $q\circ j=h\circ i$ because both composites give the Torelli map
on smooth families. Hence $j^*q^*H\simeq i^*h^*H=(h^*H)|_{\cM_{g/k}}$.

For a stable family $f_V\colon X\to V$ with generalized Jacobian
$J=\operatorname{Pic}^0_{X/V}$, the Lie algebra of $J/V$ is
$R^1(f_V)_*\mathcal O_X$ \cite[\S8.4, Thm.~1(a)]{BLR}. Relative duality gives
\[
 \operatorname{Lie}(J/V)^\vee
 \simeq (R^1(f_V)_*\mathcal O_X)^\vee
 \simeq (f_V)_*\omega_{X/V}
\]
\cite[Tag~0E6N]{Stacks}. These identifications commute with base change.
For a smooth family, $\omega_{X/V}=\Omega^1_{X/V}$, so along $j$ the
invariant differentials of $\mathcal G/\mathcal T$ identify with the
Hodge bundle. Its determinant is $i^*\lambda$. Pulling back
\eqref{eq:toroidal-hodge-identity} along $j$ therefore gives
\[
 \varphi^\circ\colon (h^*H)|_{\cM_{g/k}}
 \xrightarrow{\sim}(\lambda^{\otimes N})|_{\cM_{g/k}}.
\]
We show that this isomorphism extends across the boundary.

Choose an \'etale atlas of $\overline{\cM}_{g/k}$, and let $U$ be an
integral component of the atlas. Since $\overline{\cM}_{g/k}$ is smooth
over $k$, so is $U$. In particular, $U$ is normal. The inverse image
$U^\circ$ of $\cM_{g/k}$ is dense in $U$, since the atlas map is open
and the smooth-curve locus is dense. Thus $\varphi^\circ$ defines a
rational isomorphism on $U$. We check its order at the codimension-one
points of $U$ lying in $U\setminus U^\circ$.

Fix such a point $x$ and put $R=\mathcal O_{U,x}$, a discrete valuation
ring. The universal curve pulls back to a stable family over $R$ with
smooth generic fiber, since the generic point of $U$ lies in $U^\circ$.
Let $\widehat R$ denote the completion of $R$. By the valuative criterion
for properness of $\mathcal T$ \cite[Thm.~4.19]{DeligneMumford}, after a
finite extension $L/\operatorname{Frac}(\widehat R)$ the polarized
generic Jacobian extends to a morphism $\widetilde\tau\colon
\operatorname{Spec}R'\to\mathcal T$, where $R'$ is the integral closure
of $\widehat R$ in $L$, again a discrete valuation ring.
Write $u$ for the map to $\overline{\cM}_{g/k}$ induced by
$\operatorname{Spec}R'\to U$. Consider the diagram
\[
\begin{tikzcd}[column sep=large,row sep=large]
 \operatorname{Spec}R' \arrow[r,"\widetilde\tau"] \arrow[d,"u"']
   & \mathcal T \arrow[d,"q"] \\
 \overline{\cM}_{g/k} \arrow[r,"h"'] & P.
\end{tikzcd}
\]
The two composites agree on $\operatorname{Spec}L$, where both give the
Torelli point of the same smooth curve. Since $P$ is separated and
$\operatorname{Spec}R'$ is integral, the diagram commutes.

Let $f_{R'}\colon Y\to\operatorname{Spec}R'$ be the stable family classified
by $u$. Its generalized Jacobian $\operatorname{Pic}^0_{Y/R'}$ is
semiabelian \cite[\S9.4, Thm.~1]{BLR}. It and $\widetilde\tau^*\mathcal G$
both extend $\operatorname{Jac}(Y_L)$. The identification of their
generic fibers and its inverse extend uniquely over $R'$
\cite[Prop.~I.2.7]{FaltingsChaiDegeneration}. Pulling back
\eqref{eq:toroidal-hodge-identity} along $\widetilde\tau$ and applying
the description of invariant differentials above to $Y/R'$ now gives
\[
 (h\circ u)^*H
 \simeq \bigl(\det((f_{R'})_*\omega_{Y/R'})\bigr)^{\otimes N}
 \simeq u^*(\lambda^{\otimes N}).
\]
The last identification is base change for the Hodge bundle.
Over $L$ these are exactly the identifications defining $\varphi^\circ$,
so this isomorphism extends the prescribed generic one over $R'$.

Choose generators of the two line bundles over $R$. In these generators
the rational isomorphism is multiplication by some
$c\in\operatorname{Frac}(R)^\times$. Its extension over $R'$ is
invertible, so $\operatorname{ord}_{R'}(c)=0$. A uniformizer of $R$
remains a uniformizer of $\widehat R$. Thus, if $e$ is the ramification
index of $R'/\widehat R$, then
\[
 0=\operatorname{ord}_{R'}(c)=e\,\operatorname{ord}_R(c).
\]
Hence $c$ is a unit in $R$, and $\varphi^\circ$ extends at $x$.

Choose an affine open $\operatorname{Spec}A\subset U$ on which both line
bundles are trivial. In these trivializations the rational isomorphism is
multiplication by $c\in\operatorname{Frac}(A)^\times$. For every height-one
prime $\mathfrak p$ of $A$, we have $c\in A_{\mathfrak p}^\times$.
Indeed, $\varphi^\circ$ is already an isomorphism on $U^\circ$, and the
preceding calculation gives order zero at the codimension-one points of
$U$ lying in the boundary. Since $A$ is a Noetherian normal domain,
\cite[Tag~031T]{Stacks} gives
\[
 A=\bigcap_{\operatorname{ht}(\mathfrak p)=1}A_{\mathfrak p}
 \quad\text{inside }\operatorname{Frac}(A).
\]
Both $c$ and $c^{-1}$ belong to every ring in this intersection, so
$c\in A^\times$. The resulting isomorphism on $\operatorname{Spec}A$
is uniquely determined by its value at the generic point of $U$.
Thus these local extensions glue to an isomorphism $\varphi_U$ on $U$.

Let $U_1,U_2$ be two components of the atlas, possibly equal, and write
\[
 W=U_1\times_{\overline{\cM}_{g/k}}U_2,\qquad
 p_i\colon W\longrightarrow U_i
\]
for their overlap and its projections. The projections are \'etale, so
they are open and $W$ is reduced \cite[Tag~03PC]{Stacks}. The open subset
$W^\circ=p_1^{-1}(U_1^\circ)=p_2^{-1}(U_2^\circ)$ is dense in $W$,
since the image under $p_1$ of any nonempty open subset meets $U_1^\circ$.
Both maps $p_i^*\varphi_{U_i}$ go from the pullback of $h^*H$ to the
pullback of $\lambda^{\otimes N}$ on $W$, so their ratio is an invertible
regular function on $W$. It equals $1$ on $W^\circ$, where both maps
restrict to the pullback of the same $\varphi^\circ$. Since $W$ is
reduced, this function equals $1$ everywhere.

The isomorphisms $\varphi_U$ therefore descend to an isomorphism
$h^*H\simeq\lambda^{\otimes N}$ on $\overline{\cM}_{g/k}$
\cite[Tag~06WT]{Stacks}. The line bundle $\overline t_P^*H$ is defined
on the coarse space $\overline M_{g/k}$. Since $h=\overline t_P\circ
\epsilon$, its pullback by $\epsilon$ is $h^*H$. The isomorphism just
constructed identifies this pullback with $\lambda^{\otimes N}$.
\end{proof}

\subsection{A normal compactification of the coarse moduli space}
\label{subsec:normal-compactification}

\textcolor{black}{Over $\C$, Hain used the Satake compactification
and the rank-one Picard group in his proof of monodromy density for
divisors \cite[Prop.~4.1 and Lem.~7.2]{HainDensity}. We use the
corresponding compactification over $k$ to show that the closure of
every divisor meets the locus of toric rank one.}

Since $\overline M_{g/k}$ is projective over $k$ and $P$ is separated
over $k$, the morphism $\overline t_P\colon\overline M_{g/k}\to P$ is
proper. Its image is therefore closed. Let $W\subset P$ denote this
image with its reduced induced structure, and write
$\overline t_W\colon\overline M_{g/k}\to W$ for the induced morphism.
Consider the Stein factorization
\[
 \overline M_{g/k}\xrightarrow{\,c\,}Y_g
 \xrightarrow{\,\rho\,}W\hookrightarrow P.
\]
The morphism $c$ is proper and surjective with geometrically connected
fibers, and $\rho$ is finite. This construction identifies $Y_g$ with
the relative normalization of $W$ in $\overline M_{g/k}$
\cite[Tag~03H0]{Stacks}.

\begin{lemma}\label{lem:normal-compactification}
The scheme $Y_g$ is normal and projective, and $c$ identifies
$M_{g/k}$ with a dense open subscheme of $Y_g$. Moreover, the line
bundle $L=\rho^*(H|_W)$ is ample on $Y_g$, and
\[
 (c\circ\epsilon)^*L\simeq\lambda^{\otimes N}.
\]
\end{lemma}

\begin{proof}
Let $V\subset W$ be a nonempty affine open. In the Stein
factorization, $\rho^{-1}(V)$ is the spectrum of the ring of regular
functions on $\overline t_W^{-1}(V)$. Thus pullback gives an
isomorphism
\[
 c^*\colon\Gamma\bigl(\rho^{-1}(V),\mathcal O_{Y_g}\bigr)
 \xrightarrow{\sim}
 \Gamma\bigl(\overline t_W^{-1}(V),\mathcal O_{\overline M_{g/k}}\bigr).
\]
The ring on the right is an integrally closed domain because
$\overline t_W^{-1}(V)$ is a nonempty open subscheme of the normal
integral scheme $\overline M_{g/k}$
\cite[Tag~0358]{Stacks}. The opens $\rho^{-1}(V)$ cover $Y_g$, so
$Y_g$ is normal. It is projective because $\rho$ is finite and $W$ is
projective over $k$.

To identify the smooth-curve locus in $Y_g$, let
$Q\subset\overline M_{g/k}$ be the open locus at which $\overline t_W$
is quasi-finite \cite[Tag~01TI]{Stacks}. On $M_{g/k}$, the composite
of $\overline t_W$ with $W\hookrightarrow P$ is
\[
 M_{g/k}\xrightarrow{\,t\,} A_{g/k}
 \xrightarrow{\,\beta_g\,} (A_{g/S})_k\hookrightarrow P.
\]
Recall that the source of $\beta_g$ is the coarse space formed over $k$, and
its target is the base change to $k$ of the coarse space formed over
$S$. This comparison morphism is a universal homeomorphism and hence
induces a bijection on geometric points. By the Torelli theorem, a
smooth curve over an algebraically closed field of any characteristic
is determined by its principally polarized Jacobian
\cite[Thm.~12.1(a)]{MilneJacobian1986}. Thus each geometric fiber of
$\overline t_W|_{M_{g/k}}$ has at most one point. This restriction is
of finite type and hence quasi-finite.
Since $M_{g/k}$ is open in $\overline M_{g/k}$, it follows that
$M_{g/k}\subset Q$.

By Zariski's Main Theorem in the form of \cite[Tag~03GW]{Stacks},
$c(Q)$ is open in $Y_g$, and $c$ restricts to an isomorphism
$Q\xrightarrow{\sim}c(Q)$. Under this isomorphism, the open subscheme
$M_{g/k}\subset Q$ maps to an open subscheme of $Y_g$. Its image is
dense because $M_{g/k}$ is dense in $\overline M_{g/k}$ and $c$ is
surjective. In particular, $c$ is birational.

Finally, $H|_W$ is ample because $W$ is closed in $P$. Its pullback
$L$ along the finite morphism $\rho$ is therefore ample
\cite[Tag~0B5V]{Stacks}. By construction,
$c^*L\simeq\overline t_P^*H$. Pulling back along the coarse moduli
morphism $\epsilon$ and applying Corollary~\ref{cor:torelli-hodge-line}
gives
\[
 (c\circ\epsilon)^*L
 \simeq\epsilon^*\overline t_P^*H
 \simeq\lambda^{\otimes N}.
\]
\end{proof}

\subsection{Boundary dimensions}
\label{subsec:boundary-dimensions}

Assume $g\ge 3$, and identify $M_{g/k}$ with its image in $Y_g$ as in
Lemma~\ref{lem:normal-compactification}. Let
$\Delta_0\subset\overline M_{g/k}$ be the closure of the locus of
irreducible curves with exactly one node. For
$1\le i\le\lfloor g/2\rfloor$, let $\Delta_i$ be the closure of the
locus of curves consisting of two smooth components of genera $i$ and
$g-i$ meeting in a single node. These are the irreducible boundary
divisors of $\overline M_{g/k}$ \cite[\S~1.2]{MoonenOortTorelliLocus}.

\begin{lemma}\label{lem:boundary-dimensions}
The image $c(\Delta_0)\subset Y_g$, with its reduced induced
structure, has dimension $3g-6$. Moreover, every irreducible component
of $Y_g\setminus M_{g/k}$ has codimension at least two in $Y_g$.
\end{lemma}

\begin{proof}
Since the composite $Y_g\xrightarrow{\rho}W\hookrightarrow P$ is
finite, $c(\Delta_i)$ and $\overline t_P(\Delta_i)$ have the same
dimension for each $i$ \cite[Tag~0ECG]{Stacks}. The curves with
exactly one node form a dense open subset of $\Delta_i$, so their
images are dense in $\overline t_P(\Delta_i)$.

Let $C$ be an irreducible curve with one node. Its normalization
$\widetilde C$ is smooth of genus $g-1$. By
Proposition~\ref{prop:extended-torelli}, $\overline t_P([C])$ is the
image under $\beta_{g-1}$ of the canonically polarized Jacobian of
$\widetilde C$. Every smooth curve of genus $g-1$ occurs as such a
normalization, since we may identify any two distinct points on it.
By Torelli, the polarized Jacobian determines $\widetilde C$, and
$\beta_{g-1}$ is bijective on geometric points. It follows that
\[
 \dim c(\Delta_0)=\dim M_{g-1/k}=3(g-1)-3=3g-6.
\]

Now let $C$ consist of two smooth components of genera $i$ and $g-i$
meeting in a single node. By Proposition~\ref{prop:extended-torelli},
its image in $P$ is determined by the product of their Jacobians with
the product polarization. In particular, the image is independent of
the points at which the components are joined. An elliptic factor
varies in the one-dimensional space $A_{1/k}$, while smooth curves
of genus $h\ge 2$ vary in dimension $3h-3$. We therefore obtain
\[
 \dim c(\Delta_i)\le
 \begin{cases}
  1+\bigl(3(g-1)-3\bigr)=3g-5, & i=1,\\[2pt]
  (3i-3)+\bigl(3(g-i)-3\bigr)=3g-6, & i\ge 2.
 \end{cases}
\]

Since $c$ is surjective, every point of $Y_g\setminus M_{g/k}$ lies
in one of the images $c(\Delta_i)$. Each has dimension at most
$3g-5$. By Lemma~\ref{lem:normal-compactification}, $M_{g/k}$ is
dense and open in $Y_g$, so $\dim Y_g=3g-3$. The claimed
codimension bound follows.
\end{proof}

\subsection{Divisor classes}
\label{subsec:divisor-classes}

We continue to assume $g\ge 3$. By Lemma~\ref{lem:boundary-dimensions},
every codimension-one point of $Y_g$ lies in $M_{g/k}$. Thus every
prime divisor on $Y_g$ is the closure of a unique prime divisor on
$M_{g/k}$. The schemes $M_{g/k}$ and $Y_g$ have the same function field and the same
local rings at these points, so restriction induces an isomorphism
of Weil divisor class groups
\[
 \operatorname{Cl}(Y_g)\xrightarrow{\sim}\operatorname{Cl}(M_{g/k}).
\]

\textcolor{black}{We adapt Hain's Picard-rank-one argument for ampleness
\cite[Cor.~3.4 and Lem.~7.2]{HainDensity} to the compactification $Y_g$.
Since the closure of an effective divisor is initially only a Weil
divisor on $Y_g$, we explain why a positive multiple is Cartier and
why its class is a positive multiple of the ample Hodge class.}

\begin{proposition}\label{prop:ample-divisor-closure}
\textcolor{black}{Let $D=\sum_i n_iD_i$ be a nonzero effective Weil divisor
on $M_{g/k}$. Write $E_i$ for the reduced closure of $D_i$ in $Y_g$, and
put $D_Y=\sum_i n_iE_i$.} Then some positive multiple of $D_Y$
is an ample Cartier divisor.
\end{proposition}

\begin{proof}
The coarse space $\overline M_{g/k}$ is $\Q$-factorial, hence so is
its open subscheme $M_{g/k}$. \textcolor{black}{Thus the natural map
$\operatorname{Pic}(M_{g/k})_{\Q}\to\operatorname{Cl}(M_{g/k})_{\Q}$
is an isomorphism.}
By Moriwaki's result with $n=0$, this group is generated by the
rational Hodge class
\cite[Corollary~5.2 and its proof]{MoriwakiQPicard}.

Recall from Lemma~\ref{lem:normal-compactification} that $L$ is ample
on $Y_g$ and that $L|_{M_{g/k}}$ pulls back to $\lambda^{\otimes N}$
on $\mathcal M_{g/k}$. Thus the rational Hodge class on $M_{g/k}$ is
$\tfrac1N[L|_{M_{g/k}}]$. The restriction isomorphism above now gives
\[
 \operatorname{Cl}(Y_g)_{\Q}=\Q\cdot[L].
\]

Choose integers $m>0$ and $r$ such that $m[D_Y]-r[L]$ is torsion in
$\operatorname{Cl}(Y_g)$. Let $t$ be the order of this torsion class.
Then $(tm)[D_Y]=(tr)[L]$ in $\operatorname{Cl}(Y_g)$.
For the rest of the proof, write $m$ for $tm$ and $r$ for $tr$.
Thus $mD_Y$ is linearly equivalent to a Cartier divisor representing
$L^{\otimes r}$. In particular, $mD_Y$ is Cartier.

Put $d=\dim Y_g$ and choose a Cartier divisor $A$ with
$\mathcal O_{Y_g}(A)\cong L$. As usual,
$L^d$ denotes the top self-intersection number $(Y_g\cdot L^d)$.
Then $mD_Y\sim rA$. By linearity and
invariance of intersection numbers under linear equivalence,
\[
 m(D_Y\cdot L^{d-1})
 =(mD_Y\cdot L^{d-1})
 =(rA\cdot L^{d-1})
 =r(A\cdot L^{d-1})
 =rL^d.
\]
\textcolor{black}{A very ample tensor power of $L$ embeds $Y_g$ in
projective space over $k$. If $Z\subseteq Y_g$ is an integral closed
subvariety, its intersection with $\dim Z$ general hyperplanes is a
nonempty finite scheme. The degree of $Z$ in this embedding is
positive because it equals the degree of this nonempty zero-dimensional
scheme. Since a positive multiple of $L$ defines the embedding,
multilinearity gives $(Z\cdot L^{\dim Z})>0$}
\cite[Tag~0BEV]{Stacks}.

Taking $Z=Y_g$ gives \textcolor{black}{$L^d>0$.}
Taking $Z=E_i$ gives $(E_i\cdot L^{d-1})>0$ for each prime
component $E_i$, since $\dim E_i=d-1$. As $D_Y$ is nonzero
and effective, the sum has at least one term and every $n_i$ is positive.
Therefore
\[
 D_Y\cdot L^{d-1}=\sum_i n_i(E_i\cdot L^{d-1})>0.
\]
Since $m>0$ and $L^d>0$, the equality
$m(D_Y\cdot L^{d-1})=rL^d$ forces $r>0$.
It follows that $mD_Y$ is ample, as
$\mathcal O_{Y_g}(mD_Y)\cong L^{\otimes r}$.
\end{proof}

\subsection{Meeting the rank-one boundary}
\label{subsec:meeting-rank-one-boundary}

\textcolor{black}{Let $D$ be a nonzero effective Weil divisor on $M_{g/k}$.
By Proposition~\ref{prop:ample-divisor-closure}, some positive multiple
of its closure $D_Y$ is an ample Cartier divisor.} We use
this ampleness to show that \textcolor{black}{the support of $D_Y$} meets $c(\Delta_0)$ over the Satake
stratum of abelian dimension $g-1$.

Let $B=c(\Delta_0)$ with its reduced induced structure. Since $\Delta_0$
is irreducible and $c$ is proper, $B$ is an integral closed subvariety
of the projective scheme $Y_g$. By Lemma~\ref{lem:boundary-dimensions},
we have $\dim B=3g-6$. Recall the finite map $\rho:Y_g\to W\subseteq P$,
and put
\[
 B^\circ=\{b\in B\mid \rho(b)\in (A_{g-1/S})_k\}.
\]
The locus of irreducible curves with exactly one node has dense image
in $B$ under $c$. By Proposition~\ref{prop:extended-torelli}, its
image under $\overline t_P$ lies in the $(g-1)$-stratum. Thus $\rho(B)$
lies in the closure of this stratum. Since the stratum is open in its
closure, $B^\circ$ is open in $B$ and contains the dense image of the
one-node locus. If $C$ is a stable curve with $c([C])\in B^\circ$,
the abelian quotient of $J_C$ has dimension $g-1$. Since $\dim J_C=g$,
its toric rank is one. To show that $\operatorname{Supp}(D_Y)$ meets $B^\circ$, we use
ampleness on $B$ and a bound on the dimension of $B\setminus B^\circ$.

\begin{lemma}\label{lem:divisor-meets-rank-one}
The intersection $\operatorname{Supp}(D_Y)\cap B^\circ$ contains
a closed $k$-point.
\end{lemma}

\begin{proof}
We first bound the dimension of $B\setminus B^\circ$. The closure of
the $(g-1)$-stratum lies in the boundary of $P$, so the points of
$B\setminus B^\circ$ map to strata of abelian dimension at most $g-2$.
Let $C$ be a stable curve with $c([C])\in B\setminus B^\circ$, and
let $g_i$ denote the genera of the components of its normalization.
By Proposition~\ref{prop:extended-torelli}, the image of $C$ in $P$
is determined by the product of their Jacobians with the product
polarization. In particular, $\sum_i g_i\le g-2$, and the attaching
points make no contribution to the dimension of the image.

Fix the list of positive genera. An elliptic factor varies in dimension
one, and smooth curves of genus $g_i\ge2$ vary in dimension $3g_i-3$.
Thus each factor with $g_i>0$ contributes at most $3g_i-2$ parameters.
If at least one $g_i$ is positive, the locus of products for this list
has dimension at most
\[
 \sum_{g_i>0}(3g_i-2)
 \le 3\sum_i g_i-2
 \le 3(g-2)-2
 =3g-8.
\]
If all $g_i$ are zero, the image is the unique point of the
zero-dimensional stratum. There are only finitely many lists of
positive genera with sum at most $g-2$. Since $\rho$ is finite,
we obtain
\[
 \dim(B\setminus B^\circ)\le 3g-8=\dim B-2.
\]

If $B\subseteq\operatorname{Supp}(D_Y)$, the intersection contains
all of $B^\circ$. Otherwise, choose an ample Cartier multiple $mD_Y$
as in Proposition~\ref{prop:ample-divisor-closure}. The canonical
section of $\mathcal O_{Y_g}(mD_Y)$ restricts to a nonzero section
on $B$. The restricted section must vanish somewhere on $B$.
Otherwise, it would generate $\mathcal O_{Y_g}(mD_Y)|_B$ at every
point and hence trivialize this line bundle. The trivial line bundle has divisor class zero. Since $\dim B>0$,
this would give \textcolor{black}{$\bigl(B\cdot(mD_Y)^{\dim B}\bigr)=0$},
contradicting the positivity of the top self-intersection of an ample
divisor \cite[Tag~0BEV]{Stacks}.

Since $B$ is integral, the nonzero restricted section is locally
given by a nonzero regular function, hence by a nonzerodivisor.
Its zero locus is therefore a nonempty effective Cartier divisor
on $B$. Each irreducible component of this divisor has dimension
$\dim B-1$ \cite[Tags~01X0 and~0BCN]{Stacks}. Since
$\dim(B\setminus B^\circ)\le \dim B-2$, no such component can be
contained in $B\setminus B^\circ$. Thus the support of the divisor
meets $B^\circ$.

In either case, $\operatorname{Supp}(D_Y)\cap B^\circ$ is nonempty.
This intersection is locally closed in $Y_g$ and of finite type over
the algebraically closed field $k$, so it contains a closed $k$-point.
\end{proof}

When $D$ is prime, we can lift such a boundary point to the closure
of $D$ in $\overline M_{g/k}$ and realize the lift as a specialization
of the generic point of $D$.

\begin{lemma}[Specialization to the rank-one boundary]\label{lem:specialization-rank-one-boundary}
Let $D\subset M_{g/k}$ be a prime divisor and let $\overline D$ be its
reduced closure in $\overline M_{g/k}$. There exists a discrete valuation ring
$R$ with fraction field $k(D)$ such that the generic-point morphism
$\operatorname{Spec}k(D)\to\overline D$ extends to a morphism
$\operatorname{Spec}R\to\overline D$. The closed point of
$\operatorname{Spec}R$ maps to a $k$-point $[C_0]$ for which the
generalized Jacobian $J_{C_0}$ has toric rank one. We may take $R$ to
be the local ring at the generic point of a prime divisor on a normal
variety birational to $\overline D$.
\end{lemma}

\begin{proof}
Let $D_Y$ be the closure of $D$ in $Y_g$. Since $c$ is
proper and identifies $M_{g/k}$ with its image in $Y_g$, the image
$c(\overline D)$ is the closure of $D$ in $Y_g$. Thus
\[
c(\overline D)=\operatorname{Supp}(D_Y).
\]

By Lemma~\ref{lem:divisor-meets-rank-one}, we may choose a closed
$k$-point $y\in\operatorname{Supp}(D_Y)\cap B^\circ$. The fiber of
$\overline D\to Y_g$ over $y$ is nonempty and of finite type over $k$.
Since $k$ is algebraically closed, this fiber contains a $k$-point
$[C_0]$, represented by a stable curve $C_0$. Since $c([C_0])=y\in
B^\circ$, Proposition~\ref{prop:extended-torelli} shows that the
abelian quotient of $J_{C_0}$ has dimension $g-1$. Hence $J_{C_0}$
has toric rank one.

To construct the valuation, let $X$ be the normalization of $\overline D$. Since the normalization
is finite and surjective, we may choose a closed point $x\in X$ above
$[C_0]$. Let $\widetilde X$ be the normalization of the blowup of $X$ at $x$,
and denote by $\pi_x$ the composite
\[
\widetilde X \xrightarrow{\nu} \operatorname{Bl}_x X \longrightarrow X.
\]
The normalization map $\nu$ is finite because $\operatorname{Bl}_x X$
is a variety over $k$, and the blowup map is projective
\cite[Tags~0BXR and~02NS]{Stacks}. Both maps are therefore proper, so
their composite $\pi_x$ is proper. The blowup is an isomorphism over
the dense open subset $X\setminus\{x\}$. Since $X$ is normal, this
open subset is normal, so $\nu$ is also an isomorphism there
\cite[Tags~02OS and~0H7D]{Stacks}. Thus $\pi_x$ is birational.
Its image is closed by properness and contains $X\setminus\{x\}$,
so $\pi_x$ is surjective.

The fiber of the blowup $\operatorname{Bl}_x X\to X$ over $x$ is its
exceptional divisor $E$. The fiber $\pi_x^{-1}(x)$ of the composite is
therefore the scheme-theoretic inverse image $\nu^{-1}(E)$.
A local equation for $E$ remains a nonzerodivisor on $\widetilde X$,
since $\nu$ is dominant and both schemes are integral. Thus $\pi_x^{-1}(x)$
is an effective Cartier divisor \cite[Tags~02OS and~02OO]{Stacks}.
It is nonempty because $\pi_x$ is surjective, and each of its irreducible
components has codimension one in $\widetilde X$.
Choose the generic point $\xi$ of one such component. Since
$\widetilde X$ is normal and $\xi$ has codimension one,
$R=\mathcal O_{\widetilde X,\xi}$ is a discrete valuation ring
\cite[Tags~0BCN and~00PD]{Stacks}. By birationality,
\[
\operatorname{Frac}(R)=k(\widetilde X)=k(\overline D)=k(D).
\]
The induced morphism $\operatorname{Spec}R\to\overline D$ therefore
sends the generic point to the generic point of $D$. The closed point
maps to $\xi$ in $\widetilde X$, then to $x$ in $X$, and finally to
$[C_0]$.
\end{proof}

Over $\mathbb{C}$, for $g\geq5$, the elliptic-tail
construction of Marcucci, Naranjo, and Pirola gives a curve in the
closure of every prime divisor whose generalized Jacobian has toric
rank one: one lets exactly one of the two elliptic tails become nodal
\cite[Rem.~4.6 and Cor.~4.7]{MarcucciNaranjoPirola2016}.
In characteristic $p$, Achter and Pries show that the closure of every
component of the $p$-rank $g-1$ stratum contains an irreducible curve
with a single node \cite[Lem.~3.3(a)--(b)]{AchterPries}.
For an arbitrary prime divisor $D$ and $g\geq3$,
Lemmas~\ref{lem:divisor-meets-rank-one}
and~\ref{lem:specialization-rank-one-boundary} use ampleness and a
dimension bound to find a point of the toric-rank-one boundary and
realize it as a specialization along the valuation of $k(D)$ defined
by a prime divisor on a normal variety birational to $\overline D$.
We now apply Lemma~\ref{lem:stable-family-coarse-map} to obtain,
after a finite extension of $k(D)$, the stable family used in the
inertia calculation of Lemma~\ref{lem:divisor-rank-one-inertia}.

\begin{proposition}[Semistable reduction of toric rank one]
\label{prop:semistable-toric-rank-one}
Let $k$ be an algebraically closed field of characteristic $p>0$,
let $g\ge 3$, and let $D\subset M_{g/k}$ be a prime divisor.
There exists a discrete valuation ring $R$ of $k(D)$ whose valuation
is trivial on $k$. For some finite extension $L/k(D)$ and some
discrete valuation ring $R_L$ of $L$ extending $R$, there is a
stable family $Y^{\mathrm{st}}\to\operatorname{Spec}R_L$\textcolor{black}{.
Its generic fiber is smooth and induces the composite
$\operatorname{Spec}L\to\operatorname{Spec}k(D)\to M_{g/k}$
on the coarse moduli space.} The Jacobian of this generic fiber has semistable reduction
of toric rank one over $R_L$.
\end{proposition}

\begin{proof}
By Lemma~\ref{lem:specialization-rank-one-boundary}, the generic-point
morphism $\operatorname{Spec}k(D)\to D$ extends to a morphism
$\operatorname{Spec}R\to\overline D$ sending the closed point to a
$k$-point $[C_0]$ for which $J_{C_0}$ has toric rank one. Since $R$
contains $k$, every element of $k^\times$ is a unit in $R$.

Applying Lemma~\ref{lem:stable-family-coarse-map} to the composite with
$\overline D\hookrightarrow\overline M_{g/k}$ gives $R_L$ and the
stable family $Y^{\mathrm{st}}$ with the asserted generic fiber.
If $\Omega$ is an algebraic closure of the residue field of $R_L$,
the same lemma gives $Y^{\mathrm{st}}_\Omega\simeq C_0\times_k\Omega$.
The relative generalized Jacobian
$\operatorname{Pic}^0_{Y^{\mathrm{st}}/R_L}$ is a semiabelian scheme
extending $\operatorname{Jac}(Y^{\mathrm{st}}_L)$
\cite[\S9.4, Thm.~1]{BLR}. Its geometric special fiber is
$J_{C_0}\times_k\Omega$, which has toric rank one.
\end{proof}

\section{Proof of the monodromy theorem}
\label{sec:monodromy-theorem}

Let $k$ be an algebraically closed field of characteristic $p>0$,
let $g\geq3$, and fix a prime $\ell\neq p$. We first consider a prime
divisor $D\subset M_{g/k}$ in the unpointed coarse moduli space.
Fix $R\subset R_L$ and $Y^{\mathrm{st}}$ as in
Proposition~\ref{prop:semistable-toric-rank-one}, and put
\[
 C=Y^{\mathrm{st}}_L,\qquad J=\Jac(C).
\]
Thus $L$ is a finite extension of $k(D)$, and
\textcolor{black}{the morphism to the coarse moduli space induced by $C/L$ is the composite
$\Spec L\to\Spec k(D)\to M_{g/k}$.} The preceding two sections show
that $J_{\overline L}$ is simple and that $J$ has semistable
reduction of toric rank one over $R_L$. Together these facts
determine the geometric endomorphism ring of $J$. The same
reduction gives an element of inertia whose difference from the
identity has rank one on the Tate module. We use this element to
determine the geometric monodromy group and then apply the
reductions of Section~\ref{sec:moduli-monodromy} to pointed families.

\subsection{Endomorphisms and the action of inertia}
\label{subsec:divisor-endomorphisms-inertia}

The curve $C_{\overline L}$ represents a geometric generic point
of the reduced inverse image of $D$ \textcolor{black}{under the coarse moduli morphism $\cM_{g/k}\to M_{g/k}$}. This is again
an irreducible divisor. Indeed, the coarse moduli morphism is a
homeomorphism on underlying spaces
\cite[Thm.~1.1 and its proof]{ConradKeelMori}, so the generic
point of $D$ still has codimension one. Thus
Theorem~\ref{thm:divisor-geometric-simplicity} applies to
$J_{\overline L}$. Together with the reduction supplied by
Proposition~\ref{prop:semistable-toric-rank-one}, this gives the
following consequence.

\begin{corollary}
\label{cor:divisor-geometric-endomorphisms}
The geometric generic Jacobian of $D$ has endomorphism ring $\Z$.
\end{corollary}

\begin{proof}
This is Zarhin's criterion for an absolutely simple abelian variety
with semistable reduction of toric rank one
\cite[Cor.~3.2]{ZarhinBigMonodromy}. \textcolor{black}{We recall Zarhin's argument
\cite[proof of Thm.~3.1]{ZarhinBigMonodromy} to identify explicitly
the one-dimensional torus whose rational Tate module appears in
Lemma~\ref{lem:divisor-rank-one-inertia}.}
{\color{black}Choose a finite extension $L'/L$ over which all the
geometric endomorphisms of $J$ are defined. Such an extension exists
because $\End(J_{\overline L})$ is finitely generated as an abelian
group and each generator is defined over a finite extension of $L$.
Let $R'$ be the localization of the integral closure of $R_L$ in $L'$
at a maximal ideal above the maximal ideal of $R_L$. Then $R'$ is a
discrete valuation ring with fraction field $L'$ and $R'\cap L=R_L$.
Over $R'$, the Jacobian $J_{L'}$ still has semistable reduction of
toric rank one. Let $\mathcal N$ be its N\'eron model over $R'$, and
let $\bar s$ be a geometric point over the closed point of $\Spec R'$. Since $\mathcal N$ is smooth over $R'$, the N\'eron mapping property
\cite[\S1.2, Def.~1]{BLR} gives a unique endomorphism of $\mathcal N$
extending each $L'$-endomorphism $u$ of $J_{L'}$. On the geometric
special fiber, this extension fixes the identity and hence induces an
endomorphism $u_0$ of $\mathcal N^0_{\bar s}$.
Continue to write $Y^{\mathrm{st}}$ for the base change to $R'$
of the stable family in Proposition~\ref{prop:semistable-toric-rank-one}.
The total space $Y^{\mathrm{st}}$ is Cohen--Macaulay since the family
is a local complete intersection over the regular base $\Spec R'$.
The generic fiber is smooth and the special fiber is nodal, so the
family is smooth away from the finitely many nodes of the special
fiber. Thus $Y^{\mathrm{st}}$ is regular away from these points.
Each node has codimension two in $Y^{\mathrm{st}}$, so
$Y^{\mathrm{st}}$ is regular in codimension one and hence normal
by Serre's criterion.
By \cite[\S9.7, Cor.~2]{BLR}, the identity component
$\mathcal N^0_{\bar s}$ is consequently the generalized Jacobian
of the geometric special fiber. Thus
\[
 \mathcal N^0_{\bar s}\simeq
 \operatorname{Pic}^0\!\left(Y^{\mathrm{st}}_{\bar s}\right).
\] Since the reduction is semistable of toric
rank one, there is an exact sequence
\[
 0\longrightarrow T\xrightarrow{i}\mathcal N^0_{\bar s}
 \xrightarrow{q}A\longrightarrow 0,
\]
where $T\simeq\mathbb G_m$ and $A$ is an abelian variety.
Since every homomorphism from $T$ to $A$ is zero, the composite
\[
 T\xrightarrow{i}\mathcal N^0_{\bar s}
 \xrightarrow{u_0}\mathcal N^0_{\bar s}\xrightarrow{q}A
\]
is zero. Thus $u_0\circ i$ factors through $\ker(q)=i(T)$,
giving a unique endomorphism $u_T\colon T\to T$ such that
\[
 u_0\circ i=i\circ u_T.
\]}
Restriction therefore gives a unital homomorphism
\[
 \varphi\colon\End(J_{\overline L})\otimes_{\Z}\Q
 \longrightarrow \End(T)\otimes_{\Z}\Q=\Q,\qquad
 \textcolor{black}{u\otimes r\longmapsto u_T\otimes r},
\]
where the last equality follows from $T\simeq\mathbb G_m$.
Since $J_{\overline L}$ is simple, the source is a division
algebra. The homomorphism $\varphi$ is consequently
injective and identifies its source with $\Q$. \textcolor{black}{Let $u\in\End(J_{\overline L})$.
Since $T\simeq\mathbb G_m$, its restriction $u_T$ is given by
$t\mapsto t^n$ for a unique integer $n$. Multiplication by $n$ on
$J_{\overline L}$, denoted by $[n]$, has the same restriction to $T$.
The injectivity of $\varphi$ therefore gives
$u\otimes 1=[n]\otimes 1$. Since $\End(J_{\overline L})$ is
torsion-free as an abelian group, it follows that $u=[n]$.
Thus $\End(J_{\overline L})=\Z$.}
\end{proof}

We now consider the action of inertia on
\[
 H=V_\ell J=T_\ell J\otimes_{\Zl}\Ql.
\]
By the Kummer identification in
Subsection~\ref{subsec:symplectic-local-system}, this is the fiber
$H_{\Ql}$ of $\HH_{\Ql}$ at the geometric point represented by
$C_{\overline L}$. We use the alternating pairing
$\langle\, ,\,\rangle$ induced by the canonical principal
polarization of $J$. The choice of compatible $\ell$-power roots
of unity in $k$ identifies its values with $\Ql$.

\textcolor{black}{Since $L/k$ is finitely generated, $L$ is the
function field of an integral affine $k$-variety. As $k$ is perfect,
the smooth locus $U$ of this variety is a dense open subset
\cite[Tag~0B8X]{Stacks}. Thus $U$ contains the generic point and
has function field $L$.}
\textcolor{black}{After shrinking $U$, the coefficients of projective
equations defining $C$ are regular on $U$, giving a projective model
with smooth generic fiber $C$. By properness, the image of its nonsmooth
locus is closed in $U$ and misses the generic point, so removing this
image gives a smooth proper family.} Write $\rho_\ell$ for the action of the absolute
Galois group $G_L$ on $H$. This action factors through $\pi_1(U)$.
Moreover, $G_L\to\pi_1(U)$ is surjective, since every connected
finite \'etale cover of the normal integral scheme $U$ is integral
and hence has connected generic fiber. Thus the Zariski closure
of $\rho_\ell(G_L)$ is the group $G_{U,\ell}$ of
Subsection~\ref{subsec:symplectic-local-system}.
The following calculation with inertia is the one used in the
proof of \cite[Thm.~4.3]{ZarhinBigMonodromy}. Here it gives an
element of \textcolor{black}{$\Lie(G_{U,\ell})$}.

\begin{lemma}
\label{lem:divisor-rank-one-inertia}
The Lie algebra of $G_{U,\ell}$ contains an endomorphism $N$ of
$H$ of rank one satisfying $N^2=0$.
\end{lemma}

\begin{proof}
\textcolor{black}{The residue map embeds $k$ in $\kappa(R_L)$,
since the valuation of $R_L$ is trivial on $k$.}
Choose an extension of this valuation to $\overline L$,
and let $I\subset G_L$ be its inertia group. We may henselize
$R_L$ to compute the inertia action.
Let $\bar s$ be a geometric point above the closed
point of $\Spec R_L$, and let $\mathcal N$ be the N\'eron model of
$J$ over $R_L$. As in the proof of
Corollary~\ref{cor:divisor-geometric-endomorphisms}, write $T$ for
the torus in the identity component $\mathcal N^0_{\bar s}$ of its
geometric special fiber. By Grothendieck's fixed-part theorem \cite[Exp.~IX, (2.2.3.3) and Prop.~2.2.5]{SGASevenI}, specialization identifies $W=H^I$ with $V_\ell(\mathcal N^0_{\bar s})$. Under this
identification, the inclusion $T\hookrightarrow\mathcal N^0_{\bar s}$
gives a subspace $V_\ell T$ of $W$. For the pairing on $H$ induced by the principal polarization, Grothendieck's orthogonality theorem gives $V_\ell T=W\cap W^\perp$ \cite[Exp.~IX, Thm.~2.4 and (2.5.2)]{SGASevenI}. Since $J$ has semistable reduction, $W^\perp\subseteq W$ by \cite[Exp.~IX, Prop.~3.5(i)$\Rightarrow$(ii)]{SGASevenI}, so $V_\ell T=W^\perp$.
These results apply over an arbitrary residue field. Since $T$
has dimension one, $W^\perp$ is a line and $W$ has codimension
one in $H$.

The inertia group fixes $k\subset L$, which contains all
$\ell$-power roots of unity. In particular $I$ preserves the
pairing. {\color{black}For $\sigma\in I$, $x\in H$, and $w\in W$, we obtain by bilinearity
\[
\begin{aligned}
 \langle(\rho_\ell(\sigma)-1)x,w\rangle
 &=\langle\rho_\ell(\sigma)x,w\rangle-\langle x,w\rangle\\
 &=\langle x,\rho_\ell(\sigma)^{-1}w\rangle-\langle x,w\rangle\\
 &=\langle x,w\rangle-\langle x,w\rangle=0.
\end{aligned}
\]
In the second equality we apply $\rho_\ell(\sigma)^{-1}$ to both entries of the first pairing. The last line uses $\rho_\ell(\sigma)^{-1}w=w$, since $w\in W=H^I$.}
It follows that
\[
 (\rho_\ell(\sigma)-1)H\subset W^\perp\subset W
 \subset\ker(\rho_\ell(\sigma)-1).
\]
As $W\neq H$, some $\sigma\in I$ acts nontrivially. \textcolor{black}{For this $\sigma$, the endomorphism
$N=\rho_\ell(\sigma)-1$ is nonzero. Its image is contained in the line
$W^\perp$, so $N(H)=W^\perp$ and $N$ has rank one.
Since $W^\perp\subseteq W=H^I$, $\sigma$ fixes every vector in $N(H)$.
Thus $N$ vanishes on its image and $N^2=0$.}

\textcolor{black}{It remains to show that $N\in\Lie(G_{U,\ell})$.
Since this Lie algebra is the tangent space at the identity, it suffices to
exhibit a curve in $G_{U,\ell}$ through $\mathrm{id}_H$ with tangent vector $N$.
Set $f(t)=\mathrm{id}_H+tN$. Since $N^2=0$, we have
\[
 f(s)f(t)=\mathrm{id}_H+(s+t)N+stN^2=f(s+t).
\]
In particular, $f(t)^{-1}=f(-t)$. Since its entries are polynomial in $t$,
$f$ defines an algebraic homomorphism $\mathbb G_{a,\Ql}\to\operatorname{GL}(H)$.}

\textcolor{black}{To see that $f$ takes values in $G_{U,\ell}$, recall that
$f(1)=\mathrm{id}_H+N=\rho_\ell(\sigma)$ with $\sigma\in I\subset G_L$.
Thus for every $n\in\Z$,
\[
 f(n)=f(1)^n=\rho_\ell(\sigma^n)\in G_{U,\ell}(\Ql).
\]
As $G_{U,\ell}$ is Zariski closed, its inverse image under $f$ is a closed
subset of $\mathbb A^1_{\Ql}$ containing all integers. These are infinitely
many distinct points in $\Ql$, and a nonzero polynomial has only finitely
many roots. Hence $f$ factors through $G_{U,\ell}$. Since $f(0)=\mathrm{id}_H$,
its differential $df_0\colon T_0\mathbb G_{a,\Ql}\to\End(H)$ factors through
$\Lie(G_{U,\ell})$. The formula $f(t)=\mathrm{id}_H+tN$ then gives
$df_0(1)=N\in\Lie(G_{U,\ell})$.}
\end{proof}

\subsection{Absolute irreducibility over a finitely generated field}
\label{subsec:divisor-arithmetic-model}

Corollary~\ref{cor:divisor-geometric-endomorphisms} determines the
endomorphisms of $J$ over $\overline L$. To use the Tate
endomorphism theorem to study its Tate module, we need a field
finitely generated over $\F_p$. We descend the smooth family over
$U$ as in Section~\ref{sec:moduli-monodromy}. Thus there are a
subfield $k_0\subset k$ finitely generated over $\F_p$, a smooth
geometrically integral $k_0$-variety $U_0$, and a smooth proper
curve over $U_0$ whose base change to $k$ is the given family.
Put $F_0=k_0(U_0)$, and let $J_0$ denote the Jacobian of the
generic fiber of this descended family. Then
$J_0\times_{F_0}L\simeq J$.

Using compatible geometric generic points, we identify
$V_\ell J_0$ with $H$. {\color{black}Let $\rho_{U_0,\ell}$ denote the monodromy representation
of the descended family on $H$. We write $G_{\mathrm{arith}}\subseteq\GSp(H)$
for the Zariski closure of $\rho_{U_0,\ell}(\pi_1(U_0))$.
The multiplier character $\mu$ of $\GSp(H)$ fits into the commutative diagram
\[
\begin{tikzcd}[column sep=large,row sep=large]
 \pi_1(U_0) \arrow[r,"\rho_{U_0,\ell}"] \arrow[d]
   & \GSp(H)(\Ql) \arrow[d,"\mu"] \\
 G_{k_0} \arrow[r,"\chi_\ell"']
   & \Ql^\times
\end{tikzcd}
\]
Here the left vertical map is induced by $U_0\to\Spec k_0$, and
$\chi_\ell$ is the $\ell$-adic cyclotomic character of $G_{k_0}$.
To see why the diagram commutes, recall that the principal polarization of
$J_0$ gives a $\pi_1(U_0)$-equivariant pairing with values in $\Ql(1)$.
For $\gamma\in\pi_1(U_0)$, let $\bar\gamma$ be its image in
$G_{k_0}$. The action of $\gamma$ on $\Ql(1)$ is multiplication
by $\chi_\ell(\bar\gamma)$. Thus, for the scalar-valued form fixed
above and $x,y\in H$, equivariance gives
\[
\begin{aligned}
 \mu(\rho_{U_0,\ell}(\gamma))\langle x,y\rangle
 &=\langle\rho_{U_0,\ell}(\gamma)x,\rho_{U_0,\ell}(\gamma)y\rangle\\
 &=\chi_\ell(\bar\gamma)\langle x,y\rangle.
\end{aligned}
\]
Consequently $\mu(\rho_{U_0,\ell}(\gamma))=\chi_\ell(\bar\gamma)$.}

\begin{lemma}
\label{lem:divisor-arithmetic-irreducibility}
The geometric monodromy group $G_{U,\ell}$ is contained in $\Sp(H)$
and is normal in $G_{\mathrm{arith}}$. The identity component
$G_{\mathrm{arith}}^\circ$ acts absolutely irreducibly on $H$.
\end{lemma}

\begin{proof}
{\color{black}Choose an algebraic closure $\overline{k_0}$ of $k_0$ inside $k$,
and write $U_{0,\overline{k_0}}=U_0\times_{k_0}\overline{k_0}$.
Since $U_0$ is geometrically integral, we have the exact sequence of fundamental
groups \cite[Tag~0BTX]{Stacks}
\[
 1\longrightarrow\pi_1(U_{0,\overline{k_0}})
 \longrightarrow\pi_1(U_0)\longrightarrow G_{k_0}\longrightarrow 1.
\]
By Lemma~\ref{lem:extension-algebraically-closed-constants}, extension of
constants from $\overline{k_0}$ to $k$ induces a surjection
$\pi_1(U)\twoheadrightarrow\pi_1(U_{0,\overline{k_0}})$.
The representation on $H$ over $U$ is the pullback of the descended
representation. Thus $G_{U,\ell}$ is the Zariski closure of
$\rho_{U_0,\ell}\bigl(\pi_1(U_{0,\overline{k_0}})\bigr)$.

By exactness, the arithmetic monodromy image normalizes this geometric
image and hence also its closure $G_{U,\ell}$. The normalizer of
$G_{U,\ell}$ is Zariski closed, so it contains the closure
$G_{\mathrm{arith}}$ of the arithmetic image. This proves normality.
Moreover, the multiplier diagram preceding the lemma shows that every
element of $\pi_1(U_{0,\overline{k_0}})$ has multiplier $1$. Taking the
closure gives $G_{U,\ell}\subseteq\ker(\mu)=\Sp(H)$.

To prove absolute irreducibility for $G_{\mathrm{arith}}^\circ$, we
first realize it as the monodromy group over a finite extension of
$F_0$, following the finite-extension argument of Zarhin
\cite[\S4, before Rem.~4.1]{ZarhinBigMonodromy}. Since $U_0$ is normal
and integral with function field $F_0$, the generic point induces a
surjection $G_{F_0}\twoheadrightarrow\pi_1(U_0)$
\cite[Tag~0BQM]{Stacks}. Write $\rho$ for the composite
\[
 \rho\colon G_{F_0}\twoheadrightarrow\pi_1(U_0)
 \xrightarrow{\ \rho_{U_0,\ell}\ }G_{\mathrm{arith}}(\Ql).
\]
Its image is Zariski dense by the definition of $G_{\mathrm{arith}}$.
\textcolor{black}{Since $G_{\mathrm{arith}}$ is of finite type, it has only
finitely many irreducible components over $\overline{\Ql}$,
and hence only finitely many connected components. Each connected
component is thus Zariski open, so the component map on
$\Ql$-points is locally constant in the $\ell$-adic topology.}
The composite of $\rho$ with the map to the finite component group is
continuous, so its kernel is an open normal subgroup of $G_{F_0}$.
It is therefore $G_{F_1}$ for a finite Galois extension $F_1/F_0$.
By construction,
\[
 \rho(G_{F_1})=\rho(G_{F_0})\cap G_{\mathrm{arith}}^\circ(\Ql).
\]
Since $G_{\mathrm{arith}}^\circ$ is Zariski open and closed in
$G_{\mathrm{arith}}$, this intersection is Zariski dense in
$G_{\mathrm{arith}}^\circ$.

Put $J_1=J_0\times_{F_0}F_1$. Choose an $F_0$-embedding
$F_1\hookrightarrow\overline L$. Since $J_0\times_{F_0}L\simeq J$,
base change gives an injection
\[
 \End_{F_1}(J_1)\hookrightarrow\End_{\overline L}(J_{\overline L})=\Z,
\]
where the equality is Corollary~\ref{cor:divisor-geometric-endomorphisms}.
Every multiplication-by-$n$ endomorphism is already defined over $F_0$,
so the image contains every integer. Thus $\End_{F_1}(J_1)=\Z$.

Under our identification $V_\ell J_1=H$, the action of $G_{F_1}$ is
the restriction of $\rho$. The field $F_1$ is finitely generated over
$\F_p$, and $\ell\ne p$. The semisimplicity and Tate homomorphism
theorems of Zarhin for $p>2$ and Mori for $p=2$, in the form of
\cite[Thms.~1.2 and~1.4]{ZarhinFiniteCharacteristic}, therefore give
semisimplicity of $H$ as a $G_{F_1}$-module and an isomorphism
\[
 \End_{F_1}(J_1)\otimes_{\Z}\Ql
 \xrightarrow{\ \sim\ }\End_{G_{F_1}}(H).
\]
In particular, every $G_{F_1}$-equivariant endomorphism of $H$ is scalar.

Suppose that $W\subset H$ is a nonzero proper $G_{F_1}$-stable
subspace. Since $H$ is semisimple as a $G_{F_1}$-module, there is a
$G_{F_1}$-stable complement $W'$ such that $H=W\oplus W'$. Let
$P\colon H\to H$ be the projection given by $P(w+w')=w$ for $w\in W$
and $w'\in W'$. Both summands are stable under $G_{F_1}$, so for every
$\sigma\in G_{F_1}$ we have
\[
 P\bigl(\sigma(w+w')\bigr)
 =P(\sigma w+\sigma w')
 =\sigma w
 =\sigma P(w+w').
\]
Thus $P\in\operatorname{End}_{G_{F_1}}(H)$. Since $W$ and $W'$ are
both nonzero, $P$ acts as the identity on $W$ and as zero on $W'$.
This contradicts
$\operatorname{End}_{G_{F_1}}(H)=\Ql\,\mathrm{id}_H$.
Hence $H$ is irreducible over $\Ql$. Each element of $\rho(G_{F_1})$ is a symplectic similitude of $H$ and hence
a $\Ql$-linear endomorphism of $H$. We take the linear span of these
operators in $\End_{\Ql}(H)$. This span is an algebra, since
$\rho(\sigma)\rho(\tau)=\rho(\sigma\tau)$ for $\sigma,\tau\in G_{F_1}$.
By the finite-dimensional Jacobson density theorem, as applied in
\cite[\S1.5]{ZarhinFiniteCharacteristic}, we have
\[
 \operatorname{span}_{\Ql}\rho(G_{F_1})=\End_{\Ql}(H).
\]
Since $H$ is finite-dimensional, extension of scalars gives the natural
identification
\[
 \End_{\Ql}(H)\otimes_{\Ql}\overline{\Ql}
 \simeq \End_{\overline{\Ql}}(H\otimes_{\Ql}\overline{\Ql}).
\]
Thus every endomorphism of $H\otimes_{\Ql}\overline{\Ql}$ is a linear
combination over $\overline{\Ql}$ of the operators $\rho(\sigma)\otimes 1$
with $\sigma\in G_{F_1}$. Since
$\rho(G_{F_1})\subseteq G_{\mathrm{arith}}^\circ(\Ql)$, a
$G_{\mathrm{arith}}^\circ$-stable subspace is stable under each of these
operators and hence under every endomorphism. If the subspace contains
a nonzero vector $v$, it contains every vector in the ambient space,
since a linear endomorphism can send $v$ to any prescribed vector. This proves absolute
irreducibility.}
\end{proof}

\subsection{A rank-one endomorphism in the geometric Lie algebra}
\label{subsec:divisor-rank-one-criterion}

We now use the absolute irreducibility of the arithmetic action
together with the rank-one operator in the geometric Lie algebra.
Zarhin uses the rank-one classification of Guillemin, Quillen,
and Sternberg to determine the arithmetic monodromy Lie algebra
in this setting
\cite[Rem.~4.2 and the proof of Thm.~4.3]{ZarhinBigMonodromy}.
Let $V$ be a finite-dimensional symplectic vector space over a field of
characteristic zero, with symplectic form $\langle\, ,\,\rangle$.
Its Lie algebra $\mathfrak{sp}(V)=\Lie(\Sp(V))$ consists of the
endomorphisms $A$ of $V$ satisfying
$\langle Ax,y\rangle+\langle x,Ay\rangle=0$ for all $x,y\in V$.
\textcolor{black}{Deligne's criterion \cite[Lem.~5.11]{DeligneWeilI}
assumes that the Lie algebra generated by the rank-one operators acts
irreducibly on $V$. In our application, absolute irreducibility is
known for the connected arithmetic monodromy group, while the
rank-one operator lies in the geometric Lie algebra. We give the
proof to show how conjugation by the arithmetic group and Lie
brackets inside the geometric Lie algebra yield $\mathfrak{sp}(V)$.}

\begin{lemma}
\label{lem:divisor-symplectic-generation}
Let $V$ be a finite-dimensional symplectic vector space over a
field of characteristic zero. Let $G\subset\GSp(V)$ be a
connected algebraic subgroup acting absolutely irreducibly on
$V$. Suppose that $G$ normalizes a closed algebraic subgroup
$K\subset\Sp(V)$. If $\Lie K$ contains an endomorphism of rank
one, then $K=\Sp(V)$.
\end{lemma}

\begin{proof}
It suffices to prove the equality after extending scalars to
an algebraic closure of the ground field. Over this field $G$ is connected
and acts irreducibly on $V$.

Let $N\in\Lie K$ have rank one, and choose $v\ne0$ spanning its image.
Write $N(x)=f(x)v$ for a nonzero linear form $f$. Since
$N\in\mathfrak{sp}(V)$, we have
\[
 0=\langle Nx,y\rangle+\langle x,Ny\rangle
  =f(x)\langle v,y\rangle+f(y)\langle x,v\rangle.
\]
By nondegeneracy, we may choose $y$ with $\langle y,v\rangle=1$.
Then $f(x)=f(y)\langle x,v\rangle$ for every $x$, and $f(y)\ne0$
since $f\ne0$. Thus, writing
\[
 N_z(x)=\langle x,z\rangle z \qquad (z\in V),
\]
we obtain $N_v=f(y)^{-1}N\in\Lie K$.

We next use the action of $G$ to obtain such an operator for every vector
in the orbit $Gv$. If $b\in G$ has similitude multiplier $\mu(b)$, then
\[
 (bN_vb^{-1})(x)
 =\langle b^{-1}x,v\rangle\,bv
 =\mu(b)^{-1}\langle x,bv\rangle\,bv
 =\mu(b)^{-1}N_{bv}(x).
\]
Since $G$ normalizes $K$, conjugation by $b$ preserves $\Lie K$.
Hence $N_u\in\Lie K$ for every $u\in Gv$. By absolute irreducibility, the representation of $G$ on $V$ remains irreducible after the scalar extension above. Thus the linear span of $Gv$, being a nonzero $G$-stable subspace of $V$, equals $V$. Moreover, $Gv$ is irreducible as a variety: the connected algebraic group $G$ is irreducible, and $Gv$ is its image under the orbit morphism $b\mapsto bv$.

For $u,w\in V$, put
\[
 M_{u,w}(x)=\langle x,u\rangle w+\langle x,w\rangle u.
\]
With the convention $[N_u,N_w]=N_uN_w-N_wN_u$, we have
\begin{equation}\label{eq:divisor-rank-one-bracket}
 [N_u,N_w](x)=-\langle u,w\rangle M_{u,w}(x).
\end{equation}
Fix $u\in Gv$. Since $u\ne0$, the pairing is nondegenerate, and $Gv$
spans $V$, there is some $w\in Gv$ with $\langle u,w\rangle\ne0$.
Indeed, otherwise the linear form $w\mapsto\langle u,w\rangle$ would vanish on $\operatorname{span}(Gv)=V$, contradicting nondegeneracy and $u\ne0$. Such $w$ form a nonempty open subset of the irreducible orbit $Gv$
and hence are dense in $Gv$. Their linear span is Zariski closed and
therefore contains $Gv$, so these vectors span $V$. For each such $w$,
the bracket formula~\eqref{eq:divisor-rank-one-bracket} gives $M_{u,w}\in\Lie K$. By linearity in $w$,
this holds for every $w\in V$. Now let $u$ vary in $Gv$. Since $Gv$
spans $V$, linearity in $u$ gives $M_{u,w}\in\Lie K$ for all $u,w\in V$.

It remains to see that these operators span $\mathfrak{sp}(V)$.
Nondegeneracy of the pairing identifies an endomorphism $A$ with the
bilinear form $(x,y)\mapsto\langle Ax,y\rangle$. The condition
$A\in\mathfrak{sp}(V)$ says precisely that this form is symmetric:
\[
\begin{aligned}
 A\in\mathfrak{sp}(V)
 &\iff \langle Ax,y\rangle+\langle x,Ay\rangle=0
       \quad\text{for all }x,y\in V\\
 &\iff \langle Ax,y\rangle=\langle Ay,x\rangle
       \quad\text{for all }x,y\in V\\
 &\iff \text{the bilinear form }(x,y)\mapsto\langle Ax,y\rangle\text{ is symmetric.}
\end{aligned}
\]
For $A=M_{u,w}$, the corresponding form is
\begin{equation}\label{eq:divisor-symmetrized-products}
 \langle M_{u,w}(x),y\rangle
 =-\langle x,u\rangle\langle y,w\rangle
  -\langle x,w\rangle\langle y,u\rangle
 =\langle M_{u,w}(y),x\rangle.
\end{equation}
As $u,w$ vary, the covectors $\langle\,\cdot\,,u\rangle$ and
$\langle\,\cdot\,,w\rangle$ range over $V^\vee$. Write these covectors as $\alpha$ and $\beta$, respectively.
By~\eqref{eq:divisor-symmetrized-products}, the form corresponding to
$-M_{u,w}$ is $(x,y)\mapsto\alpha(x)\beta(y)+\beta(x)\alpha(y)$.
Choose a basis of $V$ with dual basis $e_1,\ldots,e_n$, and let
$E_{ij}$ denote the matrix with a single $1$ in position $(i,j)$
and zeros elsewhere. Taking $\alpha=e_i$ and $\beta=e_j$ gives
the matrix $E_{ij}+E_{ji}$ for $i<j$, while taking
$\alpha=\beta=e_i$ gives $2E_{ii}$. Since $2$ is invertible,
the span also contains each $E_{ii}$. Every symmetric matrix
$S=(s_{ij})$ has the expansion
\[
S=\sum_i s_{ii}E_{ii}+\sum_{i<j}s_{ij}(E_{ij}+E_{ji}).
\]
Thus the forms corresponding to the $M_{u,w}$ span all symmetric
bilinear forms on $V$. Under the identification above, the $M_{u,w}$
therefore span $\mathfrak{sp}(V)$. Since they all lie in
$\Lie K\subseteq\mathfrak{sp}(V)$, we obtain $\Lie K=\mathfrak{sp}(V)$. Algebraic groups are smooth in
characteristic zero, so $K$ and $\Sp(V)$ have the same dimension.
Since $\Sp(V)$ is irreducible, its closed subgroup $K$ must be all
of $\Sp(V)$.
\end{proof}

\subsection{Monodromy of families dominating a divisor}
\label{subsec:divisor-density-assembly}

\textcolor{black}{By Lemma~\ref{lem:divisor-arithmetic-irreducibility},
the identity component $G_{\mathrm{arith}}^\circ$ acts absolutely
irreducibly on the rational Tate module $H=V_\ell J$ and normalizes
the geometric monodromy group $G_{U,\ell}\subseteq\Sp(H)$.
Lemma~\ref{lem:divisor-rank-one-inertia} gives an endomorphism of rank
one in $\Lie G_{U,\ell}$. Applying
Lemma~\ref{lem:divisor-symplectic-generation} therefore gives}
\[
 G_{U,\ell}=\Sp(H).
\]
Thus the curve $C/L$\textcolor{black}{, whose image in $M_{g/k}$
is the generic point of $D$,} has full symplectic monodromy. We now pass to
the universal families on the moduli stack and its level covers.

\Needspace{10\baselineskip}
\begin{theorem}
\label{thm:divisor-monodromy}
Let $k$ be an algebraically closed field of characteristic $p>0$,
let $g\geq3$ and $n\geq0$, and let $\ell\neq p$ be a prime.
Let $\M$ be $\cM_{g,n/k}$ or $\cM_{g,n/k}[m]$, where $m\geq3$
is prime to $p$. Let $D\subset\M$ be a divisor and $E\to D$
a dominant morphism from a smooth integral $k$-variety.
Then the family over $E$ has full symplectic monodromy at
$\ell$. For a geometric basepoint $\bar e$ of $E$, put
$H_{\Ql}=(\HH_{\Ql})_{\bar e}$. As in
Subsection~\ref{subsec:symplectic-local-system}, write $\rho_{E,\ell}$
for the monodromy representation obtained by pulling back the
symplectic local system along $E\to D\hookrightarrow\M$. Then
\[
 G_{E,\ell}
 =\overline{\rho_{E,\ell}\bigl(\pi_1(E,\bar e)\bigr)}^{\,\mathrm{Zar}}
 =\Sp(H_{\Ql}).
\]
The same conclusion holds for a divisor in any connected finite
\'etale cover of $\cM_{g,n/k}[m]$, with the pulled-back family.
\end{theorem}

\begin{proof}

By Proposition~\ref{prop:reduction-unpointed}, it suffices first to
prove full symplectic monodromy over the normalization of every
divisor in the unpointed level scheme. Let $D_m\subset\cM_{g/k}[m]$
be such a divisor, and let $D_m^\nu$ be its normalization. The morphism
\[
 \cM_{g/k}[m]\longrightarrow M_{g/k}
\]
is finite and surjective, since the coarse moduli space is the
quotient of the level scheme by the finite group of changes of
symplectic basis. Its restriction to $D_m$ is finite onto its image.
Thus the reduced image $D\subset M_{g/k}$ is a prime divisor,
and $F:=k(D_m)$ is a finite extension of $k(D)$.

Let $C_m/F$ be the generic fiber of the universal curve over $D_m$,
and put $J_m=\operatorname{Jac}(C_m)$. We compare $C_m$ with the
curve $C/L$ constructed for $D$ in
Proposition~\ref{prop:semistable-toric-rank-one}. Recall that
$L/k(D)$ is finite and that the coarse moduli point of $C/L$ is
induced by the generic point of $D$. Choose $k(D)$-embeddings of
$F$ and $L$ into an algebraic closure of $k(D)$. Over this algebraic
closure, $C_m$ and $C$ represent the same geometric point of
$M_{g/k}$, so they are isomorphic. Since the curves are of finite
presentation, the isomorphism and its inverse descend to a finite
extension $L'$ of the compositum $LF$ \cite[Tag~01ZM]{Stacks}.
The resulting isomorphism of curves over $L'$ induces a
$G_{L'}$-equivariant isomorphism
\[
 H=V_\ell J\xrightarrow{\sim}V_\ell J_m
\]
preserving the pairings from the canonical principal polarizations.

We have proved above that the image of $G_L$ on $H$ is Zariski
dense in $\Sp(H)$. Since $L'/L$ is finite, the image of $G_{L'}$
has finite index in the image of $G_L$. Indeed, the separable part
of the extension gives a subgroup of finite index, and the purely
inseparable part induces an isomorphism on absolute Galois groups.
Write $\rho_\ell(G_L)$ as a finite union of cosets of
$\rho_\ell(G_{L'})$. Taking Zariski closures expresses $\Sp(H)$ as
a finite union of translates of the closure of $\rho_\ell(G_{L'})$.
If this closure were proper, each translate would be a proper closed
subset of $\Sp(H)$, contradicting its irreducibility. Thus
$\rho_\ell(G_{L'})$ is Zariski dense in $\Sp(H)$.

Under the displayed isomorphism, the image of $G_{L'}$ on $V_\ell J_m$
is thus Zariski dense in $\Sp(V_\ell J_m)$. The image of $G_F$
contains this subgroup and preserves the pairing, so it has the
same Zariski closure. Now $D_m^\nu$ is normal and integral with
function field $F$. Consequently, the generic-point map
$G_F\twoheadrightarrow\pi_1(D_m^\nu)$ is surjective
\cite[Tag~0BQM]{Stacks}. The action of $G_F$ on $V_\ell J_m$ is
the pullback of the monodromy representation on $D_m^\nu$ along
this map, so the two representations have the same image. We obtain
\[
 G_{D_m^\nu,\ell}=\Sp(V_\ell J_m).
\]

This verifies the hypothesis of Proposition~\ref{prop:reduction-unpointed}
for every unpointed level divisor. The proposition gives the assertion
for every smooth integral variety dominating a divisor in
$\cM_{g,n/k}[m]$, for all $n\geq 0$.
Corollary~\ref{cor:reduction-fine-level} then gives the assertion for
divisors in $\cM_{g,n/k}$.

Finally, let $X\to\cM_{g,n/k}[m]$ be a connected finite \'etale cover,
and let $D_X\subset X$ be a divisor. Its reduced image is a divisor
in the level scheme, since the finite morphism preserves the dimension
of $D_X$ and the two ambient schemes have the same dimension. If
$E\to D_X$ is dominant, then $E$ also dominates this image. The
family on $E$ is the pullback of the universal curve along the
composite $E\to X\to\cM_{g,n/k}[m]$, so the level-scheme case applies.

\end{proof}

\bibliographystyle{abbrvnat}
\bibliography{references}

\end{document}